\documentclass[10pt]{amsart}
\usepackage{amssymb}
\usepackage{amsfonts}
\usepackage{amscd}
\usepackage[dvips]{graphicx}
\usepackage[all,cmtip]{xy}
\usepackage{hyperref}
\usepackage{yhmath}
\usepackage{ulem}
\usepackage[all]{xy}
\usepackage{mathdots}
\usepackage{leftidx}
\usepackage{mathrsfs}
\usepackage{amsmath,amssymb, amsthm}
\usepackage{color}
\usepackage{xypic,color}
\usepackage{tikz}
\usepackage{mathdots}
\usepackage{picture}
\usepackage{enumerate}
\usepackage{makecell}
\usepackage{extarrows}
\usepackage{graphicx}
\usepackage{caption}
\usepackage{subcaption}
\usepackage{float}
\numberwithin{equation}{section}
\newtheorem{Them}{Theorem}[section]
\newtheorem{Lem}[Them]{Lemma}
\newtheorem{Def}[Them]{Definition}
\newtheorem{Cor}[Them]{Corollary}
\newtheorem{Prop}[Them]{Proposition}
\newtheorem{Ex}[Them]{Example}
\newtheorem{Rem}[Them]{Remark}
\newtheorem{Conj}[Them]{Conjecture}

\newcommand{\old}[1]{{\color{red}#1}}

\newcommand{\new}[1]{{\color{blue}#1}}

\newcommand{\add}{{\mathsf{add}}}

\newcommand{\stend}{\mathsf{\underline{\mathrm{End}}}}
\newcommand{\ind}{{\mathsf{ind}}}

\newcommand{\rad}{{\mathsf{rad}}}
\newcommand{\soc}{{\mathsf{soc}}}

\newcommand{\D}{{\mathsf{D}}}
\newcommand{\RS}{{\mathsf{Rsupp}}}
\newcommand{\LS}{{\mathsf{Lsupp}}}

\newcommand{\m}{\mathsf{mod}}
\newcommand{\stmod}{ \mathsf{\underline{mod}}}

\newcommand{\Hom}{{\mathsf {Hom}}}
\newcommand{\Coker}{{\mathsf {Coker}}}
\newcommand{\Ker}{{\mathsf {Ker}}}

\newcommand{\StHom}{\mathsf{\underline{Hom}}}

\newcommand{\dk}{\mathsf{dim_{_{k}}}}

\title[A construction of simple-minded systems over domestic Brauer graph algebras]{A construction of simple-minded systems over domestic Brauer graph algebras I: the 2-domestic case}
\author{Zhen Zhang}

\address{Zhen Zhang
	\newline Faculty of Arts and Sciences 
	\newline Beijing Normal University 
	\newline  Zhuhai 519087
	\newline P.R.China}
\email{zhangzhen@bnu.edu.cn}

\date{version of \today}
\newenvironment{Proof}[1][Proof]{\begin{trivlist}
\item[\hskip \labelsep {\bfseries #1}]}{\flushright$\Box$\end{trivlist}}

\newcommand{\sdp}{\times\kern-.2em\vrule height1.1ex depth-.05ex}

 \normalbaselines
\thanks{The research work is supported partially by NSFC (No.12031014 and No.12301044).}
\begin{document}
\renewcommand{\thefootnote}{\alph{footnote}}
\setcounter{footnote}{-1} 
\renewcommand{\thefootnote}{\alph{footnote}}
\setcounter{footnote}{-1}
\setcounter{footnote}{-1}

\begin{abstract}
Let $A$  be a 2-domestic Brauer graph algebra.  We present a construction for a family  of objects on $A$-$\stmod$ to be a simple-minded system and our construction provides  all simple-minded systems on $A$-$\stmod$.  As a byproduct, we provide a new proof of AR-conjecture for 2-domestic Brauer graph algebras and we prove that a weakly simple-minded system with a finite cardinality is a simple-minded system on $A$-$\stmod$. We also prove that an orthogonal system $\mathcal{S}$ which contains at least one object for an Euclidean component  extends to a simple-minded system on $A$-$\stmod$ and its extension closure $\mathcal{F}(\mathcal{S})$ is functorially finite. 
\end{abstract}

\maketitle

\section{Introduction}
Koenig and Liu \cite{KL}  introduced simple-minded systems in the stable module category of any artin algebra, which is a family of objects  satisfying orthogonality and generating condition. The authors  \cite{GLYZ} presented a sufficient and necessary condition  for an orthogonal system to be a simple-minded system over a  representation-finite self-injective algebra. They also \cite{GLYZ2} gave an explicit construction of sms's over   self-injective Nakayama algebras. For a  self-injective algebra, a few of necessary conditions of simple-minded systems were studied, including the finite cardinality of a simple-minded system,  invariance of a simple-minded system under stable equivalences 	(refer to \cite{KL}), and Chan, Liu  and Zhang \cite{CLZ} proved that any object on a quasi-tube of rank one is not in a simple-minded system over a self-injective algebra.  However, there is no full characterization of simple-minded systems for a class of representation-infinite self-injective algebras. In this paper, we shall present a construction for a family  of objects on $A$-$\stmod$ to be a simple-minded system in the stable module category of a 2-domestic Brauer graph algebra $A$ and this construction provides all simple-minded systems on $A$-$\stmod$. 
 
 Originated from the modular representation theory of finite groups,  Donovan and  Freislich \cite{DF} defined Brauer graph algebras.  Roggenkamp  \cite{Ro} showed that the class of Brauer graph algebras coincides with the class of symmetric special biserial algebras (see also \cite{MS}).   As a class of tame algebras, self-injective special biserial algebras are known pretty well, including the classification of indecomposable modules (string modules and band modules), irreducible maps between indecomposable modules, and the graph structure of stable Auslander-Reiten quivers (AR-quiver for short), etc. Bocian and Skowronski \cite{BoS} showed that a domestic Brauer graph algebra is 1-domestic or 2-domestic, and there is no $n$-domestic Brauer graph algebras for $n\geq3.$ Moreover, a Brauer graph algebra $A$ is 2-domestic if and only if the Brauer graph $G$ of $A$ is a graph with a unique cycle of even length and the multiplicity of each vertices of $G$ is one. Recently, there have been new developments for the study of Brauer graph algebras.  Adachi,  Aihara and  Chan  \cite{AAC} constructed   two-term tilting complexes over Brauer graph algebras, Opper and Zvonareva \cite{OZ} presented  a  derived equivalence classification of Brauer graph algebras and Antipov and Zvonareva \cite{AZ1} showed that Brauer graph algebras are closed under derived equivalence.  
 
 Auslander-Reiten conjecture  (AR-conjecture for short) states that stable equivalences preserve the number of isomorphism classes of non-projective simple modules.  Pogorzaly \cite{P} proved that, if AR-conjecture is true for a self-injective algebra, then it is true  for all  algebras.  Antipov and  A. Zvonareva, \cite{AZ}, Pogorzaly \cite{P1} proved that AR-conjecture holds for self-injective special biserial algebras. The study of simple-minded systems is closely related to AR-conjecture.
 If the cardinality of a simple-minded system over a finite dimensional algebra is a fixed number (the number of non-isomorphic, non-projective simple modules), then AR-conjecture holds for all finite dimensional algebras.  In this paper,  we give a new proof that AR-conjecture holds for 2-domestic Brauer graph algebras.

 Let $R$ be a self-injective algebra and $\mathcal{M}$ a family of objects on $R$-$\stmod$. In general, it is not easy to determine when the extension closure $\mathcal{F}(\mathcal{M})$ of  $\mathcal{M}$ is functorially finite on $R$-$\stmod$, even when $\mathcal{M}$ consists of  only finitely many objects. Let  $A$ be a  2-domestic Brauer graph algebra and $\mathcal{S}$  an orthogonal system which contains  at least one non-periodic module  on $A$-$\stmod$.  As a byproduct of our main result, we show that  $\mathcal{S}$ extends to a simple-minded system and  the extension closure $\mathcal{F}(\mathcal{S})$ of $\mathcal{S}$ is functorially finite on $A$-$\stmod$.  Moreover, we show that  a weakly simple-minded system with a finite cardinality is a simple-minded system on $A$-$\stmod$. 

 This paper is organized as follows. In Section 2, we recall simple-minded systems, Brauer graph algebras and some related concepts and conclusions. In Section 3, we study the stable bi-perpendicular category of a stable brick on $A$-$\stmod$ over a 2-domestic Brauer graph algebra. In Section 4,  we  construct orthogonal systems on $A$-$\stmod$.  In Section 5,  we  prove our main result Theorem \ref{sms-BGA}.
 
\section{Preliminary}	
Let $A$ be a finite-dimensional $k$-algebra and $\Gamma_{A}$  the Auslander-Reiten quiver of $A$.  An $A$-module $X$ is said to be {\bf $\tau$-periodic module} if there is a positive integer $m$ such that $X\cong\tau^{m}(X)$, otherwise, $X$ is called a {\bf non-$\tau$-periodic module}.  We say two $A$-modules $X$ and $Y$ are {\bf orthogonal} on  $A$-$\m$ if $\Hom_{A}(X,Y)=0$ and $\Hom_{A}(Y,X)=0$. $X$ and $Y$ are called  {\bf stable orthogonal} on  $A$-$\stmod$ if $\StHom_{A}(X,Y)=0$ and $\StHom_{A}(Y,X)=0$. For two families $\mathcal{X}$ and $\mathcal{Y}$ of objects  are called {\bf mutual stable orthogonal}  on  $A$-$\stmod$, if $\StHom_{A}(X,Y)=0$ and $\StHom_{A}(Y,X)=0$ for any $X\in\mathcal{X}$ and $Y\in\mathcal{Y}$.

By a {\bf component} of $\Gamma_{A}$ we mean a connected component of $\Gamma_{A}$. A component $\mathcal{C}$ of $\Gamma_{A}$ is said to be {\bf acyclic} if $\mathcal{C}$  has no oriented cycles. A subquiver $\Sigma$ of a component is said to be  {\bf right stable} (resp. {\bf left stable}) if $\tau^{-1}$ (resp. $\tau$) is defined on all modules over $\Sigma$. A subquiver $\Sigma$ is {\bf closed under successors} (resp. {\bf closed under predecessors}) on $\Gamma_A$ means that if there is a path whose start point (resp. end point) is on $\Sigma$, then the path belongs to $\Sigma$.

A connected component $\Gamma$ is said to be a {\bf stable tube} if it is of the form $\mathbb{Z}A_{\infty}/(\tau^{r})$ for some integer $r\geq1$. By a {\bf quasi-tube},  we mean a translation quiver $\Gamma$ such that its full translation subquiver formed by all vertices which are not projective-injective is a tube.  We say $\Gamma$ is a {\bf homogeneous tube} provided that  $\Gamma$ is of the form $\mathbb{Z}A_{\infty}/(\tau)$.
Skowronski proved \cite{S1} that every generalized standard family $\mathcal{C}$ of components on $\Gamma_{A}$ is almost periodic, that is, all but finitely many $\tau_{A}$-orbits on $\mathcal{C}$ are periodic. In particular, for a self-injective algebra $A$, every infinite generalized standard component $\mathcal{C}$ of $\Gamma_{A}$ is either acyclic with finitely many $\tau_{A}$-orbits or a quasi-tube. A component $\mathcal{C}$ is called an {\bf Euclidean component} if $\mathcal{C}$ is of the form $\mathbb{Z}\Delta$, where $\Delta$ is an   Euclidean quiver in  $\{\widetilde{A}_{n}(n\geq1),\widetilde{D}_{m}(m\geq4),\widetilde{E}_{6},\widetilde{E}_{7},\widetilde{E}_{8}\}$.

 A family of connected components  $\{\mathcal{C}_{i}\}_{i\in J}$ ($J$ an index set)  of $\Gamma_A$ is called  {\bf generalized standard}, if given  $X\in\mathcal{C}_{i}$ and  $Y\in\mathcal{C}_{j}$, we have $\rad^{\infty}(X,Y)=0$, where $\rad^{\infty}(-,-)$ means the intersection of $\rad^{t}(-,-)$ for $t\geq1$.  In particular, a full translation subquiver $\Sigma$ of  $\Gamma_A$  is called {\bf generalized standard} (\cite{S1}), if  $\rad^{\infty}(X,Y)=0$ for any $X$ and $Y$ on $\Sigma$.  A connected component $\mathcal{C}$ is called  {\bf stable generalized standard}, if $\underline{\rad}^{\infty}(X,Y)=0$  for all $X$ and $Y$ on $\mathcal{C}$, where $\underline{\rad}^{\infty}(X,Y)=\rad^{\infty}(X,Y)/P(X,Y)$ and $P(X,Y)$ is the subspace of $\Hom_{A}(X,Y)$ which factors through a projective module.  Note that, if a connected component $\mathcal{C}$ is generalized standard, then $\mathcal{C}$ is stable generalized standard, but the converse is not true.

Let $\Delta=(\Delta_{0},\Delta_{1})$ be a connected and acyclic quiver. Recall that the {\bf infinite translation quiver}  $(\mathbb{Z}\Delta,\tau)$ is defined as follows. The set of points of $\mathbb{Z}\Delta$ is $(\mathbb{Z}\Delta)_{0}=\mathbb{Z}\times\Delta_{0}=\{(n,z)\mid n\in\mathbb{Z}, z\in\Delta_{0}\}$. For an arrow $\alpha\colon\xymatrix{x\ar[r]^{a}& y}$ in $\Delta_{1}$, there are two arrows
\[(i,\alpha):\xymatrix{(i,x)\ar[r]&(i,y),}\ (i,\alpha'):\xymatrix{(i-1,y)\ar[r]&(i,x)}\]
in $(\mathbb{Z}\Delta)_{1}$ and those are all arrows in $(\mathbb{Z}\Delta)_{1}$. The translation $\tau$ of $\mathbb{Z}\Delta$ is defined by $\tau(i,x)=(i-1,x)$ for $(i,x)\in(\mathbb{Z}\Delta)_{0}.$ Let $\Gamma_{A}$ be AR-quiver of a representation-infinite algebra $A$. In order to distinguish different connected components of $\Gamma_{A}$, we introduce two notations to label vertices on different components of the AR-quiver. 

For an acyclic component $\mathcal{C}$ on $\Gamma_{A}$, we use  notation $(i,j,k)$ to denote a non-projective vertex on $\mathcal{C}$, where $i$ also distinguishes different acyclic components of AR-quiver, $j$ and $k$ are  determined by $\mathbb{Z}\times\Delta_{0}.$
For example, we draw a picture of  $\mathbb{Z}A^{\infty}_{\infty}$  as follows. 

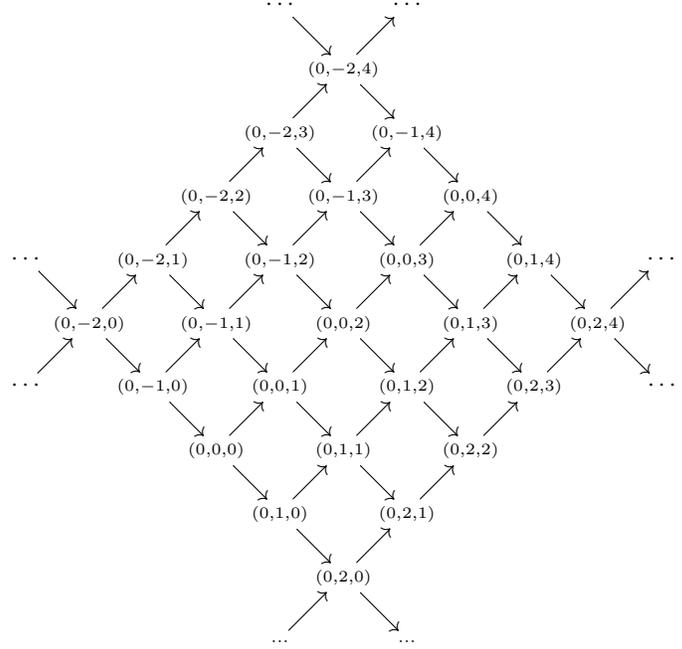
\begin{small}\[\xymatrix@dr@R=17pt@C=17pt@!0{
&&\cdots&&&&&&&&\cdots\\
&&\\
\cdots\ar[rr]&&\scriptstyle(0,-2,4)\ar[uu]\ar[rr] &&\scriptstyle(0,-1,4) \ar[rr]& &\scriptstyle (0,0,4) \ar[rr] && \scriptstyle(0,1,4)\ar[rr]&& \scriptstyle(0,2,4) \ar[rr]\ar[uu]&&\ \  \cdots\,.& \\	
&& \\	
&&\scriptstyle(0,-2,3)\ar[rr]\ar[uu] &&\scriptstyle(0,-1,3)\ar[rr] \ar[uu]&& \scriptstyle (0,0,3)\ar[rr]\ar[uu] &&\scriptstyle(0,1,3)\ar[uu]\ar[rr] &&\scriptstyle(0,2,3)\ar[uu]&&\\	
&&& \\	
&&	\scriptstyle(0,-2,2)\ar[rr]\ar[uu] &&\scriptstyle(0,-1,2)\ar[rr]\ar[uu]& & \scriptstyle(0,0,2)\ar[rr]\ar[uu] &&\scriptstyle(0,1,2)\ar[rr]\ar[uu] && \scriptstyle (0,2,2)\ar[uu]\\	
&& \\
&&\scriptstyle(0,-2,1) \ar[rr]\ar[uu]&&\ar[rr] \scriptstyle (0,-1,1)\ar[uu] && \scriptstyle(0,0,1)\ar[uu]\ar[rr]\scriptstyle &&\scriptstyle(0,1,1)\ar[uu] \ar[rr]&&\scriptstyle(0,2,1)\ar[uu]\\
&&& \\
\cdots\ar[rr]&&\scriptstyle(0,-2,0)\ar[rr]\ar[uu]  && \scriptstyle(0,-1,0)\ar[uu]\ar[rr] && \scriptstyle (0,0,0)\ar[uu]\ar[rr]  &&\scriptstyle(0,1,0) \ar[rr]\ar[uu]&&\scriptstyle(0,2,0)\ar[uu]\ar[rr]&& \scriptstyle\cdots\\ 
&&\\
&&\cdots\ar[uu]&&&&&&&&\scriptstyle\cdots\ar[uu]
}\] 
\end{small}
Note that the AR-translation $\tau$  on  $\mathcal{C}$ satisfies
$\tau(i,j,k)=(i,j-1,k-1)$.

For a quasi-tube $\mathcal{Q}$ (its stable part is of the form $\mathbb{Z}A_{\infty}/(\tau^{n})$ for a positive integer $n$) on $\Gamma_{A}$, we use the notation $Q(i,j,k)$ to denote a non-projective vertex on $\mathcal{Q}$, where $i$ distinguishes different quasi-tubes of $\Gamma_{A}$, $j$ and $k$ are  determined by $\mathbb{Z}\times\Delta_{0}.$
For example, we draw a picture of  $\mathbb{Z}A_{\infty}$ as follows.
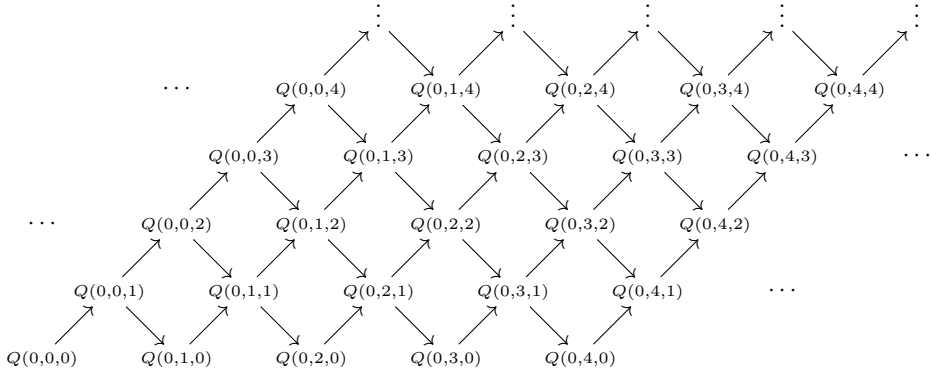
\vspace{-4cm}
\begin{small}
\[\xymatrix@dr@R=18pt@C=18pt@!0{
&&&&&&&&&&\vdots&&\\
&&\\
&&&&&&&&\vdots\ar[rr]&&\scriptstyle Q(0,4,4)\ar[uu]&&\cdots\\
&&\\
&&&&&&\vdots\ar[rr]&&\scriptstyle Q(0,3,4)\ar[rr]\ar[uu]&&\scriptstyle Q(0,4,3)\ar[uu]&&\\
&&\\
&&&&\vdots\ar[rr]&&\scriptstyle Q(0,2,4)\ar[rr]\ar[uu]&&\scriptstyle Q(0,3,3)\ar[rr]\ar[uu]&&\scriptstyle Q(0,4,2)\ar[uu]&&\cdots\\
&&\\
&&\vdots\ar[rr] &&\scriptstyle Q(0,1,4) \ar[rr]\ar[uu]& &\scriptstyle Q(0,2,3) \scriptstyle \ar[rr]\ar[uu] && \scriptstyle Q(0,3,2)\ar[rr]\ar[uu]&& \scriptstyle Q(0,4,1) \ar[uu]&&\\	
&& \\	
&&\scriptstyle Q(0,0,4) \ar[uu]\ar[rr]&&\scriptstyle Q(0,1,3)\ar[uu]\ar[rr] &&\scriptstyle Q(0,2,2) \ar[uu]\ar[rr] &&\scriptstyle Q(0,3,1)\ar[uu]\ar[rr]&&\scriptstyle Q(0,4,0)\ar[uu]\\	
&& \\	
\cdots &&\scriptstyle Q(0,0,3)\ar[uu]\ar[rr] &&\scriptstyle Q(0,1,2)\ar[uu]\ar[rr]  &&\scriptstyle Q(0,2,1)\ar[uu]\ar[rr] &&\scriptstyle Q(0,3,0)\ar[uu]  && \\	
&&& \\
&&\scriptstyle Q(0,0,2) \ar[rr]\ar[uu]  && \scriptstyle Q(0,1,1) \ar[uu]\ar[rr] &&\scriptstyle Q(0,2,0)\ar[uu]  \\
&&& \\
\cdots&&\scriptstyle Q(0,0,1) \ar[uu]\ar[rr] && \scriptstyle Q(0,1,0)\ar[uu]\\	
&& \\
&&\scriptstyle Q(0,0,0) \ar[uu] 
\\}\] 
\end{small}
Note that the AR-translation $\tau$  satisfies
$\tau Q(i,j,k)=Q(i,j-1,k)$ on the component $\mathcal{Q}$. 
 Please compare with the Euclidean component case.  We remind the reader that when we consider a general connected component $\mathcal{D}$ of a AR-quiver, we may still use a pair $(i,j)$ of integers to label a vertex on $\mathcal{D}$.
 
\subsection{Simple-minded system}
Let $\mathcal{T}$ be a  triangulated category with the shift functor $[1]$.   For two families $\mathcal{S}_{1}, \mathcal{S}_{2}$ of objects in $\mathcal{T}$, we define 
\[\mathcal{S}_{1}\star\mathcal{S}_{2}:=\{ X\in \mathcal{T}\mid \mbox{ There is a  triangle }S_{1} \longrightarrow  X \longrightarrow S_{2} \longrightarrow S_{1}[1], where \ S_{1}\in \mathcal{S}_{1}, S_{2}\in \mathcal{S}_{2}\}. \]
It is routine to check that the operator $\star$ is associative, that is,  $(\mathcal{S}_{1}\star\mathcal{S}_{2})\star\mathcal{S}_{3}=\mathcal{S}_{1}\star(\mathcal{S}_{2}\star\mathcal{S}_{3})$ for $\mathcal{S}_{1}, \mathcal{S}_{2}$ and $\mathcal{S}_{3}\subseteq \mathcal{T}$.  One denotes $(\mathcal{S})_{0}=\{0\}$ and inductively defines $(\mathcal{S})_{n}=(\mathcal{S})_{n-1}\star(\mathcal{S}\cup\{0\})$ for $n\in\mathbb{Z}^{+}$.  We say that $\mathcal{S}$ is {\it extension-closed}, if $\mathcal{S}\star\mathcal{S}\subseteq \mathcal{S}$. We denote the {\bf extension closure} of a family $\mathcal{S}$ of objects in $\mathcal{T}$ as $$\mathcal{F}(\mathcal{S}):=\bigcup_{n\geq0}(\mathcal{S})_{n},$$ which is the smallest extension-closed full subcategory of  $\mathcal{T}$  containing $\mathcal{S}$.

\begin{Def}\label{brick-orthogonal-system}
Let $\mathcal{T}$ be an additive $k$-category.  An object $M$ in $\mathcal{T}$ is a {\bf stable brick} if $\mathcal{T}(M,M)\cong k$.   Moreover, a family $\mathcal{M}$ of stable bricks in $\mathcal{T}$  is an {\bf orthogonal system} if $\mathcal{T}(M,N)=0$ for all distinct  $M, N$ in $\mathcal{M}$.
\end{Def}

\begin{Rem}\label{F-S-properties}
If $\mathcal{S}$  is an orthogonal system, then $\mathcal{F}(\mathcal{M})$ is  closed under direct summands. Please refer to \cite[Lemma 2.7]{Dugas} for more details.
\end{Rem}

\begin{Def}\label{definition-sms-right-tir} {\rm(\cite[Definition 2.1]{KL},\cite[Definition 2.4, 2.5]{Dugas})} 
Let $\mathcal{T}$ be a  triangulated category.  A family of objects $\mathcal{S}$ in $\mathcal{T}$ is a {\bf simple-minded system} if the following two conditions are satisfied$\colon$
\begin{enumerate}[$(1)$]
\item {\rm(Orthogonality)} $\mathcal{S}$ is an orthogonal system in $\mathcal{T}$. 
\item {\rm(Generating condition)} Extension closure $\mathcal{F}(\mathcal{S})$ of $\mathcal{S}$ is equal to $\mathcal{T}$.
\end{enumerate}
\end{Def}

Koenig and Liu \cite{KL} introduced a weaker concept than simple-minded system, namely weakly simple-minded system.
\begin{Def}\label{definition-wsms-right-tir} {\rm(\cite[Definition 5.3]{KL})} 
Let $\mathcal{T}$ be a triangulated category.  A family of objects $\mathcal{S}$ in $\mathcal{T}$ is a {\bf weakly simple-minded system}  if the following two conditions are satisfied$\colon$
\begin{enumerate}[$(1)$]
\item {\rm(Orthogonality)} $\mathcal{S}$ is an orthogonal system in $\mathcal{T}$. 
\item {\rm(Weakly generating condition)} For any non-zero object $X$ in $\mathcal{T}$, there is an object $S$ in $\mathcal{S}$ such that $\mathcal{T}(S,X)\ncong 0.$
\end{enumerate}
\end{Def}

\begin{Them}$($\cite[Theorem 1.3]{CLZ}$)$\label{simple-module-and-n-tube}
Let $A$ be a self-injective algebra and $\mathcal{C}$ a quasi-tube of rank $n$. Then the number of elements in a simple-minded system of $A$ lying in $\mathcal{C}$ is strictly less than $n$. In particular, none of the indecomposable modules in a simple-minded system  lie in the homogeneous tubes of the AR-quiver.
\end{Them}	

Let $\mathcal{T}$ be an additive category and $\mathcal{X}$ a full subcategory of $\mathcal{T}$. We say that a morphism $f\colon M\rightarrow X$ is a {\bf left $\mathcal{X}$-approximation} of $M$ if $X\in\mathcal{X}$, and the abelian group homomorphism 
$$\mathcal{T}(f,X')\colon \mathcal{T}(X,X')\rightarrow \mathcal{T}(M,X')$$ is surjective for each $X'\in\mathcal{X}$. Dually, a morphism $g\colon Y\rightarrow N$ is a {\bf right $\mathcal{X}$-approximation} of $N$ if $Y\in\mathcal{X}$, and the abelian group homomorphism 
$$\mathcal{T}(Y',g)\colon \mathcal{T}(Y',Y)\rightarrow \mathcal{T}(Y',N)$$ is surjective for each $Y'\in\mathcal{X}$. 

Let $\mathcal{S}$ be a full subcategory of $\mathcal{T}$ containing $\mathcal{X}$. $\mathcal{S}$ is said to be {\bf covariantly finite} (resp. {\bf contravariantly finite}) in $\mathcal{X}$ if every $S\in\mathcal{S}$ has a left (resp. right) $\add(\mathcal{X})$-approximation. We say $\mathcal{S}$  is  {\bf functorially finite} in  $\mathcal{X}$ provided that  $\mathcal{S}$ is both covariantly finite and contravariantly finite in $\mathcal{X}$. The following theorem provided us a necessary condition for an extension closure  to be functorially finite.

\begin{Them}$($\cite[Theorem 3.3]{Dugas}$)$\label{subset-of-sms}
Let $\mathcal{T}$ be a Hom-finite Krull-Schmidt triangulated category. Suppose $\mathcal{X}\subseteq\mathcal{S}$ for a simple-minded system $\mathcal{S}$ in $\mathcal{T}$. Then $(^{\perp}\mathcal{X},\mathcal{F}(\mathcal{X}))$ and $(\mathcal{F}(\mathcal{X}),\mathcal{X}^{\perp})$ are torsion pairs in $\mathcal{T}$.  In particular, $ \mathcal{F}(\mathcal{X})$ is a functorially finite subcategory of $\mathcal{T}$.
\end{Them}

\subsection{Brauer graph algebras}
We state the definition of Brauer graph algebras and some related results as follows. Please refer to \cite[Section 2]{S0} for more details. 

\begin{Def}{\rm(\cite[Section 2]{S0})}\label{Brauer-graph}
A {\bf  Brauer graph} $G=(G_{0},G_{1}, m, o)$ is a finite unoriented connected graph,  where $G_{0}$ is the set of its vertices, $G_{1}$ is the set of edges, the map $m\colon G_{0}\rightarrow \mathbb{Z}_{>0}$ is called multiplicity, $o$ is the given orientation, that is, for each vertex $v\in G_{0}$,  each edge incident to vertex  $v$ is  given  a cyclic ordering. We always assume that the cyclic orderings are clockwise. Suppose that the cyclic ordering at vertex $v$ is given by $i_1<i_2<\cdots<i_n$, then we say  $i_j$ is the predecessor of $i_{j+1}$, and  $i_{j+1}$ is the successor of $i_j$. 
\end{Def}

Given a Brauer graph $G=(G_{0},G_{1}, m, o)$, we denote  {\it valency} of $v\in G_{0}$ by $val(v)$, which is defined to be the number of edges in $G$ incident to $v$. Note that a loop is counted twice. 
Every Brauer graph determines a finite dimensional basic symmetric algebra, it is called  a {\bf Brauer graph algebra}. Let $A$ be the Brauer graph algebra determined by a Brauer graph $G$. There are a quiver $Q_{G}$ and an admissible ideal $I_{G}$ such that $A\cong kQ_{G}/I_{G}$.  We refer to \cite[Section 2]{S0} for the construction of $Q_{G}=(Q_0,Q_1)$ and $I_{G}$. We provide a concrete construction in the following example. 

\begin{Ex}\label{BGA-EX}
Take a Brauer graph  $G=(G_{0},G_{1}, m, o)$ as follows. 
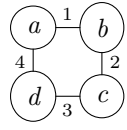
\[\xymatrix@r@R=26pt@C=26pt@!0{
*++[o][F]\txt{a} \ar@{-}[r]^{1}  \ar@{-}[d]_{4} &*++[o][F]\txt{b}\ar@{-}[d]^{2}	   \\
*++[o][F] \txt{d}\ar@{-}[r]_{3} & *++[o][F] \txt{c}  \\
}\]
The multiplicity of each vertex  $i\in G_{0}$ is $m(i)=1$. Cyclic ordering of the edge incident with vertex $a$ is $1<4<1$, with vertex $b$ is  $2<1<2$, with vertex  $c$ is  $3<2<3$, with vertex  $d$ is  $4<3<4$.  The valencies of vertices are as follows: $val(a)=val(b)=val(c)=val(d)=2$.
	
We now present the algebra $A\cong kQ_{G}/I_{G}$ determined by Brauer graph $G$ as follows. We first consider the corresponding quiver $Q_{G}=(Q_{0},Q_{1})$.  
The set of $Q_{0}$ is given by the set of edges $G_{1}$ of $G$.  We denote the vertex in $Q_{0}$ corresponding to the edges $j$ in $G_{1}$ also by $j$ for $j=1,2,3$, $4$. By the above paragraph, the cyclic orderings in $G$ determine  arrows $Q_{1}$ of $Q$ in the following way: $\alpha_{1}:1\to2$, $\beta_{1}:2\to1$,   $\alpha_{2}:2\to3$, $\beta_{2}:3\to2$, $\alpha_{3}:3\to4$, $\beta_{3}:4\to3$, and  $\alpha_{4}:4\to1$, $\beta_{4}:1\to4$.	Hence the quiver $Q_{G}$ is the following diagram. 
	
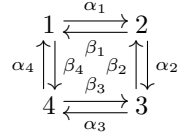
\[\xymatrix@C=.75pc@R=1.5pc{
1 \ar@<+.5ex>[rr]^{\alpha_{1}} \ar@<+.5ex>[d]^{\beta_{4}} \save[] \restore  && 2 \ar@<+.5ex>[d]^{\alpha_{2}}\ar@<+.5ex>[ll]^{\beta_{1}}\save[]  \restore \\
4 \ar@<+.5ex>[u]^{\alpha_{4}} \ar@<+.5ex>[rr]^{\beta_{3}} \save[]	\restore 
&& 3  \ar@<+.5ex>[ll]^{\alpha_{3}} \ar@<+.5ex>[u]^{\beta_{2}}  \save[]  \restore
}\]	
	
Then we consider the ideal $I_{G}$ of  $kQ_{G}$ generated by three types of relations. 

Relations of type I: $\beta_{1}\alpha_{1}-\alpha_{4}\beta_{4}$, $\beta_{2}\alpha_{2}-\alpha_{1}\beta_{1}$, $\beta_{3}\alpha_{3}-\alpha_{2}\beta_{2}$,  $\beta_{4}\alpha_{4}-\alpha_{3}\beta_{3}$.

Relations of type II:   $\alpha_{1}\beta_{1}\alpha_{1}$, $\beta_{1}\alpha_{1}\beta_{1}$,  $\alpha_{2}\beta_{2}\alpha_{2}$, $\beta_{2}\alpha_{2}\beta_{2}$,  $\alpha_{3}\beta_{3}\alpha_{3}$, $\beta_{3}\alpha_{3}\beta_{3}$,  $\alpha_{4}\beta_{4}\alpha_{4}$, $\beta_{4}\alpha_{4}\beta_{4}$.	

Relations of type III: $\alpha_{3}\alpha_{2}\alpha_{1}$, $\alpha_{4}\alpha_{3}\alpha_{2}$, $\alpha_{1}\alpha_{4}\alpha_{3}$, $\alpha_{2}\alpha_{1}\alpha_{4}$, $\beta_{3}\beta_{4}\beta_{1}$,  $\beta_{4}\beta_{1}\beta_{2}$, $\beta_{1}\beta_{2}\beta_{3}$, $\beta_{2}\beta_{3}\beta_{4}$. 
	
\medskip
	
Thus $A\cong kQ_{G}/I_{G}=~~\begin{matrix}1\\\begin{matrix}4\end{matrix}~~\begin{matrix}2\end{matrix}\\1\end{matrix}\oplus~~\begin{matrix}2\\\begin{matrix}1\end{matrix}~~\begin{matrix}3\end{matrix}\\2\end{matrix}\oplus~~\begin{matrix}3\\\begin{matrix}2\end{matrix}~~\begin{matrix}4\end{matrix}\\3\end{matrix}\oplus~~\begin{matrix}4\\\begin{matrix}3\end{matrix}~~\begin{matrix}1 \end{matrix}\\4\end{matrix}.$
\end{Ex}

We briefly recall the representation type of an algebra. From the point of view of representation type, according to Drozd's classification \cite{D}, an algebra is either tame or wild type, and not both. Let $k[x]$ be the polynomial algebra in one variable. An algebra $A$ is said to be {\it tame}, if for any dimension $d$, there exist finitely many  $k[x]$-$A$-bimodules $M_i$, $1\leq i\leq n_d$, which are finitely generated and free as left $k[x]$-module, and all but finitely many isomorphism classes of indecomposable modules of dimension $d$ in $ind\,A$ are of the form  $k[x]/(x-\lambda)\otimes_{k[x]}M_i$, for some $\lambda\in k$, and some $i\in\{1,2,\cdot\cdot\cdot,n_d\}$. Let $\mu_A(d)$ be the least number of $k[x]$-$A$-bimodules satisfying the above condition for $d$. Then $A$ is said to be of {\it polynomial growth}  (cf. \cite{SK}), if there is a positive integer $m$ such that  $\mu_A(d)<d^m$ for all $d>1$.  We say $A$  is  of {\it domestic} type (cf. \cite{CB}), if there is a positive integer  $m$ such that $\mu_A(d)<m$ for all $d>1$ and $A$ is said to be {\it d-domestic} if $d$ is the least such integer $m$. 

\begin{Them}$($\cite{BoS}$)$\label{domestic-Brauer-graph-algebra}
Let $A$ be a Brauer graph algebra with Brauer graph $G$. Then
\begin{enumerate}[$(1)$]
\item $A$  is $1$-domestic if and only if one of the following holds:
\begin{enumerate}[$(a)$]
\item $G$ is a tree with $m(i)=2$ for exactly two vertices $i_{0}, i_{1}\in G_{0}$ and $m(i)=1$ for all $i\in G_{0}, i\neq i_{0}, i_{1}$.
\item There is exactly one cycle in $G$ and this cycle has odd length, and all multiplicities of edges are one, that is, $m\equiv1$.
\end{enumerate}
\item $A$ is $2$-domestic if and only if  there is exactly one cycle in $G$ and this cycle has even length, and $m\equiv1$.
\item There is no $n$-domestic  Brauer graph algebra for $n\geq3$.
\end{enumerate}
\end{Them}

\begin{Them}$($\cite[Theorem 4.4, Corollary 4.5]{D1}$)$\label{two-domestic-Brauer-graph-algebras} 
Let $A$ be a representation-infinite domestic Brauer graph algebra with Brauer graph $G$ with $n$ edges. If $G$ has a cycle, then it is unique. Let $n_{1}$ be the number of {\rm(}additional{\rm)} edges on the inside of the cycle and $n_{2}$ the number of {\rm(}additional{\rm)} edges on the outside of the cycle.
In the notation above, 
\begin{enumerate}[$(1)$]
\item If $A$ is $1$-domestic, then $m=1$ and $p+q=2n$. Furthermore,
\begin{enumerate}[$(a)$]
\item If $G$ is a tree, then $p=q=n$. 
\item If $G$ has exactly one cycle in $G$ and this cycle has {\rm(}odd{\rm)} length $\ell$, then $p=\ell+2n_{1}$ and $q=\ell+2n_{2}$.
\end{enumerate}
\item If $A$ is $2$-domestic and the unique cycle of $Q$ has {\rm(}even{\rm)} length $\ell$, then $m=2$ and $p=\ell/2+n_{1}$ and $q=\ell/2+n_{2}$ such that $p+q=n$.
\end{enumerate}
\end{Them}

Note that the Brauer graph algebra in Example \ref{Brauer-graph} is 2-domestic.   Let $A$ be a 2-domestic Brauer graph algebra.  It follows from  \cite[Chapter 2, Theorem 7.1]{E} that  AR-quiver $\Gamma_{A}$ of $A$  consists of two stable Euclidean components $\Gamma_{0}$ and $\Gamma_{1}$ of the form $\mathbb{Z}\widetilde{A}_{p,q}$, two quasi-tubes $Q_{0}$ and $Q_{1}$ of rank $q$ and  two quasi-tubes $Q'_{0}$ and $Q'_{1}$ of rank $p$, as well as infinitely many homogeneous tubes. By Theorem  \ref{two-domestic-Brauer-graph-algebras}, $p+q=n$.  By $\widetilde{A}_{p,q}(p,q\geq1)$, we mean the quiver 
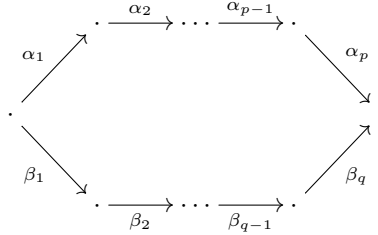
\[\xymatrix{
	& \cdot \ar[r]^-{\alpha_{2}}  & \cdots\ar[r]^-{\alpha_{p-1}} &\cdot \ar[dr]^-{\alpha_{p}} \\
	\cdot \ar[ur]^-{\alpha_{1}} \ar[dr]_-{\beta_{1}} &  & & &\cdot \\	
	&  \cdot \ar[r]_-{\beta_{2}} & \cdots\ar[r]_-{\beta_{q-1}}  &\cdot\ar[ur]_-{\beta_{q}} 
}\]

\subsection{The property of components over a 2-domestic Brauer graph algebra}
We first study the properties Euclidean components over a 2-domestic Brauer graph algebra. We recall some results about acyclic connected components  as follows.
\begin{Prop}{\rm(\cite[Corollary 5.7]{SY}}\label{symm-generali-standard} 
Let $A$ be a symmetric $k$-algebra. Then the following conditions are equivalent. 
\begin{enumerate}[$(1)$]
\item $\Gamma_A$ admits an acyclic, generalized standard, right stable full translation subquiver which is closed under successors on $\Gamma_A$. 
\item $\Gamma_A$ admits an acyclic, generalized standard, left stable full translation subquiver which is closed under predecessors on $\Gamma_A$. 
\item $A$ is  isomorphic to trivial extension algebra $T(B)$, where $B$ is a tilted algebra not of Dynkin type. 
\end{enumerate}
\end{Prop}

\begin{Them}$($\cite[Theorem 2]{BoS}$)$\label{symmetric-algebras-of-Euclidean-type-2}
Let $A$ be a basic indecomposable algebra. Then the following conditions are equivalent.
\begin{enumerate}[$(1)$]
\item $A$ is isomorphic to  a  2-domestic Brauer graph algebra. 	
\item $A$ is a symmetric algebra of Euclidean type $\widetilde{\mathbb{A}}_m$ and the Cartan matrix of $A$ is singular.
\item $A$ is isomorphic to the trivial extension $T(B)$, where $B$ is a  representation-infinite tilted algebra of Euclidean type $\widetilde{\mathbb{A}}_m$.		
\end{enumerate}
\end{Them}

\begin{Cor}\label{2-BGA-bricks}
Let $A$ be a 2-domestic Brauer graph algebra and  let  $X$ and $Y$ be two objects on  an Euclidean component.  If there is no compositions of finitely many irreducible maps from $X$ to $Y$, then $\StHom_{A}(X,Y)\cong 0$.  In particular, each Euclidean component of AR-quiver  $\Gamma_{A}$ is stable generalized standard and every indecomposable, non-projective non-periodic module is a stable brick on $A$-$\stmod$. 
\end{Cor}
\begin{proof}
Since  $A$ is 2-domestic, there are two Euclidean components, denoted by $\Gamma_{0}$ and $\Gamma_{1}$. Note that $\Gamma_{0}$ and $\Gamma_{1}$ have only finitely many $\tau$-orbits and $\Omega(_{s}\Gamma_{0})={_{s}\Gamma_{1}}$. By Proposition \ref{symm-generali-standard}  and Theorem \ref{symmetric-algebras-of-Euclidean-type-2}, there is an acyclic, generalized standard, right stable full translation subquiver which is closed under successors on $\Gamma_A$.  Without loss of generality, we assume that $\Gamma_{0}$ admits an acyclic, generalized standard, right stable full translation subquiver $\Sigma$ which is closed under successors $\Gamma_0$. Note that $\rad^{\infty}(M,N)=0$ for any $M$ and $N$ on $\Sigma$.

Let  $X$ and $Y$ be two non-projective objects on $\Gamma_{0}$. We assume that there is no compositions of finitely many irreducible maps from $X$ to $Y$ on  $\Gamma_{0}$. Thus $\Hom_{A}(X,Y)\subseteq\rad^{\infty}(X,Y)$. 
If both $X$ and $Y$ are in $\Sigma$, $\Hom_{A}(X,Y)\cong 0,$ thus $\StHom_{A}(X,Y)\cong 0.$ Otherwise, one of $X$ and  $Y$ is not  in $\Sigma$.  Since $X$ and  $Y$ are not projective and $A$ is symmetric, the $\tau$-orbits of both $X$ and  $Y$  contain no projective or injective modules. Since $\Sigma$  is closed under successors on  $\Gamma_0$, $\Sigma$ contains at least one object for each $\tau$-orbits of $\Gamma_0$.   Thus there is a large number $m$ such  that both $\tau^{-m}(X)$ and  $\tau^{-m}(Y)$ are contained in $\Sigma$. Since there is no compositions of finitely many irreducible maps from $X$ to $Y$, so are $\tau^{-m}(X)$ and  $\tau^{-m}(Y)$. Thus  $\StHom_{A}(X,Y)\cong\StHom_{A}(\tau^{-m}(X),\tau^{-m}(Y))\cong0.$
Thus our conclusion holds for $\Gamma_0$. Since  $\Omega(_{s}\Gamma_{0})={_{s}\Gamma_{1}}$ and $\Omega$ is an equivalent functor on $A$-$\stmod$, our conclusion also holds for $\Gamma_1$. 
Therefore each Euclidean component is stable generalized standard. Since both $\Gamma_0$ and $\Gamma_1$ are acyclic, every indecomposable, non-projective non-periodic module is a stable brick on $A$-$\stmod$.
\end{proof}

Then we start to study quasi-tubes of a 2-domestic Brauer graph algebra.  We recall some conclusions for quasi-tubes as follows. 
\begin{Prop}{\rm(}\cite[Theorem A]{MS}{\rm)}\label{quasi-tube}
Let $A$	be a self-injective algebra and $\Gamma$ a connected component of $\Gamma_{A}$. The following are equivalent.
\begin{enumerate}[$(1)$]
	\item $\Gamma$ is  a quasi-tube.
	\item $_{s}\Gamma$ is  a  stable tube.
	\item $\Gamma$ contains an oriented cycle.
\end{enumerate}
\end{Prop}

\begin{Them}$($\cite[Section 4, 4.17]{S2}$)$\label{2-BGA-canonical}
Let $A$ be a self-injective algebra. Then the following statements are equivalent.
\begin{enumerate}[$(1)$]
\item $A$ is a symmetric algebra of Euclidean type and  has singular Cartan matrix. 	
\item $A$ is derived equivalent to trivial extension 
$T(C)$  of a canonical algebra $C$ of Euclidean type.
\item $A$ is stably equivalent to  the trivial extension $T(C)$ of a canonical algebra of Euclidean type.		
\end{enumerate}
\end{Them}

\begin{Them}$($\cite[Corollary  2]{JKS}$)$\label{2-BGA-canonical-2}
Let $A$ be a basic, indecomposable finite dimensional symmetric $k$-algebra. Then the following statements are equivalent.
\begin{enumerate}[$(1)$]
\item $\Gamma_{A}$ admits a generalized standard family $\mathcal{C}=(\mathcal{C}_{\lambda})_{\lambda\in\Lambda}$ of quasi-tubes maximally saturated by simple modules and projective modules. 	
\item $\Gamma_{A}$ admits a generalized standard family $\mathcal{C}=(\mathcal{C}_{\lambda})_{\lambda\in\Lambda}$ of  stable tubes maximally saturated by simple modules. 	
\item $A$ is isomorphic  to  the trivial extension $T(B)$ of a canonical algebra $B$.		
\end{enumerate}
\end{Them}
\begin{Rem}\label{generalized-orthogonal}
 For a quasi-simple $X$ on a component  $\mathcal{C}_{\lambda}$ of a generalized standard family  $\mathcal{C}$, $X$ and $Y$ are  orthogonal  on $A$-$\m$ for any quasi-simple which belongs to the family $\mathcal{C}$ of quasi-tube, except $\Omega(X)$, $\Omega^{-1}(X)$ and $X$.
\end{Rem}

\begin{Cor}\label{2-BGA-quasi-tube-bricks}
Let $A$ be a 2-domestic Brauer graph algebra and $X$ a quasi-simple on a quasi-tube. Then  $X$ and $Y$ are mutual stable orthogonal on $A$-$\stmod$ for any quasi-simple $Y$, except $\Omega(X)$, $\Omega^{-1}(X)$ and $X$. In particular, each quasi-tube  is stable generalized standard. 
\end{Cor}
\begin{proof}
By Theorem \ref{symmetric-algebras-of-Euclidean-type-2} and  Theorem \ref{2-BGA-canonical}, $A$ is stably equivalent to  the trivial extension $T(C)$ of a canonical algebra $C$ of Euclidean type.  Let $F\colon T(C)$-$\stmod\rightarrow A$-$\stmod$ be a stably equivalence.   Thus stable AR-quivers $_{s}\Gamma_{A}$ and $_{s}\Gamma_{T(C)}$ are isomorphic as stable translation quivers. Note that $F$ preserves quasi-tubes and quasi-simples. 

Since AR-quiver $\Gamma_{A}$ of $A$ contains two family $\mathcal{T}_{1}$ and $\mathcal{T}_{2}$  of quasi-tubes  satisfying  $\Omega(\mathcal{T}_{1})={\mathcal{T}_{2}}$. 
By Theorem \ref{2-BGA-canonical-2}, one of $\mathcal{T}_{1}$ and $\mathcal{T}_{2}$ is the image of a generalized standard family $\mathcal{C}=(\mathcal{C}_{\lambda})_{\lambda\in\Lambda}$ of quasi-tubes in $\Gamma_{T(C)}$ under stably equivalence $F$.  Without loss of generality, we assume that  $\mathcal{T}_{1}$ is the image of $\mathcal{C}$ under $F$,  and that $X$ is in $\mathcal{T}_{1}$. Take a quasi-simple  $Y$ in $\mathcal{T}_{1}$ which  is not isomorphic to $X$. By Remark \ref{generalized-orthogonal},  $\StHom_{A}(X,Y)\cong\StHom_{T(C)}(F(X),F(Y))\cong0$ and $\StHom_{A}(Y,X)\cong\StHom_{T(C)}(F(Y),F(X))\cong0$. Since  $\Omega(\mathcal{T}_{1})={\mathcal{T}_{2}}$, $\mathcal{T}_{2}$ also satisfies the above conclusion.

Take a quasi-simple $Y$ in $\mathcal{T}_{2}$ which is not isomorphic to  $\Omega(X)$. By Serre duality,
\[\StHom_{A}(X,Y)\cong\D\StHom_{A}(Y,\nu\Omega(X))\cong\D\StHom_{A}(Y,\Omega(X))\cong0.\]
Take a quasi-simple  $Y$ in $\mathcal{T}_{2}$ which is not isomorphic to  $\Omega^{-1}(X)$. By Serre duality,
\[\StHom_{A}(Y,X)\cong\D\StHom_{A}(X,\nu\Omega(Y))\cong\D\StHom_{A}(X,\Omega(Y))\cong0.\]

Take two objects $X$ and $Y$ on a quasi-tube. If there is no compositions of finitely many irreducible maps from $X$ to $Y$, by \cite[Section 2.2]{EK}, then $\StHom_{A}(X,Y)\cong 0$. Hence each quasi-tube is stable generalized standard. 
\end{proof}

\begin{Lem}\label{sectional-non-zero}
Let $A$ be a self-injective algebra with a stable generalized standard AR-component  $\Gamma$ and  $X$ an indecomposable non-projective module on $\Gamma$. If $X$ is a stable brick on $A$-$\stmod$, then any sectional path with indecomposable non-projective modules starting at $X$ on $\Gamma$ is non-zero on $A$-$\stmod$.
\end{Lem}
\begin{proof}
Take a sectional path starting at $X$ on $\Gamma$:
\[X=X_{1}\xrightarrow{f_{1}} X_{2}\xrightarrow{} \cdots\xrightarrow{}  X_{\ell-1}\xrightarrow{f_{\ell}}  X_{\ell}.\]
Note that each $f_{t}$ is an irreducible map on $A$-$\m$ for $t=1, 2		,\cdots,\ell$. We prove that the morphism  $\underline{f_{t}\cdots f_{2}f_{1}}$ is non-zero on $A$-$\stmod$ by induction on $t$.
 Since $X_{1}$ and $X_{2}$ are indecomposable  non-projective $A$-module and the morphism $f_{1}$ is irreducible,  $\underline{f_{1}}$ is not zero on $A$-$\stmod$. Therefore our conclusion holds for  $t=1$.

We assume that our conclusion holds for $t=k$. We prove  it is true for $t=k+1$.
On the contrary, we assume that the morphism $\underline{f_{k+1}\cdots f_{2}f_{1}}$ is zero on $A$-$\stmod$. Take the triangle induced by the almost split sequence ending at $X_{k+1}$ as follows.
 \[Z\xrightarrow{(\underline{h_{1}},\underline{h_{2}})^{t}}Y\oplus X_{k} \xrightarrow{(\underline{g},\underline{f_{k}})}  X_{k+1}\rightarrow Z[1].\]
By induction,  the morphism $\underline{f_{k}\cdots f_{2}f_{1}}$ is non-zero on $A$-$\stmod$. Take a morphism 
\[\alpha=(0,\underline{f_{k}\cdots f_{2}f_{1}})^{t}\colon X\xrightarrow{}Y\oplus X_{k}.\]
Then $(\underline{g},\underline{f_{k}})\alpha=(\underline{f_{k+1}f_{k}\cdots f_{2}f_{1}},0)=(0,0)$ on  $A$-$\stmod$.
Thus there is a morphism $\beta:X\xrightarrow{}Z$ such that $(\underline{h_{1}},\underline{h_{2}})\beta=\alpha$ on $A$-$\stmod$. See the following diagram.
 \[\xymatrix{
Z\ar[r]^-{(\underline{h_{1}},\underline{h_{2}})^{t}}&Y\oplus X_{k} \ar[r]^-{(\underline{g},\underline{f_{k}})}  &X_{k+1}\ar[r]&Z[1]\\
X\ar@{-->}[u]^-{\beta}\ar[ur]^-{\alpha}\\
 }.\]

There are two cases to be considered. If $\Gamma$ has no cycle, then there is no non-zero morphism from $X$ to $Z$, since $\Gamma$ is stable generalized standard. Thus $\alpha=(0,\underline{f_{k}\cdots f_{2}f_{1}})=0$ on $A$-$\stmod$. It contradicts the induction that the morphism $\underline{f_{k}\cdots f_{2}f_{1}}$ is non-zero on $A$-$\stmod$.  If $\Gamma$ has  a cycle, then $_{s}\Gamma$ is a quasi-tube by Proposition \ref{quasi-tube}. Thus our conclusion follows from \cite[Lemma 2.6]{EK}.
\end{proof}

We also state the dual result. 
\begin{Lem}\label{sectional-non-zero'}
Let $A$ be a self-injective algebra with a stable generalized standard AR-component  $\Gamma$ and let  $X$ be  an indecomposable non-projective module on $\Gamma$. If $X$ is a stable brick on $A$-$\stmod$, then any sectional path with indecomposable non-projective modules ending at $X$ on  $\Gamma$ is non-zero on $A$-$\stmod$.
\end{Lem}

\begin{Lem}\label{Euclidean-ses}
Let $A$ be a self-injective algebra and $(i,j)$  an indecomposable non-projective module on a connected component $\Gamma$ of AR-quiver $\Gamma_{A}$. Then there is a short exact sequence for a projective module $P$ and $k,\ell\in\mathbb{Z}^{+}$ as follows.
 \begin{align}\label{short-ex-seq}
 0\rightarrow(i,j)\xrightarrow{(\alpha_{\ell},-\beta_{k},\varepsilon)^{t}}(i,j+\ell)\oplus(i+k,j)\oplus P \xrightarrow{(\gamma_{k},\delta_{\ell},\pi)}  (i+k,j+\ell)\rightarrow0,
 \end{align}
where $\alpha_{\ell}$ and $\beta_{k}$ are the compositions of irreducible maps from $(0,i,j)$ to $(0,i,j+\ell)$ and $(0,i+k,j)$ respectively,  $\gamma_{k}$ and $\delta_{\ell}$ are the compositions of irreducible maps from  $(i,j+\ell)$ and $(i+k,j)$ to $(i+k,j+\ell)$ respectively.
\end{Lem}

\begin{proof}
Take a set $\{P_{1},P_{2},\cdots,P_{t}\}$  consisting of all indecomposable projective modules  in the rectangle area with vertices $(0,i,j)$, $(0,i,j+\ell)$, $(0,i+k,j)$ and $(0,i+k,j+\ell)$ on Euclidean component $\Gamma$. Then  the projective module $P=P_{1}\oplus P_{2}\oplus\cdots\oplus P_{t}$ in sequence {\rm(}\ref{short-ex-seq}{\rm)}.
The proof is done by induction on $k$ and $\ell$. Please refer to \cite[Lemma 2.3.1]{EK}.
\end{proof}

The proof of the following result is similar with \cite[Lemma 3.4]{CLZ}. We include a proof here. 
\begin{Lem}\label{com-irrmap}
Let $A$ be a self-injective algebra with an acyclic  stable generalized standard component $\Gamma_{0}$ and let  $(0,i,j)$ be  an indecomposable non-projective module on $\mathcal{C}$. If $(0,i,j)$ is a stable brick on $A$-$\stmod$, then the compositions  $\gamma_{k}\alpha_{\ell}$ of  $\alpha_{\ell}:(0,i,j)\to(0,i,j+\ell)$ and $\gamma_{k}:(0,i,j+\ell)\to(0,i+k,j+\ell)$ is non-zero  on $A$-$\stmod$ for any positive integer $\ell$ and $k$, where $\alpha_{\ell}$ {\rm(}resp. $\gamma_{k}${\rm)} is  the compositions of irreducible maps from $(0,i,j)$ to $(0,i,j+\ell)$ {\rm(}resp. from  $(0,i,j+\ell)$ to $(0,i+k,j+\ell)${\rm)}. 
\end{Lem}
\begin{Proof}
Consider the triangle induced by the  short exact sequence \eqref{short-ex-seq} as follows.
 \begin{align}\label{short-ex-seq-1}
(0,i,j)\xrightarrow{(\underline{\alpha_{\ell}},\underline{-\beta_{k}})^{t}}(0,i,j+\ell)\oplus(0,i+k,j) \xrightarrow{(\underline{\gamma_{k}},\underline{\delta_{\ell}})}  (0,i+k,j+\ell)\rightarrow(0,i,j)[1].
\end{align}
Applying $\StHom_A((0,i,j),-)$ to the triangle \eqref{short-ex-seq-1},  we have  a long exact sequence as follows. 
\begin{align*}
\cdots \to \underline{( (0,i,j), (0,i,j))} &\xrightarrow{(\underline{\alpha_{\ell}},\underline{\beta_{k}})^t_*} \underline{((0,i,j), (0,i,j+\ell)\oplus(0,i+k,j))} \\ &\quad \xrightarrow{(\underline{\gamma_{k}},\underline{\delta_{\ell}})_*} \underline{(0,i,j),  (0,i+k,j+\ell))} \to \cdots.
\end{align*}
Take the map $\alpha:=(\underline{\alpha_{\ell}},0)^t: (0,i,j)\to (0,i,j+\ell)\oplus(0,i+k,j)$.
Then we have \[(\underline{\gamma_{k}},\underline{\delta_{\ell}})_*\big(\alpha\big) = (\underline{\gamma_{k}},\underline{\delta_{\ell}})(\underline{\alpha_{\ell}},0)^t = \underline{\gamma_{k}\alpha_{\ell}}.\] 
On the contrary, we assume that $\underline{\gamma_{k}\alpha_{\ell}}=0$. It then follows from the exactness of the above long exact sequence that there is a morphism $\varphi\in\stend_A((0,i,j))$ such that $(\underline{\alpha_{\ell}},-\underline{\beta_{k}})^t\varphi= \alpha=(\underline{\alpha_{\ell}},0)^t$.
Since $\stend_A((0,i,j))\cong k$, $\varphi$ is an isomorphism.
This means that $\underline{\beta_{k}}=0$ on $A$-$\stmod$, which contradicts Lemma \ref{sectional-non-zero}.
\end{Proof}

\begin{Rem}
\begin{enumerate}[$(1)$]
\item If $X_{n}\xrightarrow{f_{n}} \cdots \xrightarrow{f_{3}}X_{2}\xrightarrow{f_{2}} X_{1}\xrightarrow{f_{1}} X_{0}=X$ be a sectional sequence of irreducible maps between indecomposable modules, then the composition $f_{n}\cdots f_{2}f_{1}$ is non-zero on $A$-$\m$. Please  refer to \cite[VII.2, Theorem 2.4]{ARS}.
\item 	There is a similar conclusion with Corollary \ref{com-irrmap} for a quasi-tube of a self-injective algebra. Please refer to \cite[Lemma 3.4]{CLZ}.
\end{enumerate}	
\end{Rem}

We visualize the morphism in Lemma \ref{com-irrmap} as follows.
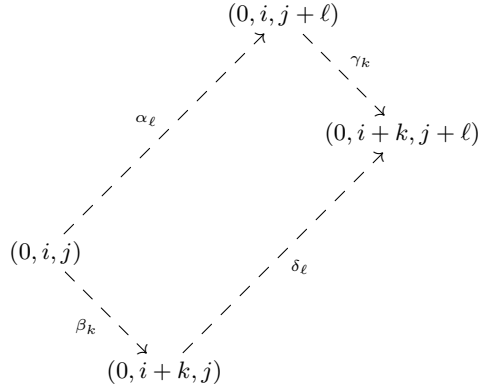
\[\small\xymatrix@R=45pt@C=45pt@!0{
	&  & & (0,i,j+\ell) \ar@{-->}[dr]^{\gamma_{k}} & \\
	& & & & (0,i+k,j+\ell)  \\
	&(0,i,j)  \ar@{-->}[dr]_{\beta_{k}} \ar@{-->}[uurr]^{\alpha_{\ell}} &  &  &  \\
 & & (0,i+k,j) \ar@{-->}[uurr]_{\delta_{\ell}} &    &  }
\]

\begin{Cor}\label{BGA-com-irrmap}
Let $A$ be a 2-domestic Brauer graph algebra and $X$ an indecomposable non-projective module on a connected component $\mathcal{C}$ of $\Gamma_{A}$. Then 
\begin{enumerate}[$(1)$]
\item If $\mathcal{C}$ is acyclic, then we may assume $X=(0,i,j)$ for integers $i$ and $j$. Then 
the compositions  $\gamma_{k}\alpha_{\ell}$ of  $\alpha_{\ell}:(0,i,j)\to(0,i,j+\ell)$ and $\gamma_{k}:(0,i,j+\ell)\to(0,i+k,j+\ell)$ is non-zero for any positive integer $\ell$ and $k$ on $A$-$\stmod$, where $\alpha_{\ell}$ {\rm (}resp. $\gamma_{k}${\rm)} is  the compositions of irreducible maps from $(0,i,j)$ to $(0,i,j+\ell)$ {\rm (}resp. from  $(0,i,j+\ell)$ to $(0,i+k,j+\ell)${\rm)}. 
\item If $\mathcal{C}$ admits a cycle, then we may assume $X=Q(0,i,j)$ for integers $i$ and $j$. Then the compositions  $\gamma_{k}\alpha_{\ell}$ of  $\alpha_{\ell}:Q(0,i,j)\to Q(0,i,j+\ell)$ and $\gamma_{k}:Q(0,i,j+\ell)\to Q(0,i+k,j+\ell)$ is non-zero for any positive integer $\ell$ and $k$ on $A$-$\stmod$, where $0<\ell<j$, and  $\alpha_{\ell}$ {\rm (}resp. $\gamma_{k}${\rm)} is  the compositions of irreducible maps from $Q(0,i,j)$ to $Q(0,i,j+\ell)$ {\rm (}resp. from  $Q(0,i,j+\ell)$ to $Q(0,i+k,j+\ell)${\rm)}. 
\end{enumerate}
\end{Cor}
\begin{proof}
The statement (1) follows from Corollary \ref{2-BGA-bricks} and the proof of Lemma \ref{com-irrmap}. The statement  (2) is a direct consequence of  \cite[Lemma 3.4]{CLZ}.
\end{proof}

\begin{Them}$($\cite[Theorem 4.17]{Z}$)$ \label{BGA-sms}
Let $A$ be a domestic Brauer graph algebra and  $\mathcal{S}$ an  orthogonal system on $A$-$\stmod$.  Then $\mathcal{S}$ is a simple-minded system if and only if $\mathcal{S}$ satisfies the following conditions.
\begin{enumerate}[$(1)$]
  \item $\mathcal{S}$ contains at least one non-periodic $A$-module.		
	\item $\Omega(\mathcal{S})$ {\rm(}or $\Omega^{-1}(\mathcal{S})${\rm)} is in  $\mathcal{F}(\mathcal{S})$.
\end{enumerate}
\end{Them}

\section{The stable bi-perpendicular category of  a stable brick}
Starting from this section, we always assume that $A$ is a 2-domestic Brauer graph algebra. Let $\Gamma_{0}$ and $\Gamma_{1}$ be the Euclidean components of the form $\mathbb{Z}\widetilde{A}_{p,q}$ and let $Q_{0}$ and $Q_{1}$ (resp. $Q'_{0}$ and $Q'_{1}$) be the quasi-tubes of rank $q$ (resp. $p$).
In this section, we determine the stable bi-perpendicular category of an object on an Euclidean component or on a quasi-tube of  rank not one. We shall take three steps to compute the stable bi-perpendicular category. 

We list some notations as follows. We use triple $(a,b,c)$ of integers to represent a  vertex for an Euclidean component.  Note that two Euclidean components $\Gamma_{0}$ and $\Gamma_{1}$  satisfy the condition $\Omega(_{s}\Gamma_{0})={_{s}\Gamma_{1}}$ and $\Omega(_{s}\Gamma_{1})={_{s}\Gamma_{0}}$. Without loss of generality, we use the triple $(0,b,c)$ to represent a vertex of $_{s}\Gamma_{0}$ (resp. the triple $(1,b,c)$ to represent a vertex of $_{s}\Gamma_{1}$) and we denote $\Omega((0,b,c))$ by $(1,b,c)$.   By the action of AR-translation $\tau$ on those triples,  $\tau(a,b,c)=(a,b-1,c-1)$.
From the shape of Euclidean components, we know that vertices $(a,b,c)$ and $(a,b-p\ell,c+q\ell)$ are equal on $A$-$\ind$ for every  integer $\ell\in\mathbb{Z}$. 
For a family of objects $\mathcal{X}$ on $A$-$\stmod$, we define 
\begin{equation}
\begin{aligned}
&\RS\mathcal{X}\colon=\{Z\in A{\text-}\stmod \mid \StHom_{A}(X,Z)\neq0, \exists  X\in\mathcal{X}\}, \nonumber\\  
&\LS\mathcal{X}\colon=\{Z\in A{\text-}\stmod\mid \StHom_{A}(Z,X)\neq0, \exists  X\in\mathcal{X}\}, \nonumber\\
&{^{\bot}\mathcal{X}^{\bot}}\colon=\{Z\in A{\text-}\stmod \mid \StHom_{A}(X,Z)=0, \StHom_{A}(Z,X)=0,\forall X\in\mathcal{X}\}.\nonumber
\end{aligned}
\end{equation}

We say $\RS\mathcal{X}$ is the {\bf right support of $\mathcal{X}$},  $\LS\mathcal{X}$ is the {\bf left support of $\mathcal{X}$} and ${^{\bot}\mathcal{X}^{\bot}}$ is the {\bf stable bi-perpendicular category of $\mathcal{X}$}. Note that ${^{\bot}\mathcal{X}^{\bot}}=A$-$\stmod\backslash(\RS\mathcal{X}\cup\LS\mathcal{X})$. 
We still use the notations $\RS X$, $\LS X$ and ${^{\bot}X^{\bot}}$ for an object $X$ on $A$-$\stmod$.

By Theorem \ref{simple-module-and-n-tube}, any object on a quasi-tube of rank one is not  in a simple-minded system. Thus we need not consider the stable bi-perpendicular category of the objects on the quasi-tubes of rank one. 

\begin{Lem}\label{Euclidean-Homogen}
Let $A$ be a 2-domestic Brauer graph algebra.  Then either $\dim_{k}\StHom_{A}(M,H)$ or $\dim_{k}\StHom_{A}(H,\\M)$ is non-zero for any non-periodic, non-projective module $M$ and any module $H$ in a homogeneous tube. 
\end{Lem}
\begin{proof}
On the contrary, we assume that there are a non-periodic, non-projective module $M_{0}$ in $\Gamma_{0}$ and a module $H_{0}$ in a homogeneous tube such that  $\StHom_{A}(M_{0},H_{0})\cong 0$ and $\StHom_{A}(H_{0},M_{0})\cong 0$. Since  $H_{0}$ is in a homogeneous tube, $\tau^{i}(H_{0})\cong H_{0}$ for every $i\in\mathbb{Z},$  \[\StHom_{A}(\tau^{i}(M_{0}),\tau^{i}(H_{0}))\cong\StHom_{A}(\tau^{i}(M_{0}),H_{0})\cong0, \  \StHom_{A}(\tau^{i}(H_{0}),\tau^{i}(M_{0}))\cong\StHom_{A}(H_{0},\tau^{i}(M_{0}))\cong0\] 
for every $i\in\mathbb{Z}.$ It is not hard to see that  $\StHom_{A}(M,H_{0})\cong 0$ and $\StHom_{A}(H_{0},M)\cong 0$ for every object $M$ in $\Gamma_{0}$. By Serre duality, $\StHom_{A}(M',H_{0})\cong 0$ and $\StHom_{A}(H_{0},M')\cong 0$ for every object $M'$ in $\Gamma_{1}$. Thus  $H_{0}$ is in the stable  bi-perpendicular category of every object on all Euclidean components. By Theorem \ref{simple-module-and-n-tube}, every simple-minded system $\mathcal{S}$ is contained in $\Gamma_{0}\cup\Gamma_{1}\cup Q_{0}\cup Q_{1}\cup Q'_{0}\cup Q'_{1}$.  By Proposition \ref{2-BGA-quasi-tube-bricks}, $H_{0}$ is in the stable  bi-perpendicular category of every object on all quasi-tubes of rank not one.  Hence $H_{0}\in{^{\bot}\mathcal{S}^{\bot}}$.  It is a contradiction.  Thus our conclusion is true. 
\end{proof}

\subsection{The intersection of  the stable bi-perpendicular category of an object (0,a,b) and an Euclidean component}
Take an object $(0,a,b)$ on $_{s}\Gamma_{0}$. By Lemma \ref{Euclidean-Homogen}, ${^{\perp}(0,a,b)^{\perp}}\cap{_{s}\mathcal{H}}$ is an empty set for every homogeneous tube $\mathcal{H}$. 
Now we  compute ${^{\perp}(0,a,b)^{\perp}}\cap{_{s}\Gamma_{0}}$, and then consider ${^{\perp}(0,a,b)^{\perp}}\cap{_{s}\Gamma_{1}}$ by Serre duality. We first give a definition as follows.
\begin{Def}\label{rectang-def}
Take  two objects $X=(i,m,n)$ and $Y=(i,j,k)$ on $\Gamma_{i}$ with $i=0,1$, $m\leq j$ and $n\geq k$  and let $\square_{X,Y}\colon=\{(i,e,f)\mid  m\leq e\leq j, k\leq f\leq n\}.$ $\square_{X,Y}$ is said to be {\bf the rectangle area of vertices $X$ and $Y$}. 
\end{Def}
Note that $\square_{X,Y}$ consists of the vertices which are contained in the rectangle area with vertices $X=(i,m,n)$, $(i,m,k)$, $(i,j,n)$ and $Y=(0,j,k)$. 

\medskip

\noindent{\bf Step one}: Compute ${^{\perp}(0,a,b)^{\perp}}\cap{_{s}\Gamma_{0}}$.

\medskip

By Corollary \ref{2-BGA-bricks} and Corollary \ref{BGA-com-irrmap}, the right (resp. left) support of vertex $(0,a,b)$ on $_{s}\Gamma_{0}$ states  as follows.	

\begin{equation}
\begin{aligned}\label{bi-perp-0}
\RS (0,a,b)\cap{_{s}\Gamma_{0}}=&\{(0,i,j)\mid  i\geq a, j\geqslant b\},\\
 \LS (0,a,b)\cap{_{s}\Gamma_{0}}=&\{(0,i,j)\mid  i\leq a, j\leq b\}.	
\end{aligned}
\end{equation}

\begin{equation}
\begin{aligned}\label{bi-perp-0'}
\RS (0,a-k,b-\ell)\cap{_{s}\Gamma_{0}}=&\{(0,i,j)\mid  i\geq a-k, j\geqslant b-\ell\},\\
 \LS  (0,a-k,b-\ell)\cap{_{s}\Gamma_{0}}=&\{(0,i,j)\mid  i\leq a-k, j\leq b-\ell\}.
\end{aligned}
\end{equation}
Thus the stable bi-perpendicular category of $(0,a,b)$ on $_{s}\Gamma_{0}$ states as follows. 
\begin{equation}
\begin{aligned}\label{bi-perp-1}
\mathcal{M}^{0}_{(0,a,b)}\colon=&{^{\perp}(0,a,b)^{\perp}}\cap{_{s}\Gamma_{0}}\\
=&{_{s}\Gamma_{0}}\backslash((\RS (0,a,b)\cup\LS (0,a,b))\cap{_{s}\Gamma_{0}})\\
=&\{(0,i,j)\mid  a-p<i<a, b<j<b+q\}\\
=&\square_{(0,a-p+1,b+q-1),(0,a-1,b+1)}.
\end{aligned}
\end{equation}

\begin{equation}
\begin{aligned}\label{bi-perp-2}
\mathcal{M}^{0}_{(0,a-k,b-\ell)}\colon=&{^{\perp}(0,a-k,b-\ell)^{\perp}}\cap{_{s}\Gamma_{0}}\\
=&{_{s}\Gamma_{0}}\backslash((\RS (0,a-k,b-\ell)\cup\LS (0,a-k,b-\ell))\cap{_{s}\Gamma_{0}}) \\
=&\{(0,i,j)\mid  a-k-p<i<a-k, b-\ell<j<b-\ell+q\}\\
=&\square_{(0,a-k-p+1,b-\ell+q-1),(0,a-k-1,b-\ell+1)}.
\end{aligned}
\end{equation}

Note that $\RS (0,a,b)=\RS (0,a-p,b+q)$ and $\LS (0,a,b)=\LS (0,a-p,b+q)$. The following diagram provides a visual depiction of $^{\bot}(0,a,b)^{\bot}\cap{_{s}\Gamma_{0}}$ on $_{s}\Gamma_{0}$ as follows. 
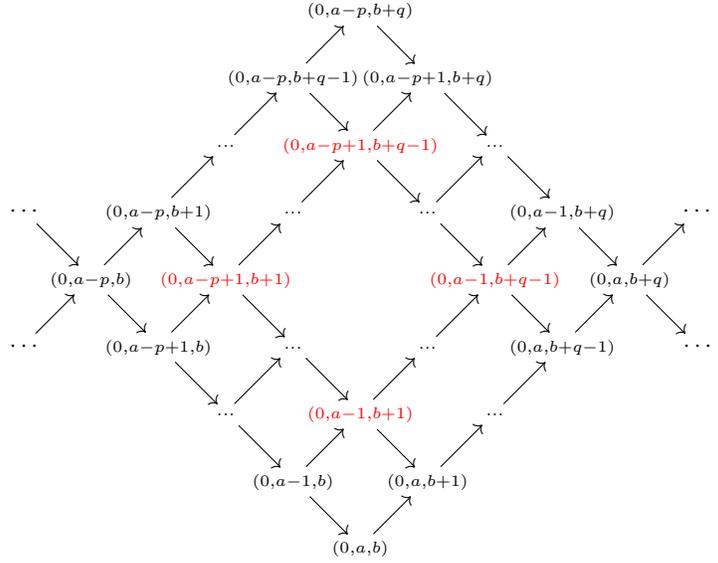
\vspace{-1.5cm}
\begin{small}
\[\xymatrix@dr@R=18pt@C=18pt@!0{
&&&&&&&&&&\cdots\\
&&\\
&&\scriptstyle(0,a-p,b+q)\ar[rr] &&\scriptstyle\scriptstyle(0,a-p+1,b+q) \ar[rr]& &\scriptstyle \cdots\scriptstyle \ar[rr] && \scriptstyle \scriptstyle(0,a-1,b+q)\ar[rr]&& \scriptstyle\scriptstyle(0,a,b+q) \ar[rr]\ar[uu]&&\ \  \cdots\,.& \\	
&& \\	
&&\scriptstyle(0,a-p,b+q-1)\ar[rr]\ar[uu] &&\scriptstyle{\color{red}(0,a-p+1,b+q-1)}\ar[rr] \ar[uu]&& 	\scriptstyle \cdots\ar[rr]\ar[uu] &&\scriptstyle\old{(0,a-1,b+q-1)}\ar[uu]\ar[rr] &&\scriptstyle\scriptstyle(0,a,b+q-1)\ar[uu]&&\\	
&&& \\	
&&	\scriptstyle \cdots \ar[uu] &&\scriptstyle\cdots \ar[uu]& & \scriptstyle &&\scriptstyle\cdots\ar[uu] && \scriptstyle \cdots\ar[uu]\\	
&&	& \\
&&\scriptstyle \scriptstyle(0,a-p,b+1) \ar[rr]\ar[uu]&&\ar[rr] \scriptstyle \old{(0,a-p+1,b+1)}\ar[uu] && \scriptstyle\cdots \ar[rr]\scriptstyle &&\scriptstyle\old{(0,a-1,b+1)}\ar[uu] \ar[rr]&&\scriptstyle(0,a,b+1)\ar[uu]\\
&&& \\
\cdots\ar[rr]&&\scriptstyle(0,a-p,b)\ar[rr]\ar[uu]  && \scriptstyle\scriptstyle(0,a-p+1,b)\ar[uu]\ar[rr] && \scriptstyle \cdots\ar[uu]\ar[rr]  &&\scriptstyle(0,a-1,b) \ar[rr]\ar[uu]&&\scriptstyle(0,a,b)\ar[uu] \\ 
&&\\
&&\cdots \ar[uu] 
}\] 
\end{small}
	
\medskip
	
\noindent{\bf Step two}: compute $^{\perp}(0,a-k,b-\ell)^{\perp}\cup{_{s}\Gamma_{1}}$ for $\ell\in\mathbb{Z}$.	

\medskip

By Serre duality, 
\[\StHom_{A}(X,Y)\cong\D\StHom_{A}(Y,\nu\Omega(X))\cong\D\StHom_{A}(Y,\Omega(X))\cong\D\StHom_{A}(\Omega^{-1}(Y),X).\]
Thus the right support of vertex $(0,a,b)$ on  $_{s}\Gamma_{1}$ is as follows.
\begin{equation}
\begin{aligned}\label{bi-perp-3}
\RS (0,a,b)\cap{_{s}\Gamma_{1}}&=\LS \Omega(0,a,b)\cap{_{s}\Gamma_{1}}\\
&=\LS (1,a,b)\cap{_{s}\Gamma_{1}}\\
&=\{(1,i,j)\mid  i\leq a, j\leq b\}. 
\end{aligned}
\end{equation}
Similarly, we have the left support of vertex $(0,a,b)$ on $_{s}\Gamma_{1}$ as follows.
\begin{equation}
	\begin{aligned}	\label{bi-perp-4}
		\LS (0,a,b)\cap{_{s}\Gamma_{1}}&=\RS \Omega^{-1}(0,a,b)\cap{_{s}\Gamma_{1}}\\
		&=\RS (1,a+1,b+1)\cap{_{s}\Gamma_{1}}\\
		&=\{(1,i,j)\mid  i\geq a+1, j\geq b+1\}.
	\end{aligned}
\end{equation}
In general, 
\begin{equation}
\begin{aligned}\label{bi-perp-5}
\RS(0,a-k,b-\ell)\cap{_{s}\Gamma_{1}}&=\LS\Omega(0,a-k,b-\ell)\cap{_{s}\Gamma_{1}}\\
&=\LS(0,a-k,b-\ell)\cap{_{s}\Gamma_{1}}\\
&=\{(1,i,j)\mid  i\leq a-k, j\leq b-\ell\}. 
\end{aligned}
\end{equation}
	
\begin{equation}
\begin{aligned}\label{bi-perp-6}
\LS (0,a-k,b-\ell)\cap{_{s}\Gamma_{1}}=&\RS \Omega^{-1}(0,a-k,b-\ell)\cap{_{s}\Gamma_{1}}\\
&=\RS (1,a-k+1,b-\ell+1)\cap{_{s}\Gamma_{1}}\\
&=\{(1,i,j)\mid  i\geq a-k+1, j\geq b-\ell+1\}.	
\end{aligned}
\end{equation}
Thus the stable bi-perpendicular category of $(0,a,b)$ on $_{s}\Gamma_{1}$ states as follows. 
By equations (\ref{bi-perp-3}) and (\ref{bi-perp-4}),
\begin{equation}
\begin{aligned}\label{bi-perp-7}
\mathcal{M}^{1}_{(0,a,b)}\colon=&{^{\perp}(0,a,b)^{\perp}}\cap{_{s}\Gamma_{1}}\\
=&{_{s}\Gamma_{1}}\backslash((\RS (0,a,b)\cup\LS (0,a,b))\cap{_{s}\Gamma_{1}})\\
=&\{(1,i,j)\mid  a-p+1\leq i\leq a, b+1\leq j\leq b+q\}\\
=&\square_{(1,a-p+1,b+q),(1,a,b+1)}.
\end{aligned}
\end{equation}
Note that  $\mathcal{M}^{1}_{(0,a,b)}={^{\perp}(1,a+1,b)^{\perp}}\cap{_{s}\Gamma_{1}}$.

Thus By equations (\ref{bi-perp-5}) and (\ref{bi-perp-6}),	
\begin{equation}
\begin{aligned}\label{bi-perp-8}
\mathcal{M}^{1}_{(0,a-k,b-\ell)}\colon&={^{\perp}(0,a-k,b-\ell)^{\perp}}\cap{_{s}\Gamma_{1}}\\
&={_{s}\Gamma_{1}}\backslash((\RS (0,a-k,b-\ell)\cup\LS (0,a-k,b-\ell))\cap{_{s}\Gamma_{1}})\\
&=\{(1,i,j)\mid  a-k-p+1\leq i\leq a-k, b-\ell+1\leq j\leq b-\ell+q\}\\
&=\square_{(1,a-k-p+1,b-\ell+q),(1,a-k,b-\ell+1)}.
\end{aligned}
\end{equation}
The  following diagram shows us the position of $\mathcal{M}^{1}_{(0,a,b)}$ (the rectangle area with red vertices) on $\Gamma_{1}$.

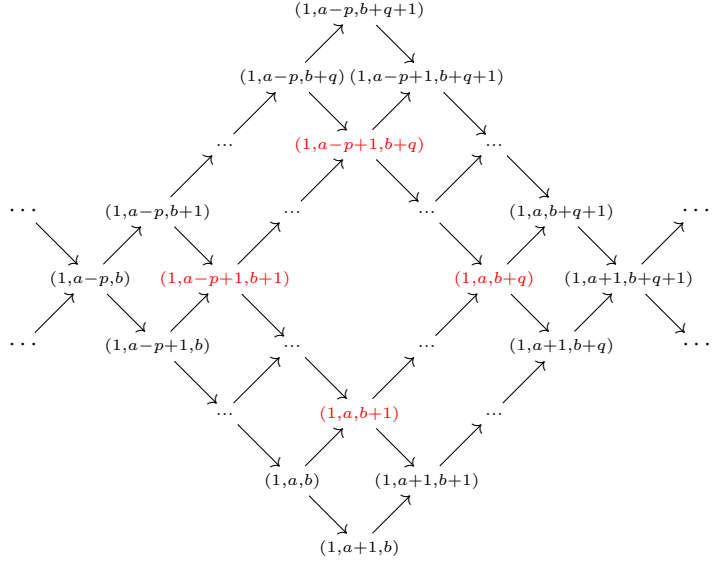
\vspace{-1.7cm}
\begin{small}
\[\xymatrix@dr@R=18pt@C=18pt@!0{	
&&&&&&&&&&\cdots&&\\
&&&\\
&&\scriptstyle(1,a-p,b+q+1)\ar[rr] &&\scriptstyle\scriptstyle(1,a-p+1,b+q+1) \ar[rr]& &\scriptstyle \cdots\scriptstyle \ar[rr] && \scriptstyle \scriptstyle(1,a,b+q+1)\ar[rr]&& \scriptstyle\scriptstyle(1,a+1,b+q+1) \ar[rr]\ar[uu]&&\ \  \cdots\,.& \\	
&&\\	
&&\scriptstyle(1,a-p,b+q)\ar[rr]\ar[uu] &&\scriptstyle\old{(1,a-p+1,b+q)}\ar[rr] \ar[uu]&& 	\scriptstyle \cdots\ar[rr]\ar[uu] &&\scriptstyle\old{(1,a,b+q)}\ar[uu]\ar[rr] &&\scriptstyle\scriptstyle(1,a+1,b+q)\ar[uu]&&\\	
&&&\\	
&&	\scriptstyle \cdots \ar[uu] &&\scriptstyle\cdots \ar[uu]& & \scriptstyle &&\scriptstyle\cdots\ar[uu] && \scriptstyle \cdots\ar[uu]\\	
&&\\
&&\scriptstyle \scriptstyle(1,a-p,b+1) \ar[rr]\ar[uu]&&\ar[rr] \scriptstyle \old{(1,a-p+1,b+1)}\ar[uu] && \scriptstyle\cdots \ar[rr]\scriptstyle &&\scriptstyle\old{(1,a,b+1)}\ar[uu] \ar[rr]&&\scriptstyle(1,a+1,b+1)\ar[uu]\\
&&& \\
\cdots\ar[rr]&&\scriptstyle(1,a-p,b)\ar[rr]\ar[uu]  && \scriptstyle\scriptstyle(1,a-p+1,b)\ar[uu]\ar[rr] && \scriptstyle \cdots\ar[uu]\ar[rr]  &&\scriptstyle(1,a,b) \ar[rr]\ar[uu]&&\scriptstyle(1,a+1,b)\ar[uu] \\ 
&&&&&&&&&&\\
&&\cdots \ar[uu]   
}\] 
\end{small}

\subsection{The intersection of the stable bi-perpendicular category of $(0,a,b)$  and a quasi-tube}
In this subsection, we determine the intersection of ${^{\perp}(0,a-k,b-\ell)^{\perp}}$ and quasi-tube of  rank not one. We need  some preparations as follows.

\begin{Prop}{\rm(\cite[Corollary 3.5]{K})} \label{quasi-simple-1}Let $A$ be a self-injective artin algebra, and let  $0\xrightarrow{}X\xrightarrow{}Y\xrightarrow{f}Z\xrightarrow{} 0$ be a short exact sequence on $A$-$\m$,  which is not almost split. Suppose that one of the modules  $Y$ and $Z$ is indecomposable and $f$ is irreducible. Then the following conditions are equivalent for a natural number $n>1:$
\begin{enumerate}[$(1)$]
\item $\alpha(\tau^{-1}(X))=n.$	
\item $X$ is simple and $\rad P/\soc P$ decomposes into $n$ summands, where $P$ denotes the projective cover of $X$.
\item $X$ is simple and $\rad I/\soc I$ decomposes into $n$ summands, where $I$ denotes the injective envelope  of $X$.
\end{enumerate}
Moreover, if $X$ satisfies these conditions, then $Y$ is either indecomposable projective or isomorphic to the radical of some indecomposable projective module.
\end{Prop}	

Let $A$ be an algebra and $\mathcal{C}$ a component of  AR-quiver of $A$. A {\bf coray point} of $\mathcal{C}$ is defined to be a point $X$ of $\mathcal{C}$ such that there exists an infinite sectional path in $\mathcal{C}$ 
\[\cdots \rightarrow [r+1]X \rightarrow [r]X \rightarrow \cdots \rightarrow [2]X \rightarrow [1]X=X.
\] ending with $X$ and containing all sectional paths ending with $X.$ The corresponding $A$-module $X$ is called the {\bf coray module.}

\begin{Prop}{\rm(\cite[Proposition 4.12]{Z})} \label{sms-quasi-simple}
Let $A$ be representation-infinite self-injective algebra. Given a triangle 	
\begin{align}\label{tri-1}
	X\xrightarrow{}Y\xrightarrow{\underline{f}}Z\xrightarrow{} X[1]
\end{align}
in $A$-$\stmod$, where both $Y$ and $Z$ are indecomposable, non-projective, and $f:Y\rightarrow Z$ is an irreducible map on $A$-$\m$. Then  $X$ {\rm(}resp. $X[1]${\rm)} is a coray point of some stable connected AR-component. 
\end{Prop}	
\begin{proof}Note that a point is a ray point if and only it is a coray point in the stable	 AR-quiver of a self-injective algebra. Let $X$ be a point in the stable connected component $_{s}\Gamma_{1}$ of stable AR-quiver of $A$. We only  prove that $X$  is a  coray point of $_{s}\Gamma_{1}$, since $[1]: A$-$\stmod\rightarrow A$-$\stmod$ is an equivalence.  
We assume that the triangle  (\ref{tri-1}) is induced by a short exact sequence as follows.
\begin{align}\label{ses-1}
	0\xrightarrow{}X\xrightarrow{}Y\oplus P\xrightarrow{(f,h)}Z\xrightarrow{} 0.
\end{align}
Where $P$ is a projective $A$-module. If the sequence   (\ref{ses-1}) is  an almost split sequence, then $Y$ is  isomorphic to $\rad(P)/\soc(P)$ and $P$ is indecomposable. Since $Y$ is indecomposable,  $\alpha(\tau^{-1}(X))=2$. Thus $X$ is a coray point on $_{s}\Gamma_{1}$. 
	
Now we assume that the short exact sequence  (\ref{ses-1}) is not an almost split sequence.  
	
\noindent{\bf Case one:} $P$ is isomorphic to zero.
	
We shall show that $X$ is a coray point on $_{s}\Gamma_{1}$. On the contrary, we assume that $X$ is not  a  coray point of $_{s}\Gamma_{1}$.  Take  $\alpha(\tau^{-1}(X))=n>1$.   By Proposition \ref{quasi-simple-1}, $Y$ is either indecomposable projective or isomorphic to the radical of some indecomposable projective module. The former case contradicts the condition that $Y$ is non-projective. For the latter case, $Y$ is isomorphic to the radical $\rad P$ of an indecomposable projective $P$. By the equivanent conditions of Proposition \ref{quasi-simple-1},  $X$ is  a simple module and $X\cong\soc P$. Therefore $Z$ is isomorphic to $\rad P/\soc P$ and it is decomposable (with $n$ direct summands). This contradicts the condition that $Z$ is indecomposable. Hence $\alpha(\tau^{-1}(X))=1$ and then $X$  is a  coray point of $_{s}\Gamma_{1}$.
	
\noindent{\bf Case two}: $P$ is not isomorphic to zero.

Since $f$ is an irreducible map, $f$ is a monomorphism or an epimorphism.  
If $f$ is a monomorphism,  then we take a short exact sequence as follows.
\begin{align}\label{ses-3}
	0\xrightarrow{}Y\xrightarrow{f}Z\xrightarrow{} \Coker(f) \xrightarrow{}0.
\end{align}
It can be proved dually that $\Coker(f)$ is a ray point on $_{s}\Gamma_{1}$. Considering the triangle induced by short exact sequence (\ref{ses-3}), we have $X[1]\cong\Coker(f)$ on $A$-$\stmod$. Thus $X$ is a coray  point on $\_{s}\Gamma_{1}$.
	
If f is an epimorphism, then we take a short exact sequence as follows.
\begin{align}\label{ses-4}
0\xrightarrow{}\Ker(f)\xrightarrow{}Y\xrightarrow{f}Z\xrightarrow{} 0.
\end{align}
It is similar to the Case one that $\Ker(f)$ is a coray  point on $_{s}\Gamma_{1}$. Considering the triangle induced by short exact sequence (\ref{ses-4}), we have $X\cong \Ker(f)$ on $A$-$\stmod$. Thus $X$ is a coray point on $_{s}\Gamma_{1}$.
\end{proof}
	
By the shape  of Euclidean component, there are precise two irreducible maps ended at  vertex $(0,a,b)$ on ${_{s}\Gamma_{0}}$. 
\[(0,a,b-1)\rightarrow(0,a,b),\ \  (0,a-1,b)\rightarrow(0,a,b).\]
Extending them to triangles as follows.
\[(0,a,b-1)\xrightarrow{}(0,a,b)\xrightarrow{}Z_{0}\xrightarrow{} (0,a,b-1)[1],\]
\[(0,a-1,b)\xrightarrow{}(0,a,b)\xrightarrow{}Z_{1}\xrightarrow{} (0,a-1,b+1)[1].\]
Dually, there are precise two irreducible maps started at  vertex $(0,a,b)$	in ${_{s}\Gamma_{0}}$. 
\[(0,a,b)\rightarrow(0,a,b+1),\ \  (0,a,b)\rightarrow(0,a+1,b).\]
Extending them to triangles as follows.
\[Z'_{0}\xrightarrow{}(0,a,b)\xrightarrow{}(0,a,b+1)\xrightarrow{} Z'_{0}[1],\]	
\[Z'_{1}\xrightarrow{}(0,a,b)\xrightarrow{}(0,a+1,b)\xrightarrow{}Z'_{1}[1].\]	
By Proposition \ref{sms-quasi-simple}, both $Z_{0}$ and $Z_{1}$ (resp. $Z'_{0}$ and $Z'_{1}$) are quasi-simples. Note that one of them is on a quasi-tube of rank $p$ and the other one is on a quasi-tube of rank $q$. Without loss of generality, we assume $Z'_{0}$ (resp. $Z'_{1}$) is contained on the quasi-tube $Q_{0}$ (resp. $Q'_{0}$) of rank $q$ (resp. $p$). Note that $Z_{0}=\Omega(Z'_{0})$ and $Z_{1}=\Omega(Z'_{1})$.
We denote $Z'_{0}$ (resp. $Z'_{1}$) by $Q(0,0,0)$ (resp. $Q(0,0,0)'$). Thus $Z_{0}=\Omega(Z'_{0})=Q(1,0,0)$ and $Z_{1}=\Omega(Z'_{1})=Q(1,0,0)'$. Take  $\tau Q(0,i,0)=Q(0,i-1,0)$ and $\tau Q(0,i,0)'=Q(0,i-1,0)'$ for any integer $i$. 

Thus every non-projective vertex on a quasi-tube is labeled by $Q(i,j,k)$ or  $Q(i,j,k)'$, where $i=0,1$ and $j,k\in\mathbb{Z}$.
Note that $Q(i,j,k)=Q(i,j,k+q\ell)$ and $Q(i,j,k)'=Q(i,j,k+p\ell)'$ for $\ell\in\mathbb{Z}$.
\begin{Lem}{\rm(\cite[Lemma 3.8]{Z0})}\label{irreducible-quasi-simples}
Let $A$ be a  self-injective algebra and let $M$, $N$ be two indecomposable, non-projective modules. Let $M_{n}\xrightarrow{f_{n}} \cdots \xrightarrow{f_{3}}M_{2}\xrightarrow{f_{2}} M_{1}\xrightarrow{f_{1}} M_{0}=M$ be a sectional path ending at $M$, and define $S_{i}$  by extending   $M_{i}\xrightarrow{f_{i}} M_{i-1}$ to a triangle 
\[M_{i}\xrightarrow{f_{i}}M_{i-1} \xrightarrow{}  S_{i}\rightarrow M_{i}[1].\]
Then we have $\tau(S_{i-1})=S_{i}$ for $i=2,\cdots, n$. 
\end{Lem}

\begin{Prop}{\rm(}\cite[I. Section 5, Corollary 5.7]{ARS}{\rm)}\label{seq-two-isos}
Let $A$ be an artin algebra.  Take a short exact sequence on $A$-$\m$ as follows. 
\[0\xrightarrow{}X\xrightarrow{(f_{1},f_{2})^{t}}Y_{1}\oplus Y_{2} \xrightarrow{(g_{1},g_{2})}  Z\rightarrow 0.\]
Then the following conclusion holds.     
\begin{enumerate}[$(1)$]
\item \[\xymatrix{X \ar[r]^-{f_{1}}\ar[d]_-{-f_{2}} & Y_{1}\ar[d]^-{g_{1}} \\
Y_{2} \ar[r]_-{g_{2}}& Z}\] 
is both  a pushout and a pullback diagram. 
\item $\Coker(f_{1})\cong\Coker(g_{2})$, $\Coker(f_{2})\cong\Coker(g_{1})$. 
\end{enumerate}
\end{Prop}

\begin{Prop}{\rm(\cite[Proposition 2.2]{Ap})}\label{hom-quasi-sim=2}
Let  $A$ be a symmetric special biserial algebra, $M$ a string module on $A$-$\m$ and let $\Phi$ be the set of quasi-simples which are not band modules. Then 
\[\sum_{X\in\Phi}\dk\StHom_A(M, X)=\sum_{X\in\Phi}\dk\StHom_A(X,M)=2.\] 	
\end{Prop}
	
\begin{Def}
Let $A$ be a self-injective algebra and $X=(j,k)$ a $\tau$-periodic $A$-module. The {\bf wing} of $X$ is the set of isoclasses of indecomposable non-projective modules on a quasi-tube given by
\begin{align*}
\mathcal{W}_{X} &:= \{(m,\ell) \mid m\leq j,\;\; j+k\leq m+\ell\}. 
\end{align*}
\end{Def}
	
According to the equation of dimension of stable Hom-space provided by Erdmann and Kerner \cite[2.2]{EK}, 
\begin{equation}
\begin{aligned}\label{eq-0}	
\RS(0,a,b)\cap({_{s}Q_{1}}\cup{_{s}Q'_{1}})&=\mathcal{W}_{Q(1,0,0)}\cup\mathcal{W}_{Q(1,0,0)'},\\
\LS(0,a,b)\cap({_{s}Q_{1}}\cup{_{s}Q'_{1}})&=\{0\}.
\end{aligned}	
\end{equation}
\begin{equation}
 \begin{aligned}\label{eq-0'}
\LS(0,a,b)\cap({_{s}Q_{0}}\cup{_{s}Q'_{0}})&=\mathcal{W}_{Q(0,0,0)}\cup\mathcal{W}_{Q(0,0,0)'},\\
\RS(0,a,b)\cap({_{s}Q_{0}}\cup{_{s}Q'_{0}})&=\{0\}.
\end{aligned}	
\end{equation}
Where $\mathcal{W}_{Q(i,j,0)}$ is the wing of quasi-simple $Q(i,j,0)$ for $i=0,1$ and $j\in\mathbb{Z}.$ See also \cite[Lemma 5.12]{CLZ}.
	
\begin{equation}
\begin{aligned}\label{eq-1}
\RS (0,a-\ell,b)\cap{_{s}Q_{1}}=\mathcal{W}_{Q(1,0,0)},\ \ 
\RS(0,a-\ell,b)\cap{_{s}Q'_{1}}=\mathcal{W}_{\tau^{\ell}Q(1,0,0)'}.
\end{aligned}	
\end{equation}
Thus
\begin{align}\label{eq-2}
\RS(0,a-\ell,b)\cap({_{s}Q_{1}}\cup{_{s}Q'_{1}})=\mathcal{W}_{Q(1,0,0)}\cup\mathcal{W}_{\tau^{\ell}Q(1,0,0)'}.  
\end{align}
\begin{equation}
	\begin{aligned}\label{eq-3}
\RS(0,a,b-\ell)\cap{_{s}Q_{1}}=\mathcal{W}_{\tau^{\ell}Q(1,0,0)},\  \ 
\RS(0,a,b-\ell)\cap{_{s}Q'_{1}}=\mathcal{W}_{Q(1,0,0)'}.
\end{aligned}	
\end{equation}

\begin{align}\label{eq-4}
\RS(0,a,b-\ell)\cap({_{s}Q_{1}}\cup{_{s}Q'_{1}})=\mathcal{W}_{\tau^{\ell}Q(1,0,0)}\cup\mathcal{W}_{Q(1,0,0)'}. 
\end{align}
Thus 
\begin{equation}
\begin{aligned}\label{eq-5}	
\RS(0,a-k,b-\ell)\cap({_{s}Q_{1}}\cup{_{s}Q'_{1}})&=\mathcal{W}_{\tau^{\ell}Q(1,0,0)}\cup\mathcal{W}_{\tau^{k}Q(1,0,0)'},\\
\LS(0,a-k,b-\ell)\cap({_{s}Q_{1}}\cup{_{s}Q'_{1}})&=\{0\}.
\end{aligned}  
\end{equation}
Similarly, we have description of the left support as follows. 
\begin{equation}
\begin{aligned}\label{eq-6}	
\LS (0,a-\ell,b)\cap{_{s}Q_{0}}=\mathcal{W}_{Q(0,0,0)}, \ \
\LS(0,a-\ell,b)\cap{_{s}Q'_{0}}=\mathcal{W}_{\tau^{\ell}Q(0,0,0)'}.
\end{aligned}	
\end{equation}
Thus
\begin{align}\label{eq-7}	
\LS(0,a-\ell,b)\cap({_{s}Q_{0}}\cup{_{s}Q'_{0}})=\mathcal{W}_{Q(0,0,0)}\cup\mathcal{W}_{\tau^{\ell}Q(0,0,0)'}.
\end{align}		
\begin{equation}
	\begin{aligned}\label{eq-8}	
\LS(0,a,b-\ell)\cap{_{s}Q_{0}}=\mathcal{W}_{\tau^{\ell}Q(0,0,0)},\ \
\LS (0,a,b+\ell)\cap{_{s}Q'_{0}}=\mathcal{W}_{Q(0,0,0)'},
\end{aligned}	
\end{equation}
Thus
\begin{align}\label{eq-9}	
\LS(0,a,b-\ell)\cap({_{s}Q_{0}}\cup{_{s}Q'_{0}})=\mathcal{W}_{\tau^{\ell}Q(0,0,0)}\cup\mathcal{W}_{Q(0,0,0)'}.
\end{align}
Thus 
\begin{equation}
\begin{aligned}\label{eq-10}	
\LS(0,a-k,b-\ell)\cap({_{s}Q_{0}}\cup{_{s}Q'_{0}})&=\mathcal{W}_{\tau^{\ell}Q(0,0,0)}\cup\mathcal{W}_{\tau^{k}Q(0,0,0)'},\\
\RS(0,a-k,b-\ell)\cap({_{s}Q_{0}}\cup{_{s}Q'_{0}})&=\{0\}.
\end{aligned}
\end{equation}
Thus by equation (\ref{eq-0}), 	
\begin{equation}
\begin{aligned}\label{eq-13}	
\mathcal{Q}^{1}_{(0,a,b)}\colon&={^{\perp}(0,a,b)^{\perp}}\cap({_{s}Q_{1}}\cup{_{s}Q'_{1}})\\
&=({_{s}Q_{1}}\cup{_{s}Q'_{1}})\backslash((\RS (0,a,b)\cup\LS (0,a,b))\cap({_{s}Q_{1}}\cup{_{s}Q'_{1}}))\\
&=\{Q(1,i,j)\mid 1\leq i\leq q-1, 0\leq j\leq q-2, i+j\leq q-1\}\\
&\cup\{Q(1,i,j)'\mid 1\leq i\leq p-1,0\leq j\leq p-2, i+j\leq p-1\}.
\end{aligned}
\end{equation}
Thus by equation (\ref{eq-0'}), 
\begin{equation}
\begin{aligned}\label{eq-11}	
\mathcal{Q}^{0}_{(0,a,b)}\colon&={^{\perp}(0,a,b)^{\perp}}\cap({_{s}Q_{0}}\cup{_{s}Q'_{0}})\\
&=({_{s}Q_{0}}\cup{_{s}Q'_{0}})\backslash((\RS (0,a,b)\cup\LS (0,a,b))\cap({_{s}Q_{0}}\cup{_{s}Q'_{0}}))\\
&=\{Q(0,i,j)\mid  1\leq i\leq q-1, 0\leq j\leq q-2, i+j\leq q-1\}\\
&\cup\{Q(0,i,j)'\mid  1\leq i\leq p-1, 0\leq j\leq p-2, i+j\leq p-1\}.
\end{aligned}
\end{equation}
Thus by equation (\ref{eq-5}), 	
\begin{equation}
\begin{aligned}\label{eq-14}	
\mathcal{Q}^{1}_{(0,a-k,b-\ell)}\colon=&^{\perp}(0,a-k,b-\ell)^{\perp}\cap({_{s}Q_{1}}\cup{_{s}Q'_{1}})\\
=&({_{s}Q_{1}}\cup{_{s}Q'_{1}})\backslash((\RS (0,a-k,b-\ell)\cup\LS (0,a-k,b-\ell))\cap({_{s}Q_{1}}\cup{_{s}Q'_{1}}))\\
=&\{Q(1,i,j)\mid  -\overline{\ell+1}\leq i\leq-\overline{\ell+q-1}, 0\leq j\leq q-2, i+j\leq q-1\}\\
&\cup\{Q(1,i,j)'\mid  -\overline{k+1}\leq i\leq -\overline{k+p-1}, 0\leq j\leq p-2, i+j\leq p-1\}.
\end{aligned}
\end{equation}
Thus by equation (\ref{eq-10}), 	
\begin{equation}
\begin{aligned}\label{eq-12}	
\mathcal{Q}^{0}_{(0,a-k,b-\ell)}\colon&={^{\perp}(0,a-k,b-\ell)^{\perp}}\cap({_{s}Q_{0}}\cup{_{s}Q'_{0}})\\
&=({_{s}Q_{0}}\cup{_{s}Q'_{0}})\backslash((\RS (0,a-k,b-\ell)\cup\LS(0,a-k,b-\ell))\cap({_{s}Q_{0}}\cup{_{s}Q'_{0}}))\\
&=\{Q(0,i,j)\mid -\overline{\ell+1}\leq i\leq-\overline{\ell+q-1}, 0\leq j\leq q-2, i+j\leq q-1\}\\
&\cup\{Q(0,i,j)'\mid-\overline{k+1}\leq i\leq -\overline{k+p-1}, 0\leq j\leq p-2, i+j\leq p-1\}.
\end{aligned}
\end{equation}
$^{\bot}(0,a,b)^{\bot}\cap{_{s}Q_{0}}$ is depicted as the triangle area with red  vertices in the following diagram.

\begin{small}
\[\xymatrix@dr@R=16pt@C=16pt@!0{\scriptstyle Q(0,0,p)\ar[rr] &&\scriptstyle Q(0,1,p-1) \ar[rr]& &\scriptstyle \cdots \scriptstyle \ar[rr] && \scriptstyle \cdots\ar[rr]&& \scriptstyle Q(0,p-1,1) \ar[rr]&&\scriptstyle\  Q(0,p,0)\,.& \\	
&& \\	
\scriptstyle Q(0,0,p-1) \ar[uu]\ar[rr]&&\scriptstyle \old{Q(0,1,p-2)}\ar[uu]\ar[rr] &&\scriptstyle \cdots \ar[uu]\ar[rr] &&\scriptstyle Q(0,p-2,1)\ar[uu]\ar[rr]&&\scriptstyle \old{Q(0,p-1,0)}\ar[uu]\\	
&& \\	
\scriptstyle\cdots \ar[uu]\ar[rr] &&\scriptstyle\cdots\ar[uu] & & \cdots  && \\	
&&& \\
\scriptstyle \cdots\ar[rr]\ar[uu]  && \scriptstyle Q(0,1,1) \ar[uu]&&  \\
&&& \\
\scriptstyle Q(0,0,1) \ar[uu]\ar[rr] && \scriptstyle \old{Q(0,1,0)}\ar[uu] \\	
&& \\
\scriptstyle Q(0,0,0) \ar[uu] \\}\] 
\end{small}

Finally, 
\begin{align} 
^{\perp}(0,a,b)^{\perp}\cap{_{s}\Gamma_{A}}&=\mathcal{M}^{0}_{(0,a,b)}\cup\mathcal{M}^{1}_{(0,a,b)}\cup\mathcal{Q}^{0}_{(0,a,b)}\cup\mathcal{Q}^{1}_{(0,a,b)},\\	
^{\perp}(0,a-k,b-\ell)^{\perp}\cap{_{s}\Gamma_{A}}&=\mathcal{M}^{0}_{(0,a-k,b-\ell)}\cup\mathcal{M}^{1}_{(0,a-k,b-\ell)}\cup\mathcal{Q}^{0}_{(0,a-k,b-\ell)}\cup\mathcal{Q}^{1}_{(0,a-k,b-\ell)}.	
\end{align}

\subsection{The stable bi-perpendicular category of an object on a  quasi-tube of rank not 1}
By Theorem \ref{simple-module-and-n-tube}, 
Take a stable brick $Q(0,c,d)$ on the quasi-tube $_{s}Q_{0}$. By Theorem \ref{simple-module-and-n-tube}, if $d\geq q-1$, then $Q(0,c,d)$ is not contained in a simple-minded system. 
Without loss of generality, we assume the non-negative integer $d\leq q-2$.
In this subsection, we shall determine the intersection of ${^{\bot}Q(0,c,d)^{\bot}}$ and quasi-tubes of rank not $1$ or Euclidean components.

\medskip
\noindent{\bf Step one}: compute ${^{\bot}Q(0,c,d)^{\bot}}\cap {_{s}Q_{0}}$.	

\medskip

\begin{equation}
\begin{aligned}\label{equ-0}
\RS Q(0,c,d)\cap {_{s}Q_{0}}&=\{Q(0,i,j)\mid c\leq i\leq c+d, c+d\leq i+j\},\\
\LS Q(0,c,d)\cap {_{s}Q_{0}}&=\{Q(0,i,j)\mid i\leq c, c\leq i+j\leq c+d\}.
\end{aligned}
\end{equation}

\begin{equation}
\begin{aligned}\label{equ-1}
\mathcal{N}^{0}_{Q(0,c,d)}\colon=&{^{\bot}Q(0,c,d)^{\bot}}\cap {_{s}Q_{0}}\\
=&{_{s}Q_{0}}\backslash((\RS Q(0,c,d)\cup\LS Q(0,c,d))\cap{_{s}Q_{0}})\\
=&\{Q(0,i,j)\mid c+1\leq i\leq c+d-1, i+j\leq c+d-1\}\\
&\cup\{Q(0,i,j)\mid c+d-q+1\leq i\leq c-1, i+j\leq c-1\}\\
&\cup\{Q(0,i,j)\mid c+d-q+1\leq i\leq c-1,c+d+1\leq i+j\leq c+q-1\}.
\end{aligned}
\end{equation}
$\mathcal{N}^{0}_{Q(0,c,d)}$ is a set of  vertices which consists of  
\begin{enumerate}[$(1)$]
	\item The rectangle with vertices $Q(1,c+d-q+1,q-d)$, $Q(0,c+d-q+1,2q-d-2)$, $Q(0,c-1,q)$ and $Q(0,c-1,d+2)$.
	\item The triangle with vertices $Q(0,c+d-q+1,0)$, $Q(0,c+d-q+1,q-d-2)$ and $Q(0,c-1,0)$.
	\item The triangle with vertices $Q(0,c+1,0)$, $Q(0,c+1,d-2)$ and $Q(0,c+d-1,0)$. 
\end{enumerate}

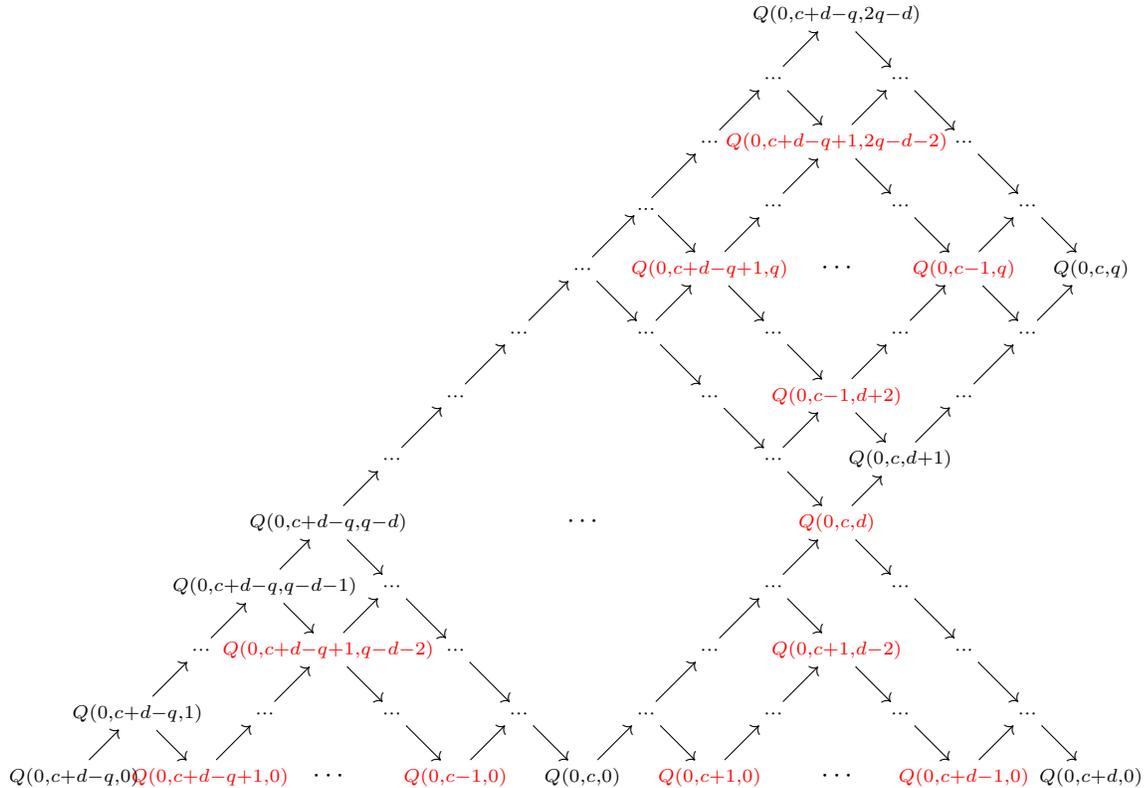
\[\xymatrix@dr@R=17pt@C=17pt@!0{	
	\scriptstyle Q(0,c+d-q,2q-d) \ar[rr]&&\scriptstyle \cdots\ar[rr]&&\scriptstyle \cdots\ar[rr]&&\scriptstyle \cdots\ar[rr]&&\scriptstyle Q(0,c,q)&&&&\\	
	&&&&&\\			
	\scriptstyle \cdots\ar[uu]\ar[rr]	&&\scriptstyle \scriptstyle \old{Q(0,c+d-q+1,2q-d-2)}\ar[uu]\ar[rr]&&\scriptstyle \cdots\ar[rr]&&\scriptstyle\scriptstyle\old{ Q(0,c-1,q)}\ar[rr]\ar[uu]&& \scriptstyle\cdots\ar[uu] &&\\	
	&&&&\\			
	\scriptstyle \cdots\ar[uu]	&&\scriptstyle \cdots\ar[uu]&&\cdots&&\scriptstyle \cdots\ar[uu]&& \scriptstyle\cdots \ar[uu] &&\\	
	&&&&\\		
	\scriptstyle\dots\ar[rr]\ar[uu] &&\scriptstyle \old{Q(0,c+d-q+1,q)}\ar[uu] \ar[rr]& &\scriptstyle \cdots \scriptstyle \ar[rr] && \scriptstyle \old{Q(0,c-1,d+2)}\ar[rr]\ar[uu]&& \scriptstyle Q(0,c,d+1)\ar[uu] &&\\	
	&&&& \\	
	\scriptstyle\cdots\ar[uu]\ar[rr]&&\scriptstyle \cdots\ar[uu]\ar[rr] &&\scriptstyle \cdots\ar[rr] &&\scriptstyle\cdots\ar[uu]\ar[rr]&&\scriptstyle \old{Q(0,c,d)}\ar[uu]\ar[rr]&&\scriptstyle\cdots\ar[rr]&&\scriptstyle\cdots\ar[rr]&&\scriptstyle\cdots\ar[rr]&&\scriptstyle Q(0,c+d,0)\\	
	&&&&&\\	
	\scriptstyle\cdots \ar[uu]&& &&&&&&\scriptstyle\cdots\ar[uu]\ar[rr]&&\scriptstyle \old{Q(0,c+1,d-2)}\ar[rr]&&\scriptstyle\cdots\ar[rr]&&\scriptstyle \old{Q(0,c+d-1,0)}\ar[uu]\\	
	&&&& \\
	\scriptstyle \cdots\ar[uu]  && &&\cdots&&&&\scriptstyle \cdots\ar[uu]&&\scriptstyle\cdots\ar[uu]&&\cdots&&  \\
	&&&&\\
	\scriptstyle \cdots \ar[uu]&& \scriptstyle &&&&&&\scriptstyle\cdots\ar[rr]\ar[uu]&&\scriptstyle \old{Q(0,c+1,0)}\ar[uu] \\	
	&&&&&\\
	\scriptstyle Q(0,c+d-q,q-d) \ar[uu] \ar[rr]&&\scriptstyle \cdots\ar[rr]&&\scriptstyle\cdots\ar[rr]&&\scriptstyle \cdots\ar[rr]&&\scriptstyle Q(0,c,0)\ar[uu]\\
	&&&&\\
	\scriptstyle Q(0,c+d-q,q-d-1)\ar[rr]\ar[uu]&&\scriptstyle \old{Q(0,c+d-q+1,q-d-2)}\ar[rr]\ar[uu]&&\scriptstyle\cdots\ar[rr]&&\scriptstyle \old{Q(0,c-1,0)}\ar[uu]&& \\
	&&&&&\\
	\scriptstyle\cdots\ar[uu]&&\scriptstyle\cdots\ar[uu]&&\cdots\\
	&\\
	\scriptstyle Q(0,c+d-q,1)\ar[uu]\ar[rr]&& \ \ \scriptstyle \old{Q(0,c+d-q+1,0)}\ar[uu]\\
	&\\
	\scriptstyle Q(0,c+d-q,0)\ar[uu]
}\] 
Note that, by Theorem \ref{simple-module-and-n-tube}, any object in the subset $\{Q(0,i,j)\mid c+d-q+1\leq i\leq c-1,q\leq j, c+d+1\leq p+q\leq c-1+q\}$ of (1) is not contained in a simple-minded system on $A$-$\stmod.$ Please see the location of $\mathcal{N}^{0}_{Q(0,c,d)}$ in the above diagram. 

\noindent{\bf Step two}: compute ${^{\bot}Q(0,c,d)^{\bot}}\cap {_{s}Q_{1}}$.	

\medskip

By Serre duality, 
\[\StHom_{A}(X,Y)\cong\D\StHom_{A}(Y,\nu\Omega(X))\cong\D\StHom_{A}(Y,\Omega(X))\cong\D\StHom_{A}(\Omega^{-1}(Y),X).\]
Thus 
\begin{equation}
\begin{aligned}\label{equ-2}
\RS Q(0,c,d)\cap {_{s}Q_{1}}=&\LS \Omega Q(0,c,d)\cap {_{s}Q_{1}}\\
=&\{Q(1,i,j)\mid i\leq c, c\leq i+j\leq c+d\},\\
\LS Q(0,c,d)\cap {_{s}Q_{1}}=&\RS \Omega^{-1} Q(0,c,d)\cap {_{s}Q_{1}}\\
=&\{Q(1,i,j)\mid c+1\leq i\leq c+d+1, c+d+1\leq i+j\}.
\end{aligned}
\end{equation}

\begin{equation}
\begin{aligned}\label{equ-3}
\mathcal{N}^{1}_{Q(0,c,d)}\colon=&
{^{\bot}Q(0,c,d)^{\bot}}\cap {_{s}Q_{1}}\\
=&{_{s}Q_{1}}\backslash((\RS Q(0,c,d)\cup\LS Q(0,a,b))\cap{_{s}Q_{1}})\\
=&{_{s}Q_{1}}\backslash ((\LS\Omega (0,c,d)\cup\RS\Omega^{-1} Q(0,c,d))\cap {_{s}Q_{1}})\\
=&\{Q(1,i,j)\mid c+1\leq i\leq c+d, i+j\leq c+d\}\\
&\cup\{Q(1,i,j)\mid c+d-q+2\leq i\leq c-1, i+j\leq c-1\}\\
&\cup\{Q(1,i,j)\mid c+d-q+2\leq i\leq c-1, c+d+1\leq i+j\leq c+q-1\}.
\end{aligned}
\end{equation}

\medskip

\noindent{\bf Step three}: compute ${^{\bot}Q(0,c,d)^{\bot}}\cap {_{s}\Gamma_{0}\cup{_{s}\Gamma_{1}}}$.	

\medskip

Now we determine ${^{\bot}Q(0,c,d)^{\bot}}\cap {_{s}\Gamma_{0}}$
for $Q(0,c,d)$. Take a set of quasi-simples \[\Sigma_{Q(0,c,d)}\colon=\{Q(0,c,0), Q(0,c+1,0),\cdots, Q(0,c+d,0)\}.\]
Take an object $(0,i,j)$ on $_{s}\Gamma_{0}$.  By \cite[2.2]{EK},
\begin{align}\label{equation-0}
\dim_{k}\StHom_{A}((0,i,j),Q(0,c,d))=\Sigma^{d}_{\ell=0}\dim_{k}\StHom_{A}((0,i,j),Q(0,c+\ell,0)).
\end{align}
Thus 
\begin{align}\label{LS-eq}
\LS Q(0,c,d)\cap {_{s}\Gamma_{0}}=(\bigcup^{d}_{\ell=0}\LS Q(0,c+\ell,0))\cap {_{s}\Gamma_{0}}.
\end{align}

By equation (\ref{eq-10}),  $Q(0,c+\ell,0)\in\RS(0,a-t,b-c-\ell)$ 
for every $\ell, t\in\mathbb{Z}$. Thus \[(0,a-t,b-c-\ell)\in\LS Q(0,c+\ell,0)\] for every $\ell, t\in\mathbb{Z}$. 
For an integer $h$, if  $\overline{h}\equiv h(mod\ p)\notin\{\overline{b-c-d}, \overline{b-c-d+1}, \cdots, \overline{b-c}\}$, then $\StHom_{A}((0,t,h),Q(0,c+\ell,0))=0$ for $\ell\in\{1,2,\cdots,d\}$ and $t\in\mathbb{Z}$ by equations (\ref{equation-0}) and (\ref{eq-10}). Therefore, 

\begin{equation}
\begin{aligned}	\label{equ-4}
\RS Q(0,c,d)\cap {_{s}\Gamma_{0}}&=\{(0,i,b+c+\ell)\mid i\in\mathbb{Z},\ell=0,1,2,\cdots, d\},\\
\LS Q(0,c,d)\cap {_{s}\Gamma_{0}}&=\RS \Omega Q(0,c,d)\cap {_{s}\Gamma_{0}}\\
&=\RS  Q(1,c,d)\cap {_{s}\Gamma_{0}}=\{0\}.
\end{aligned}
\end{equation}

\begin{equation}
\begin{aligned}	\label{equ-5}
\LS Q(0,c,d)\cap {_{s}\Gamma_{1}}&=\{(1,i,b+c-1+\ell)\mid i\in\mathbb{Z},\ell=0,1,2,\cdots, d\},\\
\RS Q(0,c,d)\cap {_{s}\Gamma_{1}}&=\LS \Omega Q(0,c,d)\cap {_{s}\Gamma_{1}}\\
&=\LS Q(1,c,d)\cap {_{s}\Gamma_{1}}=\{0\}.
\end{aligned}
\end{equation}
Thus 
\begin{equation}
\begin{aligned}\label{equ-6}
\RS Q(0,c,d)\cap({_{s}\Gamma_{0}}\cup{_{s}\Gamma_{1}})&=\{(0,i,b+c+\ell)\mid i\in\mathbb{Z},\ell=0,1,2,\cdots, d\},\\
\LS Q(0,c,d)\cap({_{s}\Gamma_{0}}\cup{_{s}\Gamma_{1}})&=\{(1,i,b+c-1+\ell)\mid i\in\mathbb{Z},\ell=0,1,2,\cdots, d\}.
\end{aligned}
\end{equation}

\begin{equation}
\begin{aligned}\label{equ-7}	
\mathcal{P}_{Q(0,c,d)}\colon=&	
{^{\bot}Q(0,c,d)^{\bot}}\cap({_{s}\Gamma_{0}}\cup{_{s}\Gamma_{1}})\\
=&{_{s}\Gamma_{0}}\cup{_{s}\Gamma_{1}}\backslash((\LS Q(0,c,d)\cup\RS Q(0,c,d))\cap({_{s}\Gamma_{0}}
\cup{_{s}\Gamma_{1}}))\\
=&\{(0,i,\overline{b+c+\ell})\mid i\in\mathbb{Z},\ell=d+1,d+2,\dots, p-1\}\\
\cup&\{(1,i,\overline{b+c+\ell})\mid i\in\mathbb{Z},\ell=d,d+1,\dots, p\}.
\end{aligned}
\end{equation}
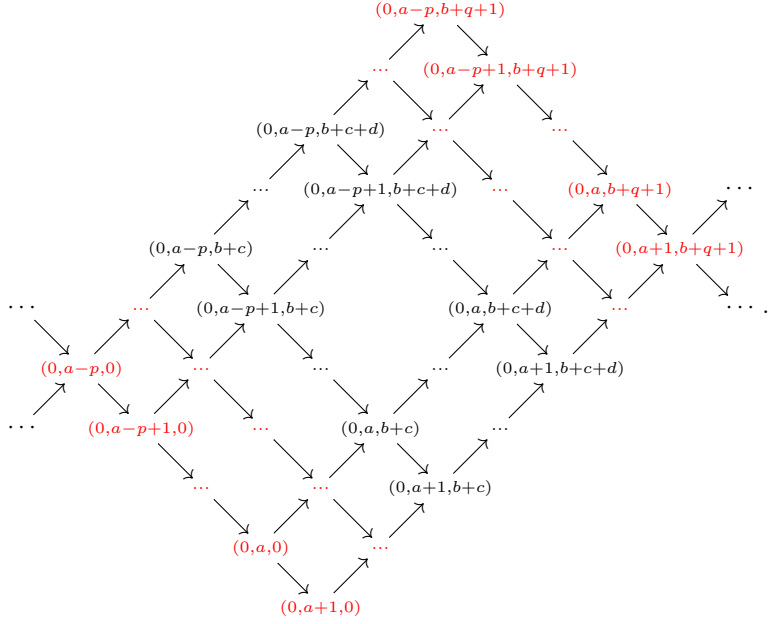
\vspace{-1cm}
\begin{small}\[\xymatrix@dr@R=16pt@C=16pt@!0{	
&&&&&&&&&&\cdots&&\\
&&\\
&&\scriptstyle\old{(0,a-p,b+q+1)}\ar[rr] &&\scriptstyle\scriptstyle\old{(0,a-p+1,b+q+1)} \ar[rr]& &\scriptstyle \old{\cdots}\scriptstyle \ar[rr] && \scriptstyle \scriptstyle\old{(0,a,b+q+1)}\ar[rr]&& \scriptstyle\scriptstyle\old{(0,a+1,b+q+1)} \ar[rr]\ar[uu]&&\ \  \cdots\,.& \\	
&& \\	
&&\scriptstyle\old{\cdots}\ar[rr]\ar[uu] &&\scriptstyle\old{\cdots}\ar[rr] \ar[uu]&& 	\scriptstyle \old{\cdots}\ar[rr] &&\scriptstyle\old{\cdots}\ar[uu]\ar[rr] &&\scriptstyle\old{\cdots}\ar[uu]&&\\	
&& \\	
&&\scriptstyle(0,a-p,b+c+d)\ar[rr]\ar[uu] &&\scriptstyle(0,a-p+1,b+c+d)\ar[rr] \ar[uu]&& \scriptstyle \cdots\ar[rr] &&\scriptstyle(0,a,b+c+d)\ar[uu]\ar[rr] &&\scriptstyle(0,a+1,b+c+d)\ar[uu]&&\\	
&&&& \\	
&&	\scriptstyle \cdots \ar[uu] &&\scriptstyle\cdots \ar[uu]& & \scriptstyle &&\scriptstyle\cdots\ar[uu] && \scriptstyle \cdots\ar[uu]\\	
&&& \\
&&\scriptstyle \scriptstyle(0,a-p,b+c) \ar[rr]\ar[uu]&&\ar[rr] \scriptstyle (0,a-p+1,b+c)\ar[uu] && \scriptstyle\cdots \ar[rr]\scriptstyle &&\scriptstyle(0,a,b+c)\ar[uu] \ar[rr]&&\scriptstyle(0,a+1,b+c)\ar[uu]\\
&&& \\
&&\scriptstyle \scriptstyle\old{\cdots}\ar[rr]\ar[uu]&&\ar[rr] \scriptstyle \old{\cdots}\ar[uu] && \scriptstyle\old{\cdots} \ar[rr] &&\scriptstyle\old{\cdots}\ar[uu] \ar[rr]&&\scriptstyle\old{\cdots}\ar[uu]\\
&&&\\
\cdots\ar[rr]&&\old{\scriptstyle(0,a-p,0)}\ar[rr]\ar[uu]  && \scriptstyle\old{\scriptstyle(0,a-p+1,0)}\ar[uu]\ar[rr] && \scriptstyle\old{ \cdots}\ar[rr]  &&\old{\scriptstyle(0,a,0)} \ar[rr]\ar[uu]&&\old{\scriptstyle(0,a+1,0)}\ar[uu] \\ 
&&&&&&&&&&\\
&&\cdots \ar[uu]   
}\] 
\end{small}
Thus ${^{\bot}Q(0,c,d)^{\bot}}\cap{_{s}\Gamma_{A}}=\mathcal{N}^{0}_{Q(0,c,d)}\cup\mathcal{N}^{1}_{Q(0,c,d)}\cup\mathcal{P}_{Q(0,c,d)}.$

\section{the construction of orthogonal systems}
The following conclusion provides us a necessary condition for an orthogonal system to be a simple-minded system on $A$-$\stmod.$

\begin{Lem}\label{necessary-condition-sms}
Let $A$	 be a 2-domestic Brauer graph algebra and $\mathcal{S}$  a simple-minded system on $A$-$\stmod$. Then $\mathcal{S}$  contains at least one non-periodic module.
\end{Lem}
\begin{proof}
By Theorem \ref{simple-module-and-n-tube},  any object in a homogeneous tube is not in $\mathcal{S}$. We show that $\mathcal{S}$  contains at least one non-periodic module. On the contrary,  if the set  $\mathcal{S}\cap(_{s}\Gamma_{0}\cup{_{s}\Gamma_{1}})$ is empty, then $\mathcal{S}$ is contained in $_{s}Q_{0}\cup{_{s}Q_{1}}\cup{_{s}Q'_{0}}\cup{_{s}Q'_{1}}$. By Corollary \ref{2-BGA-quasi-tube-bricks}, 
any object on a homogeneous tube is contained in the stable bi-perpendicular category of  $_{s}Q_{0}\cup{_{s}Q_{1}}\cup{_{s}Q'_{0}}\cup{_{s}Q'_{1}}$. Thus every  object on a homogeneous tube is in $^{\bot}\mathcal{S}^{\bot}$. It is a contradiction. Hence our conclusion holds. 
\end{proof}
According to Lemma \ref{necessary-condition-sms}, we shall mainly consider the orthogonal systems which contain at least one object for an Euclidean component. 
\subsection{The construction  of orthogonal systems on Euclidean components}
Without loss of generality, we always assume $p\leq q$ for any Euclidean component  in this section. Now we shall  present all orthogonal systems on $_{s}\Gamma_{0}\cup{_{s}\Gamma_{1}}$ containing the vertex  $(0,a,b)$ as follows. 

Take $(0,a_{1},b_{1})\colon=(0,a,b)$. By equation (\ref{bi-perp-1}),  one vertex  $(0,a_{2},b_{2})\in\mathcal{M}^{0}_{(0,a_{1},b_{1})}$ if and only if   $a_{1}-p<a_{2}<a_{1}, b_{1}<b_{2}<b_{1}+q$. 
Thus the set  $\{(0,a_{1},b_{1}),(0,a_{2},b_{2}),\cdots,(0,a_{k},b_{k})\}$ is an orthogonal system on $A$-$\stmod$ if and only if $a_{1}-p<a_{i}<a_{1}, b_{1}<b_{i}<b_{1}+q$,  $i=2,3,\cdots,k$, and 
if $a_{i}<a_{j}$, then $b_{j}<b_{i}$ for any $i,j=1,2,\cdots, k$.
Thus 
\begin{equation}
	\begin{aligned}\label{orthogonal-system-1}
		\mathcal{O}^{0}_{(0,a_{1},b_{1})}\colon=\{(0,a_{i},b_{i})\mid\ &a_{1}-p<a_{i}<a_{1}, b_{1}<b_{i}<b_{1}+q\ \text{for}\  i=2,\cdots, k, \text{and}, 
		\\&\text{if}\  a_{i}<a_{j}, \text{then}\  b_{j}<b_{i}\    \text{for\ any}\  i,j=1,2,\cdots, k\}.
	\end{aligned}
\end{equation}
is an orthogonal system on $A$-$\stmod$ containing $(0,a,b)$. 

Similarly, the set  $\{(1,c_{1},d_{1}),(1,c_{2},d_{2}),\cdots,(1,c_{\ell},d_{\ell})\}$ is an orthogonal system on $A$-$\stmod$ if and only if $c_{1}-p<c_{i}<c_{1}$ for $i=2,3,\cdots,\ell$, $d_{1}<d_{i}<d_{1}+q$ for any $i=1,2,\cdots, \ell$, and 
if $c_{i}<c_{j}$, then $d_{j}<d_{i}$ for any $i,j=1,2,\cdots, \ell$. Thus 
\begin{equation}
	\begin{aligned}\label{orthogonal-system-2}
		\mathcal{O}^{1}_{(1,c_{1},d_{1})}\colon=\{(1,c_{i},d_{i})\mid\ &\ c_{1}-p<c_{i}<c_{1}, 
		d_{1}<d_{i}<d_{1}+q\  \text{for}\  i=2,\cdots, \ell, \\
		&\text{and}\  \text{if}\ c_{i}<c_{j},\  \text{then}\  d_{j}<d_{i}\ \text{for\ any}\  i,j=1,2,\cdots, \ell\}.
	\end{aligned}
\end{equation}
is also an orthogonal system on $A$-$\stmod$.

By equation (\ref{bi-perp-7}), one vertex $(1,c,d)\in\mathcal{M}^{1}_{(0,a,b)}$ if and only if $(1,c,d)\in{^{\perp}(1,a+1,b)^{\perp}}\cap{_{s}\Gamma_{1}} $ if and only if $a_{1}-p+1\leq c\leq a_{1}, b_{1}+1\leq d\leq b_{1}+q$. Therefore  the set $\mathcal{O}^{0}_{(0,a_{1},b_{1})}\cup\{(1,c_{1},d_{1})\}$ is an orthogonal system on $A$-$\stmod$ if  and only if $a_{i}-p+1\leq c_{1}\leq a_{i}, b_{i}+1\leq d_{1}\leq b_{i}+q$
for any $i\in\{1,2,\cdots,k\}$. 
Therefore  the set $\mathcal{O}^{0}_{(1,a_{1},b_{1})}\cup\mathcal{O}^{1}_{(0,c_{1},d_{1})}$ is an orthogonal system on $A$-$\stmod$ if  and only if $a_{j}-p+1\leq c_{i}\leq a_{j}, b_{j}+1\leq d_{i}\leq b_{j}+q$. Thus we have the following proposition.

\begin{Prop}\label{Euclidean-orthogonal}
Let $A$ be a 2-domestic Brauer graph algebra and let  $\Gamma_{0}$  and $\Gamma_{1}$ be two Euclidean components. Take a subset $\mathcal{M}_{0}=\{(0,a_{i},b_{i})\mid 1\leq i\leq k\}$ on  $_{s}\Gamma_{0}$ and a subset $\mathcal{M}_{1}=\{(1,c_{i},d_{i})\mid 1\leq i\leq \ell\}$ on  $_{s}\Gamma_{1}$, respectively. Then
\begin{enumerate}
\item $\mathcal{M}_{0}$ is an orthogonal system if and only if $a_{1}-p<a_{i}<a_{1}, b_{1}<b_{i}<b_{1}+q$ for $i=2,3,\cdots,k$, and 
if $a_{i}<a_{j}$, then $b_{j}<b_{i}$ for any $i,j=1,2,\cdots, k$.
\item $\mathcal{M}_{1}$ is an orthogonal system if and only if $c_{1}-p<c_{i}<c_{1}$, $d_{1}<d_{i}<d_{1}+q$ for $i=2,3,\cdots,\ell$, and if $c_{i}<c_{j}$, then $d_{j}<d_{i}$ for any $i,j=1,2,\cdots, \ell$. 
\item We assume that  $\mathcal{M}_{0}$ and $\mathcal{M}_{1}$  are two orthogonal systems.
$\mathcal{M}_{0}\cup\mathcal{M}_{1}$ is an orthogonal system if and only if $a_{j}-p+1\leq c_{i}\leq a_{j}, b_{j}+1\leq d_{i}\leq b_{j}+q$.
for any $i,j=1,2,\cdots,k.$ 
\end{enumerate}
\end{Prop}

Note that the cardinality of a maximal orthogonal system on an Euclidean component $\mathbb{Z}\widetilde{A}_{p,q}$ may be $2,3,\cdots,p$.  We present the following example for an  Euclidean component $\mathbb{Z}\widetilde{A}_{3,3}$ over a $2$-domestic Brauer graph algebra.
\begin{Ex}\label{ortho-sys-Euclidean-comp}
Take the following  two Euclidean components  $_{s}\Gamma_{0}$ and $_{s}\Gamma_{1}$ of the form $\mathbb{Z}\widetilde{A}_{3,3}$  over a $2$-domestic Brauer graph algebra. Note that $(0,-2,3)$ and $(0,1,0)$ are identical.
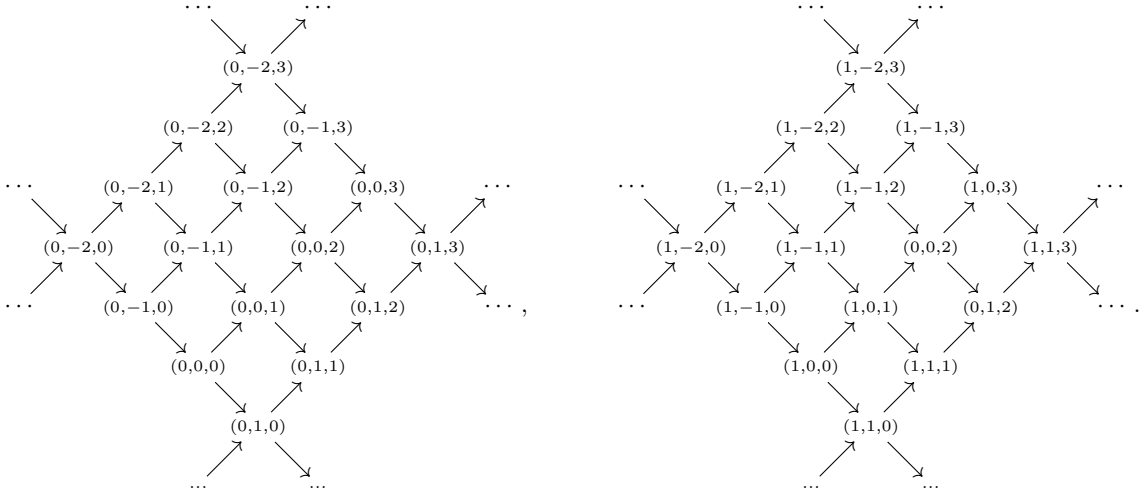
\vspace{-0.7cm}
\begin{small}\[\xymatrix@dr@R=16pt@C=16pt@!0{
	&&\cdots&&&&&&\cdots&&\\
	&&\\
	\cdots\ar[rr]&&\scriptstyle(0,-2,3)\ar[uu]\ar[rr] &&\scriptstyle(0,-1,3) \ar[rr]&& \scriptstyle(0,0,3)\ar[rr]&& \scriptstyle(0,1,3) \ar[rr]\ar[uu]&&\ \  \cdots\,,& \\	
	&& \\	
	&&\scriptstyle(0,-2,2)\ar[rr]\ar[uu] &&\scriptstyle(0,-1,2)\ar[rr] \ar[uu]&&\scriptstyle(0,0,2)\ar[uu]\ar[rr] &&\scriptstyle(0,1,2)\ar[uu]&&\\	
	&&& \\	
	&&\scriptstyle(0,-2,1) \ar[rr]\ar[uu]&&\ar[rr] \scriptstyle (0,-1,1)\ar[uu]  &&\scriptstyle(0,0,1)\ar[uu] \ar[rr]&&\scriptstyle(0,1,1)\ar[uu]\\
	&&& \\
	\cdots\ar[rr]&&\scriptstyle(0,-2,0)\ar[rr]\ar[uu]  && \scriptstyle(0,-1,0)\ar[uu]\ar[rr]  &&\scriptstyle(0,0,0) \ar[rr]\ar[uu]&&\scriptstyle(0,1,0)\ar[uu]\ar[rr]&& \scriptstyle\cdots\\ 
	&&\\
	&&\cdots\ar[uu]&&&&&&\scriptstyle\cdots\ar[uu]}
\xymatrix@dr@R=16pt@C=16pt@!0{
	&&\cdots&&&&&&\cdots&&\\
	&&\\
	\cdots\ar[rr]&&\scriptstyle(1,-2,3)\ar[uu]\ar[rr] &&\scriptstyle(1,-1,3) \ar[rr]&& \scriptstyle(1,0,3)\ar[rr]&& \scriptstyle(1,1,3) \ar[rr]\ar[uu]&&\ \  \cdots\,.& \\	
	&& \\	
	&&\scriptstyle(1,-2,2)\ar[rr]\ar[uu] &&\scriptstyle(1,-1,2)\ar[rr] \ar[uu]&&\scriptstyle(0,0,2)\ar[uu]\ar[rr] &&\scriptstyle(0,1,2)\ar[uu]&&\\	
	&&& \\	
	&&\scriptstyle(1,-2,1) \ar[rr]\ar[uu]&&\ar[rr] \scriptstyle (1,-1,1)\ar[uu]  &&\scriptstyle(1,0,1)\ar[uu] \ar[rr]&&\scriptstyle(1,1,1)\ar[uu]\\
	&&& \\
	\cdots\ar[rr]&&\scriptstyle(1,-2,0)\ar[rr]\ar[uu]  && \scriptstyle(1,-1,0)\ar[uu]\ar[rr]  &&\scriptstyle(1,0,0) \ar[rr]\ar[uu]&&\scriptstyle(1,1,0)\ar[uu]\ar[rr]&& \scriptstyle\cdots\\ 
	&&\\
	&&\cdots\ar[uu]&&&&&&\scriptstyle\cdots\ar[uu]
}
\] 
\end{small}

By (1) of Proposition \ref{Euclidean-orthogonal}, the sets $\mathcal{U}_{0}=\{(0,1,0),(0,-1,1)\}$ and $\mathcal{U}_{1}=\{(0,1,0),(0,-1,2),(0,0,1)\}$ are  all maximal orthogonal systems containing $(0,1,0)$ on $_{s}\Gamma_{0}$. 
Now we list all maximal orthogonal systems containing $(0,1,0)$ on ${_{s}\Gamma_{0}}\cup{_{s}\Gamma_{1}}$.

\noindent{\bf Case one:} All maximal orthogonal systems containing $\mathcal{U}_{1}$ on ${_{s}\Gamma_{0}}\cup{_{s}\Gamma_{1}}$.

By Serre duality, \[\mathcal{M}_{1}\colon={^{\bot}\mathcal{U}_{0}^{\bot}}\cap{_{s}\Gamma_{1}}=\{(1,-1,2),(1,-1,3),(1,0,1), (1,1,1)\}.\] 
By (2) of Proposition \ref{Euclidean-orthogonal}, there are four maximal orthogonal systems on $\mathcal{M}$ as follows.
\[\mathcal{V}_{1}=\{(1,-1,3),(1,0,1)\},\mathcal{V}_{2}=\{(1,-1,3),(1,1,1)\}, \mathcal{V}_{3}=\{(1,-1,2),(1,0,1)\}, \mathcal{V}_{4}=\{(1,-1,2),(1,1,1)\}.\]
Thus any maximal orthogonal system containing $\mathcal{U}_{0}$ on ${_{s}\Gamma_{0}}\cup{_{s}\Gamma_{1}}$ is of the form $\mathcal{U}_{0}\cup\mathcal{V}_{i}$ for $i=1,2,3,4.$ Note that the cardinality of $\mathcal{U}_{0}\cup\mathcal{V}_{i}$ is $4$ for each $i=1,2,3,4$. 

\noindent{\bf Case two:} All maximal orthogonal systems containing $\mathcal{U}_{1}$ on ${_{s}\Gamma_{0}}\cup{_{s}\Gamma_{1}}$.

By Serre duality, \[\mathcal{M}'_{1}\colon={^{\bot}\mathcal{U}_{1}^{\bot}}\cap{_{s}\Gamma_{1}}=\{(1,-1,3),(1,0,2), (1,1,1)\}.\] 
By (2) of Proposition \ref{Euclidean-orthogonal}, $\mathcal{M}'_{1}$ is an   orthogonal system. Hence  $\mathcal{U}_{1}\cup\mathcal{M}'_{1}$ is the only maximal orthogonal system containing $\mathcal{U}_{1}$ on ${_{s}\Gamma_{0}}\cup{_{s}\Gamma_{1}}$. Note that the cardinality of $\mathcal{U}_{1}\cup\mathcal{M}'_{1}$ is $6$. 
\end{Ex}

\subsection{The construction of orthogonal systems on  quasi-tubes}
In this subsection, we shall construct  orthogonal systems on  stable generalized standard quasi-tubes. 

\begin{Def}\label{triangle-areas}
Let $A$ be a self-injective algebra, and let $Q$ be a quasi-tube and $(a,b)$ a vertex on the quasi-tube $Q$. Take the set 
\[\bigtriangleup_{(a,b)}\colon=\{(i,j)\mid a\leq i\leq a+b, a\leq i+j\leq a+b\}.\]
$\bigtriangleup_{(a,b)}$  is called {\bf triangle area of vertices $(a,b)$} and we call the number $b$ the {\bf height} of $\bigtriangleup_{(a,b)}.$
\end{Def}
Note that $\bigtriangleup_{(a,b)}$  is  the set of vertices contained in the triangle with vertices $(a,0)$, $(a,b)$ and $(a+b,0)$ and that, if $b\leq-1$, then we always assume that $\bigtriangleup_{(a,b)}=\varnothing.$  
We shall show that the intersection between stable bi-perpendicular category of an orthogonal system on  an Euclidean component and a quasi-tube is a union of  triangle areas of some vertices.

Take a vertex $(0,a,b)\in{_{s}Q_{0}}$. By equation (\ref{eq-12}),  $^{\bot}(0,a,b)^{\bot}\cap{_{s}Q_{0}}=\bigtriangleup_{Q(0,1,q-2)}$. First, we shall construct all orthogonal systems on triangle area  $\bigtriangleup_{Q(0,1,q-2)}$ by induction on the height $j$ of triangle area  $\bigtriangleup_{Q(0,i,j)}$.  
Without loss of generality, We shall always fix the orthogonal system appeared in this subsection with the second coordinate in an ascending order. For example, the orthogonal system $\{Q(0,1,0), Q(0,2,0),\cdots,Q(0,q-1,0)\}$ on $\bigtriangleup_{Q(0,1,q-2)}$ satisfies the rule. 

\begin{Ex}\label{ortho-example}
We list all orthogonal systems on triangle area $\bigtriangleup_{Q(0,1,q-2)}$ for $q-2=0,1,2,3$ as follows. 
\begin{enumerate}[$(1)$]
\item The case $q-2=0$. There is only one vertex $Q(0,1,0)$ on triangle area $\bigtriangleup_{Q(0,1,0)}$.
Thus $\{Q(0,1,0)\}$ is the only orthogonal system on triangle area  $\bigtriangleup_{Q(0,1,0)}$.

\item The case $q-2=1$.   There are four orthogonal systems on triangle area  $\bigtriangleup_{Q(0,1,1)}$ as follows.
\[\{Q(0,1,0)\}, \{Q(0,2,0)\}, \{Q(0,1,0),Q(0,2,0)\}, \{Q(0,1,1)\}.\]

\item The case $q-2=2$. We  list all orthogonal systems on triangle area  $\bigtriangleup_{Q(0,1,2)}$ as follows.
\begin{enumerate}[$(a)$]
\item Orthogonal systems with one object: Any  subset of $\bigtriangleup_{Q(0,1,1)}$ containing only one object.
\item Orthogonal systems with two objects: \[ \{Q(0,1,0),Q(0,2,0)\}, \{Q(0,1,0),Q(0,2,1)\}, \{Q(0,1,0),Q(0,3,0)\},\] 
\[\{Q(0,1,1),Q(0,3,0)\},\{Q(0,1,2),Q(0,2,0)\}, \{Q(0,2,0),Q(0,3,0)\}.\]
\item Orthogonal systems with three objects: \[  \{Q(0,1,0),Q(0,2,0),Q(0,3,0)\}.\]
\end{enumerate}

\item The case $q-2=3$. We  only list all maximal orthogonal systems on triangle area  $\bigtriangleup_{Q(0,1,q-2)}$, since there are a lot of orthogonal systems in this case. Note that every orthogonal system is a subset of a  maximal orthogonal system.

\begin{enumerate}[$(a)$]
\item Maximal orthogonal systems with two objects:
\[\{Q(0,1,1),Q(0,3,1)\},\{Q(0,1,3),Q(0,2,1)\}.\]
\item Maximal orthogonal systems with three objects: 
\[\{Q(0,1,0),Q(0,2,0),Q(0,3,1)\}, \{Q(0,1,0),Q(0,2,1),Q(0,4,0)\},\{Q(0,1,0),Q(0,2,2),Q(0,3,0)\},\] 
\[ \{Q(0,1,1),Q(0,3,0),Q(0,4,0)\},\{Q(0,1,2),Q(0,2,0),Q(0,4,0)\},\{Q(0,1,3),Q(0,2,0),Q(0,3,0)\}.\]
\item Maximal orthogonal systems with four objects: \[\{Q(0,1,0),Q(0,2,0),Q(0,3,0),Q(0,4,0)\}.\]
\end{enumerate}
\end{enumerate}
\end{Ex}

Now we construct orthogonal systems on the triangle area $\bigtriangleup_{Q(0,1,q-2)}$ by induction as follows. 
Take a vertex $Q(0,1,i_{1})$ on $\bigtriangleup_{Q(0,1,q-2)}$, where $0\leq i_{1}\leq q-2$.  By equation (\ref{equ-1}), we consider two cases as follows. 

\begin{enumerate}[$(1)$]
\item If $i_{1}=0$, then  ${^{\bot}Q(0,1,i_{1})^{\bot}}\cap\bigtriangleup_{Q(0,1,q-2)}=\bigtriangleup_{Q(0,2,q-3)}$. 	
\item  If $i_{1}\geq1$, then  ${^{\bot}Q(0,1,i_{1})^{\bot}}\cap\bigtriangleup_{Q(0,1,q-2)}=\bigtriangleup_{Q(0,2,i_{1}-2)}\cup\bigtriangleup_{Q(0,i_{1}+2,q-i_{1}-3)}$. 
\end{enumerate}
Note that, the triangle area $\bigtriangleup_{Q(0,2,i_{1}-2)}=\varnothing$ when $i_{1}=1$. Note that   both $i_{1}-2$ and $q-i_{1}-3$ are less than $p-2$ strictly, and that the triangle areas $\bigtriangleup_{Q(0,2,i_{1}-2)}$ and $\bigtriangleup_{Q(0,i_{1}+2,q-i_{1}-3)}$ are disjoint and they are mutual stable orthogonal on $A$-$\stmod$. 

We add one more object on the set $\{Q(0,1,i_{1})\}$ from  ${^{\bot}Q(0,1,i_{1})^{\bot}}\cap\bigtriangleup_{Q(0,1,q-2)}$ as follows. 

\begin{enumerate}[$(a)$]
\item For (1), we  take a vertex $(0,2,i_{2})$ on triangle area  $\bigtriangleup_{Q(0,2,q-3)}$. By equation (\ref{equ-1}),  \[{^{\bot}\{Q(0,1,0), Q(0,2,i_{2}) \}^{\bot}}\cap\bigtriangleup_{Q(0,1,q-2)}=\bigtriangleup_{Q(0,3,i_{2}-2)}\cup\bigtriangleup_{Q(0,i_{2}+3,q-i_{2}-4)}.\]

\item For $i_{1}\geq2$ in (2), we  take a vertex $(0,2,i_{2})$ on  triangle area  $\bigtriangleup_{Q(0,2,i_{1}-2)}$.
 By equation (\ref{equ-1}),  \[{^{\bot}\{Q(0,1,i_{1}), Q(0,2,i_{2}) \}^{\bot}}\cap\bigtriangleup_{Q(0,1,q-2)}=\bigtriangleup_{Q(0,3,i_{2}-2)}\cup\bigtriangleup_{Q(0,i_{2}+3,i_{1}-i_{2}-3)}\cup\bigtriangleup_{Q(0,i_{1}+2,q-i_{1}-3)}.\]
Note that $i_{2}<i_{1}$ in this case.

\item For $i_{1}=1$ in (2), we take a vertex $(3,i_{3})$ on triangle area  $\bigtriangleup_{Q(0,3,q-4)}$. 
\begin{enumerate}[$(i)$]
\item If $i_{3}=0$, then  ${^{\bot}\{Q(0,1,i_{1}), Q(0,3,0) \}^{\bot}}\cap\bigtriangleup_{Q(0,1,q-2)}=\bigtriangleup_{Q(0,4,q-5)}$. 
\item If $i_{3}\geq1$, then  \[{^{\bot}\{Q(0,1,i_{1}), Q(0,3,i_{3}) \}^{\bot}}\cap\bigtriangleup_{Q(0,1,q-2)}=\bigtriangleup_{Q(0,4,i_{3}-2)}\cup\bigtriangleup_{Q(0,i_{3}+4,q-i_{3}-5)}.\]
\end{enumerate}
\end{enumerate}

Note that  the heights of triangle areas for above cases are less than the height of triangle areas on ${^{\bot}Q(0,1,i_{1})^{\bot}}\cap\bigtriangleup_{Q(0,1,q-2)}$, respectively. They 
are disjoint union of triangle areas and  those  triangle areas are mutual stable orthogonal on $A$-$\stmod$. 

Take a vertex $Q(0,\ell,i_{\ell})$ on triangle area  $\bigtriangleup_{Q(0,1,q-2)}$, where $0\leq i_{\ell}\leq q-2$.  Let $\mathcal{S}_{\ell}$ be an orthogonal system containing the vertex $Q(0,\ell,i_{\ell})$ such that $i\geq\ell$ for each $Q(0,i,j)$ in $\mathcal{S}_{\ell}$. 
By equation (\ref{equ-1}), $\mathcal{S}_{\ell}$ is contained  in the disjoint of two triangle areas of vertices $Q(0,\ell+1,i_{\ell}-2)$ and  $Q(0,\ell+i_{\ell},q-\ell-i_{\ell}-2)$. If the stable bi-perpendicular category ${^{\bot}\mathcal{S}_{\ell}^{\bot}}\cap\bigtriangleup_{Q(0,1,q-2)}$ is non-zero, then by induction we choose a non-zero element on  ${^{\bot}\mathcal{S}_{\ell}^{\bot}}\cap\bigtriangleup_{Q(0,1,q-2)}$, which is contained in a union triangle areas of some vertices. Note that the height of those vertices is lower than $i_{\ell}-2$ or $q-\ell-i_{\ell}-2$. Continuing this operation, with finitely many steps, we  get a maximal orthogonal system with the second coordinate is larger than or equals to $\ell$.
In this way, we can construct all maximal orthogonal systems, and 
in particular, we  get all orthogonal systems on the triangle area  $\bigtriangleup_{Q(0,1,q-2)}$.

Now we determine orthogonal systems on  the triangle areas  $\bigtriangleup_{Q(0,1,q-2)}\cup\bigtriangleup_{Q(1,1,q-2)}$.
By  equation (\ref{eq-13}), 	$^{\bot}(0,a,b)^{\bot}\cap{_{s}Q_{1}}=\bigtriangleup_{Q(1,1,q-2)}$. 
Thus 
\[^{\bot}(0,a,b)^{\bot}\cap({_{s}Q_{0}}\cup{_{s}Q_{1}} )=\bigtriangleup_{Q(0,1,q-2)}\cup\bigtriangleup_{Q(1,1,q-2)}.\]
Note that $\bigtriangleup_{Q(1,1,q-2)}=\Omega(\bigtriangleup_{Q(0,1,q-2)})$. Using a similar method with the above, we may construct all orthogonal systems
on $^{\bot}(0,a,b)^{\bot}\cap{_{s}Q_{1}}$.

By equation (\ref{equ-3}), \[^{\bot}Q(0,i,j)^{\bot}\cap{_{s}Q_{1}}={\Omega(Q(0,i,j))^{\bot}\cap{^{\bot}\Omega^{-1}(Q(0,i,j))} }\cap{_{s}Q_{1}}={Q(1,i,j)^{\bot}\cap{^{\bot}Q(1,i+1,j)}}\cap{_{s}Q_{1}}.\]
Take a vertex $Q(0,1,i_{1})$ on triangle areas  $\bigtriangleup_{Q(0,1,q-2)}\cup\bigtriangleup_{Q(1,1,q-2)}$. Then 
\begin{equation}
\begin{aligned}\label{biper-orth-sys-1}	
^{\bot}Q(0,1,i_{1})^{\bot}\cap(\bigtriangleup_{Q(0,1,q-2)}\cup\bigtriangleup_{Q(1,1,q-2)})=&(\bigtriangleup_{Q(0,2,i_{1}-2)}\cup\bigtriangleup_{Q(0,2+i_{1},q-i_{1}-3)})\\
&\cup(\bigtriangleup_{Q(1,2,i_{1}-1)}\cup\bigtriangleup_{Q(1,3+i_{1},q-i_{1}-4)})\\
&\cup\{Q(1,1,j)\mid i_{1}<j<q-1\}.
\end{aligned}
\end{equation}
We shall construct an orthogonal system $\mathcal{S}$ containing the vertex $Q(0,1,i_{1})$. According to equation (\ref{biper-orth-sys-1}), we consider two cases. 

\medskip

\noindent{\bf Case one:}  If $\mathcal{S}$ contains an object of the form $Q(1,1,j)$ for $i_{1}<j<q-1$. Without loss of generality, we assume that $Q(1,1,i_{h})$($i_{1}<i_{h}<q-1$) is in $\mathcal{S}$. Thus 
\begin{equation}
\begin{aligned}\label{biper-orth-sys-2}	^{\bot}\{Q(0,1,i_{1}),Q(1,1,i_{h})\}^{\bot}&\cap(\bigtriangleup_{Q(0,1,q-2)}\cup\bigtriangleup_{Q(1,1,q-2)})\\
&=(\bigtriangleup_{Q(0,2,i_{1}-2)}\cup\bigtriangleup_{Q(0,2+i_{1},i_{h}-i_{1}-2)}\cup\bigtriangleup_{Q(0,2+i_{h},q-i_{h}-3)})\\
&\cup(\bigtriangleup_{Q(1,2,i_{1}-1)}\cup\bigtriangleup_{Q(1,3+i_{1},i_{h}-i_{1}-3)}\cup\bigtriangleup_{Q(1,2+i_{h},q-i_{h}-3)}).
\end{aligned}
\end{equation}

\noindent{\bf Case two:}  If $\mathcal{S}$ contains no object of the form $Q(1,1,j)$ for $i_{1}<j<q-1$.
By equation (\ref{biper-orth-sys-1}), $\mathcal{S}$ is contained in the following triangle areas of some vertices:
\[(\bigtriangleup_{Q(0,2,i_{1}-2)}\cup\bigtriangleup_{Q(0,2+i_{1},q-i_{1}-3)})
\cup(\bigtriangleup_{Q(1,2,i_{1}-1)}\cup\bigtriangleup_{Q(1,3+i_{1},q-i_{1}-4)}).\]

In conclusion, it is not hard to know that the stable bi-perpendicular category of an orthogonal systems on the triangle areas  $\bigtriangleup_{Q(0,1,q-2)}\cup\bigtriangleup_{Q(1,1,q-2)}$  is contained in triangle areas of some vertices.  By induction on the height of vertices,  we may construct all maximal  orthogonal systems on triangle areas $\bigtriangleup_{Q(0,1,q-2)}\cup\bigtriangleup_{Q(1,1,q-2)}$, thus all orthogonal systems on  $\bigtriangleup_{Q(0,1,q-2)}\cup\bigtriangleup_{Q(1,1,q-2)}$.
We shall compute the cardinality of a maximal orthogonal system on triangle areas $\bigtriangleup_{Q(0,1,q-2)}\cup\bigtriangleup_{Q(1,1,q-2)}$. 

\begin{Them}\label{cardinality-of-orth-sys}
Let $A$ be a symmetric algebra  and let  $Q_{0}$ and $Q_{1}$ be  generalized standard quasi-tubes  of rank $q$ satisfying $\Omega(_{s}Q_{0})={_{s}Q_{1}}$. 
\begin{enumerate}[$(1)$]
\item Let $\mathcal{S}_{1}$ be a maximal orthogonal system on triangle areas $\bigtriangleup_{Q(0,x,y-1)}\cup\bigtriangleup_{Q(1,x,y)}$ over $_{s}Q_{0}\cup{_{s}Q_{1}}$ for $y\leq q-2$.  Then the cardinality of $\mathcal{S}_{1}$ is $y+1$  for $y=0,1,\cdots,q-2$. 
\item Let $\mathcal{S}_{2}$ be a maximal orthogonal system on the triangle areas  $\bigtriangleup_{Q(0,v,w)}\cup\bigtriangleup_{Q(1,v+1,w-1)}$ over $_{s}Q_{0}\cup{_{s}Q_{1}}$  for $w\leq q-2$.  Then the cardinality of $\mathcal{S}_{2}$ is $w+1$  for $w=0,1,\cdots,q-2$. 	
\item Let $\mathcal{S}$ be a maximal orthogonal system on the triangle areas $\bigtriangleup_{Q(0,j,i)}\cup\bigtriangleup_{Q(1,j,i)}$ over  $_{s}Q_{0}\cup{_{s}Q_{1}}$  for $i\leq q-2$. 
 Then the cardinality of $\mathcal{S}$ is $i+1$ for $i=0,1,\cdots,q-2$.
\end{enumerate}
\end{Them}
\begin{proof}
We prove (1), (2) and (3) by induction on the heights $y$, $w$ and $i$, respectively. 
For (1) and (2), we only prove (1), since the proof of (2) is dual. For (3), we proof only the case for $j=1$, since $\tau(Q(0,j,i))=Q(0,j-1,i)$ and $\tau$ is a stable self-equivalence on $A$-$\stmod$. Therefore we assume that $j=1.$

We first induct on the height $y$ for conclusion (1). If $y=0$, then $\bigtriangleup_{Q(0,x,-1)}\cup\bigtriangleup_{Q(1,x,0)}=\{Q(1,x,0)\}$ is  a maximal orthogonal system on triangle areas  $\bigtriangleup_{Q(0,x,-1)}\cup\bigtriangleup_{Q(1,x,0)}.$ Thus conclusion (1) is true for $y=0$.
If  $y=1$, then \[\bigtriangleup_{Q(0,x,0)}\cup\bigtriangleup_{Q(1,x,1)}=\{Q(0,x,0),Q(1,x,0),Q(1,x,1), Q(1,x+1,0)\}.\] It is easy to know that there are two maximal orthogonal systems on triangle areas  $\bigtriangleup_{Q(0,x,0)}\cup\bigtriangleup_{Q(1,x,1)}$  as follows.
\[\{Q(0,x,0),Q(1,x,1)\}, \{Q(1,x,0),Q(1,x+1,0)\}.\]
Thus  conclusion (1) is true for $y=1$. 

Then we  induct on the height $i$ for conclusion (3).  If $i=0$, then \[\bigtriangleup_{Q(0,1,0)}\cup\bigtriangleup_{Q(1,1,0)}=\{Q(0,1,0),Q(1,1,0)\}\] and $Q(1,1,0)=\Omega(Q(0,1,0))$. It is easy to know that  $\{Q(0,1,0)\}$ and $\{Q(1,1,0)\}$ are all maximal orthogonal systems on triangle areas $\bigtriangleup_{Q(0,1,0)}\cup\bigtriangleup_{Q(1,1,0)}$. Therefore conclusion (3) is true for $i=0$.

If $i=1$, then \[\bigtriangleup_{Q(0,1,1)}\cup\bigtriangleup_{Q(1,1,1)}=\{Q(0,1,0),Q(0,1,1), Q(0,2,0), Q(1,1,0), Q(1,1,1),Q(1,2,0)\}.\] There are two maximal orthogonal systems $\{Q(0,1,0),Q(0,2,0)\}$ and $\{Q(0,1,1)\}$ on triangle area $\bigtriangleup_{Q(0,1,1)}$.
According to equation (\ref{equ-3}), the set $\{Q(0,1,0),Q(0,2,0)\}$ is still maximal on  triangle areas $\bigtriangleup_{Q(0,1,1)}\cup\bigtriangleup_{Q(1,1,1)}$. Thus the set $\{Q(1,1,0),Q(1,2,0)\}$ is maximal on triangle areas  $\bigtriangleup_{Q(0,1,1)}\cup\bigtriangleup_{Q(1,1,1)}$. By equation (\ref{equ-3}), \[{^{\bot}\{Q(0,1,1)\}^{\bot}}\cap(\bigtriangleup_{Q(0,1,1)}\cup\bigtriangleup_{Q(1,1,1)})=\{Q(1,2,0)\}.\]
Thus the set  $\{Q(0,1,1),Q(1,2,0)\}$ is a maximal orthogonal system on  triangle areas $\bigtriangleup_{Q(0,1,1)}\cup\bigtriangleup_{Q(1,1,1)}$.  Using similar method, we know that the sets $\{Q(0,1,0),Q(1,1,1)\}$ and 
 $\{Q(0,2,0),Q(1,1,0)\}$ are  also maximal orthogonal systems on $\bigtriangleup_{Q(0,1,1)}\cup\bigtriangleup_{Q(1,1,1)}$. Hence there are five maximal orthogonal systems on triangle areas $\bigtriangleup_{Q(0,1,1)}\cup\bigtriangleup_{Q(1,1,1)}$, which have  the same cardinality $2$. Therefore conclusion (3) is true for $i=1$.
 
 Now we assume that  conclusions (1) and (2) hold for $y=k$ and $w=k$ respectively, and that  conclusion (3) is true for $i=k$.  We prove that (1), (2) and (3) are true for $y=k+1, w=k+1$, and $i=k+1,$ respectively. 
Let $\mathcal{S}$ be a maximal orthogonal system on triangle areas $\bigtriangleup_{Q(0,1,k+1)}\cup\bigtriangleup_{Q(1,1,k+1)}.$
Take \[\mathcal{L}_{1}\colon=\{Q(0,1,j)\mid 0\leq j\leq k+1\}\cup\{Q(1,1,j)\mid 0\leq j\leq k+1\}.\]

\medskip

\noindent{\bf Claim one:} $\mathcal{L}_{1}\cap\mathcal{S}\neq\varnothing$.

\medskip

Otherwise, $\mathcal{S}\subseteq\bigtriangleup_{Q(0,2,k)}\cup\bigtriangleup_{Q(1,2,k)}$ and it is a maximal orthogonal system on triangle areas  $\bigtriangleup_{Q(0,2,k)}\cup\bigtriangleup_{Q(1,2,k)}$. Take \[\mathcal{L}_{2}\colon=\{Q(0,2,j)\mid 0\leq j\leq k\}\cup\{Q(1,2,j)\mid 0\leq j\leq k\}.\] By induction, the cardinality of $\mathcal{S}$ is $k+1$. We claim that $\mathcal{L}_{2}\cap\mathcal{S}\neq\varnothing$. Otherwise, $\mathcal{S}\subseteq\bigtriangleup_{Q(0,3,k-1)}\cup\bigtriangleup_{Q(1,3,k-1)}$. By induction, the cardinality of $\mathcal{S}$ is $k$. It is a contradiction. Thus there is an object in $\mathcal{S}$ which is contained in $\mathcal{L}_{2},$ denoted by $S$. Then we assume that $S=(s,2,t)$ for $s=0, 1$ and $0\leq t\leq k$. By equations (\ref{equ-1}) and (\ref{equ-3}), we always have  $Q(1,1,0)\in{^{\bot}\mathcal{S}^{\bot}}$. Therefore $\mathcal{S}$ is not maximal on triangle areas  $\bigtriangleup_{Q(0,1,k+1)}\cup\bigtriangleup_{Q(1,1,k+1)}$. It is a contradiction. Therefore $\mathcal{L}_{1}\cap\mathcal{S}\neq\varnothing$. Thus Claim one is proved. 

Let $\mathcal{S}_{1}$ be a maximal orthogonal system on triangle areas  $\bigtriangleup_{Q(0,x,k)}\cup\bigtriangleup_{Q(1,x,k+1)}$.
 
 \medskip
 
\noindent{\bf Claim two:} $\mathcal{S}_{1}\cap\{Q(1,r,x+k+1-r)\mid x\leq r\leq x+k+1\}\neq\varnothing$.

\medskip

Otherwise, $\mathcal{S}_{1}\cap\{Q(1,r,x+k+1-r)\mid x\leq r\leq x+k+1\}=\varnothing$, then $\mathcal{S}_{1}$ is a maximal orthogonal  system on triangle areas $\bigtriangleup_{Q(0,x,k)}\cup\bigtriangleup_{Q(1,x,k)}$. Take
\[\mathcal{H}_{1}=\{Q(0,r,x+k-r)\mid x\leq r\leq x+k\}\cup\{Q(1,r,x+k-r)\mid x\leq r\leq x+k\}.\]
Dual with Claim one, $\mathcal{S}_{1}\cap\mathcal{H}_{1}\neq\varnothing.$ 
Take $M\in\mathcal{S}_{1}\cap\mathcal{H}_{1}\neq\varnothing.$ 
If $M=Q(1,x,k)$, by equations (\ref{equ-1}) and (\ref{equ-3}), then 
$Q(0,x+k+1,0)\in{^{\bot}\mathcal{S}_{1}^{\bot}}$. It is a contradiction. If $M=\mathcal{H}_{1}\backslash\{Q(1,x,k)\}$, then $Q(0,x,k+1)\in{^{\bot}\mathcal{S}_{1}^{\bot}}$. It contradicts the condition that $\mathcal{S}_{1}$ is maximal on triangle areas $\bigtriangleup_{Q(0,x,k)}\cup\bigtriangleup_{Q(1,x,k+1)}$. Therefore $\mathcal{S}_{1}\cap\{Q(1,r,x+k+1-r)\mid x\leq r\leq x+k+1\}\neq\varnothing$. Thus Claim two is proved. 

Take $M_{1}=Q(1,s,x+k+1-s)\in\mathcal{S}_{1}\cap\{Q(1,r,x+k+1-r)\mid x\leq r\leq x+k+1\}$.
Then
\begin{equation}
\begin{aligned}
&{^{\bot}\mathcal{S}_{1}^{\bot}}\cap(\bigtriangleup_{Q(0,x,k)}\cup\bigtriangleup_{Q(1,x,k+1)}\\
&=(\bigtriangleup_{Q(0,x,s-x-2)}\cup\bigtriangleup_{Q(0,s,x+k-s)})
\cup(\bigtriangleup_{Q(1,x,s-x-1)}\cup\bigtriangleup_{Q(1,s+1,x+k-s-1)})\\
&=(\bigtriangleup_{Q(0,x,s-x-2)}\cup\bigtriangleup_{Q(1,x,s-x-1)})\cup(\bigtriangleup_{Q(0,s,x+k-s)}\cup\bigtriangleup_{Q(1,s+1,x+k-s-1)}).\\
\end{aligned}
\end{equation}
Note that $s-x-1\leq k$ and $x+y-s-2\leq k-1$. 
By induction of (1) and (2), the cardinality of a maximal orthogonal system on triangle areas $\bigtriangleup_{Q(0,x,s-x-2)}\cup\bigtriangleup_{Q(1,x,s-x-1)}$ is $(s-x-1)+1=s-x$ and the cardinality of a maximal orthogonal system on triangle areas $\bigtriangleup_{Q(0,s,x+k-s)}\cup\bigtriangleup_{Q(1,s+1,x+k-s-1)}$ is $(x+k-s)+1=x+k-s+1$. Since triangle areas  $\bigtriangleup_{Q(0,x,s-x-2)}\cup\bigtriangleup_{Q(1,x,s-x-1)}$ and $\bigtriangleup_{Q(0,s,x+k-s)}\cup\bigtriangleup_{Q(1,s+1,x+k-s-1)}$ are mutual orthogonal, the cardinality of $\mathcal{M}$ is $1+(s-x)+(x+k-s+1)=k+2.$ Thus conclusion (1) is proved. 

Now we prove (3). By Claim one, $\mathcal{L}_{1}\cap\mathcal{S}\neq\varnothing$.
We consider two cases as follows.

\noindent{\bf Case one:} $\mathcal{S}\cap\{Q(0,1,j)\mid 0\leq j\leq k+1\}\neq\varnothing.$
We assume that $S_{1}=Q(0,1,j_{1})$ in $\mathcal{S}$.
By equation ( \ref{equ-1}), 
\begin{equation}
\begin{aligned}
&{^{\bot}S_{1}^{\bot}}\cap(\bigtriangleup_{Q(0,1,k+1)}\cup\bigtriangleup_{Q(1,1,k+1)}\\
&=(\bigtriangleup_{Q(0,2,j_{1}-2)}\cup\bigtriangleup_{Q(0,j_{1}+2,k-j_{1})})
\cup(\bigtriangleup_{Q(1,2,j_{1}-1)}\cup\bigtriangleup_{Q(0,j_{1}+3,k-j_{1}-1)})\\
&=(\bigtriangleup_{Q(0,2,j_{1}-2)}\cup\bigtriangleup_{Q(1,2,j_{1}-1)})\cup(\bigtriangleup_{Q(0,j_{1}+2,k-j_{1})}\cup\bigtriangleup_{Q(0,j_{1}+3,k-j_{1}-1)}).\\
\end{aligned}
\end{equation}
Note that $j_{1}-1\leq k$ and $k-j_{1}\leq k$.  By induction of (1) and (2), the cardinality of  a maximal orthogonal system on triangle areas $\bigtriangleup_{Q(0,2,j_{1}-2)}\cup\bigtriangleup_{Q(1,2,j_{1}-1)}$ is $(j_{1}-1)+1=j_{1}$ and the cardinality of  a maximal orthogonal system on triangle areas  $\bigtriangleup_{Q(0,j_{1}+2,k-j_{1})}\cup\bigtriangleup_{Q(0,j_{1}+3,k-j_{1}-1)}$ is $(k-j_{1})+1$. Since triangle areas  $\bigtriangleup_{Q(0,2,j_{1}-2)}\cup\bigtriangleup_{Q(1,2,j_{1}-1)}$ and $\bigtriangleup_{Q(0,j_{1}+2,k-j_{1})}\cup\bigtriangleup_{Q(0,j_{1}+3,k-j_{1}-1)}$ are mutual orthogonal, the cardinality of $\mathcal{S}$ is $1+j_{1}+(k-j_{1}+1)=k+2.$
Thus our conclusion (3) holds in this case. 

\noindent{\bf Case two:} $\mathcal{S}\cap\{Q(1,1,j)\mid 0\leq j\leq k+1\}\neq\varnothing.$
We assume that $S_{1}=Q(0,1,k_{1})$ in $\mathcal{S}$.
By equation (\ref{equ-1}), 
\begin{equation}
\begin{aligned}
&{^{\bot}S_{1}^{\bot}}\cap(\bigtriangleup_{Q(0,1,k+1)}\cup\bigtriangleup_{Q(1,1,k+1)}\\
&=(\bigtriangleup_{Q(1,2,k_{1}-2)}\cup\bigtriangleup_{Q(1,k_{1}+2,k-k_{1})})
\cup(\bigtriangleup_{Q(0,1,k_{1}-1)}\cup\bigtriangleup_{Q(0,k_{1}+2,k-k_{1})})\\
&=(\bigtriangleup_{Q(1,2,k_{1}-2)}\cup\bigtriangleup_{Q(0,1,k_{1}-1)})\cup(\bigtriangleup_{Q(1,k_{1}+2,k-k_{1})}\cup\bigtriangleup_{Q(0,k_{1}+2,k-k_{1})}).\\
\end{aligned}
\end{equation}
Note that $k_{1}-1\leq k$ and $k-k_{1}\leq k$. 
By induction of (1) and (2), the cardinality of  a maximal orthogonal system on triangle areas  $\bigtriangleup_{Q(1,2,k_{1}-2)}\cup\bigtriangleup_{Q(0,1,k_{1}-1)}$ is $(k_{1}-1)+1=k_{1}$ and the cardinality of  a maximal orthogonal system on triangle areas  $\bigtriangleup_{Q(1,k_{1}+2,k-k_{1})}\cup\bigtriangleup_{Q(0,k_{1}+2,k-k_{1})}$ is $(k-k_{1})+1$. Since triangle areas  $\bigtriangleup_{Q(1,2,k_{1}-2)}\cup\bigtriangleup_{Q(0,1,k_{1}-1)}$ and $\bigtriangleup_{Q(1,k_{1}+2,k-k_{1})}\cup\bigtriangleup_{Q(0,k_{1}+2,k-k_{1})}$ are mutual orthogonal, the cardinality of $\mathcal{S}$ is $1+k_{1}+(k-k_{1}+1)=k+2.$
Thus our conclusion (3) holds in this case. 

Hence conclusion (3) is proved by summarizing the above proof.
\end{proof}

\subsection{A maximal orthogonal systems on $A$-$\stmod$}	
Let $\mathcal{M}$ be an orthogonal system on the stable Euclidean component $_{s}\Gamma_{0}$ over $A$-$\stmod$.  Let $\mathcal{P}_{i}$ be the set of quasi-simples on $Q_{i}$ and $\mathcal{Q}_{i}$ the set of quasi-simples on $Q'_{i}$ for $i=0,1$. Without loss of generality, by equation (\ref{orthogonal-system-1}), we assume that $\mathcal{M}$ equals to 
\begin{equation}
\begin{aligned}\label{orthogonal-system-gamma-1}
\mathcal{O}^{0}_{(0,a_{1},b_{1})}\colon=\{(0,a_{i},b_{i})\mid\ &a_{1}-p<a_{i}<a_{1}, b_{1}<b_{i}<b_{1}+q\ \text{for}\  i=2,\cdots, k, 
\\&\text{and}, \text{if}\  a_{i}<a_{j}, \text{then}\  b_{j}<b_{i}\    \text{for\ any}\  i,j=1,2,\cdots, k\}.
\end{aligned}
\end{equation}
and we assume that $(0,a_{1},b_{1})=(0,a,b)$.  

\begin{Prop}\label{cardinality-of-ortho-sys-2}
Let $\mathcal{M}$, $_{s}Q_{0}$ and $_{s}Q_{1}$ {\rm(}resp. $_{s}Q'_{0}$ and $_{s}Q'_{1}${\rm)} be as above. Then the cardinality of a maximal orthogonal system on ${^{\bot}\mathcal{M}^{\bot}}\cap{_{s}Q_{0}}\cup{_{s}Q_{1}}$ {\rm(}resp. ${^{\bot}\mathcal{M}^{\bot}}\cap{_{s}Q'_{0}}\cup{_{s}Q'_{1}}${\rm)} is $q-k$ {\rm(}resp. $p-k${\rm)}.
\end{Prop}
\begin{proof}
We only prove that the cardinality of a maximal orthogonal system on ${^{\bot}\mathcal{M}^{\bot}}\cap{_{s}Q_{0}}\cup{_{s}Q_{1}}$ is $q-k$, since the other case is similar to deal with. By Proposition \ref{sms-quasi-simple} and the construction of an orthogonal system on an Euclidean component,  if $i\neq j,$ then  $(\LS(0,a_{i},b_{i})\cap\mathcal{P}_{0})\cap(\LS(0,a_{j},b_{j})\cap\mathcal{P}_{0})=\varnothing$.
Thus $\mid\LS(\mathcal{M})\cap\mathcal{P}_{0}\mid=k$. Without loss of generality, we may assume that 
\[\LS(\mathcal{M})\cap\mathcal{P}_{0}=\{Q(0,b_{i},0)\mid i=1,2,\cdots,k\}\] 
with the second coordinate in an ascending way. Note that $\RS(\mathcal{M})\cap\mathcal{P}_{1}=\Omega(\LS(\mathcal{M})\cap\mathcal{P}_{0})$ and that  $(\RS(\mathcal{M})\cap\mathcal{P}_{1})\cap(\LS(\mathcal{M})\cap\mathcal{P}_{0})=\varnothing$.
	
Take \[\mathcal{R}=\{Q(j,b_{i}+1,b_{i+1}-b_{i}-2)\mid i=1,\cdots, k-1, j=0,1\} \cup\{Q(j,\overline{b_{k}+1},b_{1}+p-b_{k}-2)\mid j=0,1\},\]
where $\overline{b_{k}+1}$ is the congruence class of modulo $q$. By equation (\ref{eq-12}), ${^{\bot}\mathcal{M}^{\bot}}\cap({_{s}Q_{0}}\cup{_{s}Q_{1}})$
is a disjoint union of $\triangle_{v}$ for $v\in\mathcal{R}$.

By Subsection 3.2, we may construct a maximal orthogonal system for each triangle area as follows. \[\triangle_{Q(0,b_{i}+1,b_{i+1}-b_{i}-2)}\cup\triangle_{Q(1,b_{i}+1,b_{i+1}-b_{i}-2)}\  \text{for}\  i=1,\cdots,k,\  \text{and}\  \triangle_{Q(0,\overline{b_{k}+1},b_{1}+q-b_{k}-2)}\cup\triangle_{Q(1,\overline{b_{k}+1},b_{1}+q-b_{k}-2)}.\] 
Note that if $b_{i+1}=b_{i}+1$, then $b_{i+1}-b_{i}-2=-1$. In this case, we assume that $\triangle_{Q(0,b_{i}+1,b_{i+1}-b_{i}-2)}=\varnothing.$
Put those maximal orthogonal systems together, denoted by $\mathcal{S}_{1}$. Therefore $\mathcal{S}_{1}$ is a maximal orthogonal system on ${^{\bot}\mathcal{M}^{\bot}}\cap({_{s}Q_{0}}\cup{_{s}Q'_{0}})$. 
	
Now we determine the cardinality of $\mathcal{S}_{1}$. Since $(0,a,b)$ is in $\mathcal{M}$, $\mathcal{S}_{1}$ is contained in triangle areas  $\bigtriangleup_{Q(0,1,q-2)}\cup\bigtriangleup_{Q(1,1,q-2)}$.  Since $\LS(\mathcal{S}_{1})\cap\mathcal{P}_{0}$ is an orthogonal system on $\bigtriangleup_{Q(0,1,q-2)}\cup\bigtriangleup_{Q(1,1,q-2)}$. By the construction of $\mathcal{S}_{1}$, $\LS(\mathcal{M})\cap\mathcal{P}_{0}$ are mutual orthogonal  with $\mathcal{S}_{1}$ on $A$-$\stmod$. Thus $(\LS(\mathcal{M})\cap(\mathcal{P}_{0}\backslash\{Q(0,0,0)\}))\cup\mathcal{S}_{1}$ is a maximal orthogonal system on triangle areas  $\bigtriangleup_{Q(0,1,q-2)}\cup\bigtriangleup_{Q(1,1,q-2)}$.
By Theorem \ref{cardinality-of-orth-sys}, $\mid\mathcal{S}_{1}\mid=q-1-(k-1)=q-k.$
By the similar method, we may also construct  a maximal orthogonal system  $\mathcal{S}_{2}$ on ${^{\bot}\mathcal{M}^{\bot}}\cap({_{s}Q'_{0}}\cup{_{s}Q'_{1}})$, and the cardinality of $\mathcal{S}_{2}$ is $p-m$.
\end{proof}

Let $\mathcal{M}$,  $\mathcal{S}_{1}$ and  $\mathcal{S}_{2}$ be as in the proof of Proposition \ref{cardinality-of-ortho-sys-2}.  By equations (\ref{bi-perp-8}) and (\ref{equ-7}), ${^{\bot}(\mathcal{M}\cup\mathcal{S}_{1}\cup\mathcal{S}_{2})^{\bot}}\cap{_{s}\Gamma_{1}}$	is a disjoint union of some rectangle areas of some vertices. By Subsection 3.1, we may construct a maximal orthogonal system for each rectangle area. Put those maximal orthogonal systems together, denoted by  $\mathcal{M}'$. $\mathcal{M}'$ is  a maximal orthogonal system on ${^{\bot}(\mathcal{M}\cup\mathcal{S}_{1}\cup\mathcal{S}_{2})^{\bot}}\cap{_{s}\Gamma_{1}}$. Take  $\mathcal{S}=\mathcal{M}\cup\mathcal{M}'\cup\mathcal{S}_{1}\cup\mathcal{S}_{2}$. Hence $\mathcal{S}$ is a maximal orthogonal system on $A$-$\stmod$.
Since $\mathcal{S}_{1}\cup\mathcal{S}_{2}$ is a maximal orthogonal system on ${^{\bot}\mathcal{M}^{\bot}}\cap({_{s}Q_{0}}\cup{_{s}Q_{1}}\cup{_{s}Q'_{0}}\cup{_{s}Q'_{1}})$, $\mathcal{M}$ is maximal on  ${^{\bot}(\mathcal{S}_{1}\cup\mathcal{S}_{2})^{\bot}}\cap{_{s}\Gamma_{0}}$  and $\mathcal{M}'$ is maximal on  ${^{\bot}(\mathcal{S}_{1}\cup\mathcal{S}_{2})^{\bot}}\cap{_{s}\Gamma_{1}}$.  $\mathcal{S}_{1}\cup\mathcal{S}_{2}$ is a maximal orthogonal system on ${^{\bot}\mathcal{M}'^{\bot}}\cap({_{s}Q_{0}}\cup{_{s}Q_{1}}\cup{_{s}Q'_{0}}\cup{_{s}Q'_{1}})$.  The cardinality of a maximal orthogonal system on ${^{\bot}\mathcal{M}^{\bot}}\cap({_{s}Q_{0}}\cup{_{s}Q_{1}}\cup{_{s}Q'_{0}}\cup{_{s}Q'_{1}})$ equals to the cardinality of a maximal orthogonal system on ${^{\bot}\mathcal{M}'^{\bot}}\cap({_{s}Q_{0}}\cup{_{s}Q_{1}}\cup{_{s}Q'_{0}}\cup{_{s}Q'_{1}})$. 

\begin{Prop}\label{cardinality-of-ortho-sys-3}
Let  $\mathcal{M}$ and $\mathcal{M}'$ be as above. Then the cardinality of $\mathcal{M}$ equals to the cardinality of $\mathcal{M}'$.
\end{Prop}
\begin{proof}
Otherwise, ${\mid\mathcal{M}\mid}\neq{\mid\mathcal{M}'\mid}$. Without loss of generality, we assume that ${\mid\mathcal{M}\mid}>{\mid\mathcal{M}'\mid}$. 
Then ${\mid\LS(\mathcal{M})\cap\mathcal{P}_{0}\mid}>{\mid\LS(\mathcal{M}')\cap\mathcal{P}_{0}\mid}$ and ${\mid\RS(\mathcal{M})\cap\mathcal{P}_{1}\mid}>{\mid\RS(\mathcal{M}')\cap\mathcal{P}_{1}\mid}$. Similar with the above construction, we know that  ${^{\bot}{\mathcal{M}'}^{\bot}}\cap({_{s}Q_{0}}\cup{_{s}Q_{1}})$
is a disjoint union of $\triangle_{v'}$ for $v'$ in some set $\mathcal{R}'$(compared with $\mathcal{R}$ in Proposition \ref{cardinality-of-ortho-sys-2}) and  ${\mid\mathcal{R}\mid}>{\mid\mathcal{R'}\mid}$.  Thus it is not hard to know that the cardinality of  ${^{\bot}\mathcal{M}'^{\bot}}\cap({_{s}Q_{0}}\cup{_{s}Q_{1}}\cup{_{s}Q'_{0}}\cup{_{s}Q'_{1}})$ is strictly less than the cardinality of  ${^{\bot}\mathcal{M}^{\bot}}\cap({_{s}Q_{0}}\cup{_{s}Q_{1}}\cup{_{s}Q'_{0}}\cup{_{s}Q'_{1}})$. It is a  contradiction. 
\end{proof}

\begin{Cor}\label{cardinality-of-ortho-sys-4}
Let $A$ be a $2$-domestic Brauer graph algebra and $\mathcal{S}$ a maximal orthogonal system on $A$-$\stmod$ containing at least one object on an Euclidean component. Then the cardinality of $\mathcal{S}$ is the number of non-projective simple $A$-modules.
\end{Cor}
\begin{proof}
By Proposition \ref{cardinality-of-ortho-sys-3},  ${\mid\mathcal{S}\cap{_{s}\Gamma_{0}}\mid}={\mid\mathcal{S}\cap{_{s}\Gamma_{1}}\mid}=k$. By Proposition   \ref{cardinality-of-ortho-sys-2}, ${\mid\mathcal{S}\cap ({_{s}Q_{0}}\cup{_{s}Q_{1}})\mid}=q-k$ and ${\mid\mathcal{S}\cap ({_{s}Q'_{0}}\cup{_{s}Q'_{1}})\mid}=p-k$.
Thus the cardinalty of  $\mathcal{S}$ is $2k+(q-k)+(p-k)=q+p=n$ by Theorem \ref{two-domestic-Brauer-graph-algebras}, the number of non-projective simple $A$-modules.  
\end{proof}
Thus every maximal orthogonal system containing at least one object for an Euclidean component on $A$-$\stmod$ has the same cardinality, the number of  non-projective simple $A$-modules.

\section{A sufficient and necessary condition for an orthogonal system to be a simple-minded system}
Let $A$ be a 2-domestic Brauer graph algebra and $\mathcal{S}$  a maximal orthogonal system on $A$-$\stmod$ containing at least one object for an Euclidean component. Our main task in this section is to prove that $\mathcal{S}$ is a simple-minded system
on $A$-$\stmod$, see Theorem \ref{sms-BGA}. We assume that $\mathcal{M}=\{M_{1},\cdots,M_{\ell}\}$ is the subset of $\mathcal{S}$ which is contained in the Euclidean component $_{s}\Gamma_{0}$ and we may assume that  $\mathcal{M'}=\{M'_{1},\cdots,M'_{\ell}\}$ is the subset of $\mathcal{S}$ contained in $_{s}\Gamma_{1}$ by Proposition \ref{cardinality-of-ortho-sys-3}. Without loss of generality, we assume that the object $M_{1}=(0,a_{1},b_{1})=(0,a,b)$ (resp. $M'_{1}=(1,a'_{1},b'_{1})$) and $M_{i}=(0,a_{i},b_{i})$ (resp. $M'_{i}=(1,a'_{i},b'_{i})$)  for $\i=2,\cdots,\ell.$ By equations (\ref{bi-perp-1}) and (\ref{bi-perp-2}),   ${^{\bot}\mathcal{M}^{\bot}}\cap{_{s}\Gamma_{0}}$  consists of  a disjoint of  $\ell$ rectangle areas. By equation (\ref{eq-12}),  ${^{\bot}\mathcal{M}^{\bot}}\cap({_{s}Q'_{0}}\cup{_{s}Q'_{1}})$ is a disjoint of several  triangle areas and $\mathcal{S}\cap({_{s}Q'_{0}}\cup{_{s}Q'_{1}})\subseteq{^{\bot}\mathcal{M}^{\bot}}\cap({_{s}Q'_{0}}\cup{_{s}Q'_{1}})$.

We first consider a special case that every object in $\mathcal{S}\cap({_{s}Q'_{0}}\cup{_{s}Q'_{1}})$ and  $\mathcal{S}\cap({_{s}Q_{0}}\cup{_{s}Q_{1}})$  is a quasi-simple.  
Since ${^{\bot}\mathcal{M}^{\bot}}\cap({_{s}Q'_{0}}\cup{_{s}Q'_{1}})$ is a disjoint of several  triangle areas and they are mutual stable orthogonal on $A$-$\stmod$, without loss of generality, we  need only figure out how can we construct a maximal orthogonal system which contains only quasi-simples on one of the triangle areas, such as $^{\bot}({^{\bot}\{M_{1},M_{2}\}^{\bot}\cap\square_{M_{1},M_{2}}})^{\bot}\cap({_{s}Q'_{0}}\cup{_{s}Q'_{1}})$.
Note that $^{\bot}({^{\bot}\{M_{1},M_{2}\}^{\bot}\cap\square_{M_{1},M_{2}}})^{\bot}\cap({_{s}Q'_{0}}\cup{_{s}Q'_{1}})=\bigtriangleup_{Q(0,1,a_{2}-a_{1}-2)'}\cup\bigtriangleup_{Q(1,1,a_{2}-a_{1}-2)'}$ and 
all quasi-simples on triangle areas  $\bigtriangleup_{Q(0,1,a_{2}-a_{1}-2)'}\cup\bigtriangleup_{Q(1,1,a_{2}-a_{1}-2)'}$ state as follows.
\[\mathcal{Q}\colon=\{Q(0,i,0)'\mid 1\leq i\leq a_{2}-a_{1}-1\}\cup\{Q(1,i,0)'\mid 1\leq i\leq a_{2}-a_{1}-1\}.\]
Let $\mathcal{W}\subseteq\mathcal{Q}$ be a maximal orthogonal system on triangle areas  $\bigtriangleup_{Q(0,1,a_{2}-a_{1}-2)'}\cup\bigtriangleup_{Q(1,1,a_{2}-a_{1}-2)'}$. By Theorem \ref{cardinality-of-orth-sys},  the cardinality of $\mathcal{W}$ is $a_{2}-a_{1}-1$. The following proposition characterizes the maximal orthogonal system $\mathcal{W}$.

\begin{Prop}\label{ortho-quasi-1}
Let $\mathcal{Q}$ and $\mathcal{W}$ be as the above. 
Then $\mathcal{W}$ is of the form $\{Q(1,i,0)'\mid 1\leq i\leq\ell\}\cup\{Q(0,i,0)'\mid \ell+1\leq i\leq a_{2}-a_{1}-1\}$ for an integer  $0\leq\ell\leq a_{2}-a_{1}-1.$
\end{Prop}
\begin{proof}
there are two  cases to be considered as follows.

\noindent{\bf Case one:} $Q(0,1,0)'$ is in $\mathcal{W}$.\\ 
Then $Q(1,1,0)'$ and $Q(1,2,0)'$ are not in $\mathcal{W}$ by equation (\ref{equ-1}). Consider a subset of $\mathcal{Q}$ as follows.
\[\mathcal{Q}_{1}\colon=\{Q(0,i,0)'\mid 2\leq i\leq a_{2}-a_{1}-1\}\cup\{Q(1,i,0)'\mid 3\leq i\leq a_{2}-a_{1}-1\}.\]
By equation (\ref{equ-1}), $\mathcal{W}\subseteq\mathcal{Q}_{1}$.   
We claim that $Q(0,2,0)'$ is in $\mathcal{W}$. Otherwise, \[\mathcal{W}=\{Q(0,1,0)'\}\cup(\mathcal{W}\cap(\bigtriangleup_{Q(0,3,a_{2}-a_{1}-4)'}\cup\bigtriangleup_{Q(1,3,a_{2}-a_{1}-4)'})).\]
By Theorem \ref{cardinality-of-orth-sys}, the cardinality of $\mathcal{W}\cap(\bigtriangleup_{Q(0,3,a_{2}-a_{1}-4)'}\cup\bigtriangleup_{Q(1,3,a_{2}-a_{1}-4)'})$ is $a_{2}-a_{1}-4+1=a_{2}-a_{1}-3$. Then the cardinality of $\mathcal{W}$ is $a_{2}-a_{1}-3+1=a_{2}-a_{1}-2<a_{2}-a_{1}-1.$ It is a contradiction. Thus $Q(0,2,0)'$ is in $\mathcal{W}$. Proceeding this process, $\mathcal{W}=\{Q(0,i,0)'\mid 1\leq i\leq a_{2}-a_{1}-1\}$.

\noindent{\bf Case two:} $Q(0,1,0)'$ is not in $\mathcal{W}$.

Since $Q(0,1,0)'$ is not in $\mathcal{W}$, by Theorem  \ref{cardinality-of-orth-sys}, $Q(1,1,0)'$ is in $\mathcal{W}$.  By equation (\ref{equ-1}), \[\mathcal{W}=\{Q(1,1,0)'\}\cup(\mathcal{W}\cap(\bigtriangleup_{Q(0,2,a_{2}-a_{1}-3)'}\cup\bigtriangleup_{Q(1,2,a_{2}-a_{1}-3)'})).\] Then $Q(1,2,0)'\in\mathcal{W}$ or $Q(0,2,0)'\in\mathcal{W}$. If $Q(0,2,0)'\in\mathcal{W}$, by Case one, \[\mathcal{W}=\{Q(1,1,0)'\}\cup\{Q(0,k,0)'\mid 2\leq k\leq a_{2}-a_{1}-1\}.\] 
If $Q(1,2,0)'\in\mathcal{W},$ then  \[\mathcal{W}=\{Q(1,1,0)',Q(1,2,0)'\}\cup(\mathcal{W}\cap(\bigtriangleup_{Q(0,3,a_{2}-a_{1}-4)'}\cup\bigtriangleup_{Q(1,3,a_{2}-a_{1}-4)'})).\]
Proceeding this process, there exists an positive integer $r$ such that  \[\mathcal{W}=\{Q(1,k,0)'\mid 1\leq k\leq r\}\cup\{Q(0,k,0)'\mid r+1\leq k\leq a_{2}-a_{1}-1\}.\]

By combining Case one and Case two, our conclusion holds. 
\end{proof}

We state the general form of Proposition \ref{ortho-quasi-1} as follows and we will not give a  proof here. 
\begin{Cor}\label{quasi-tube-simple}
Let $\mathcal{W}$ be a maximal orthogonal system on $A$-$\stmod$ which consists of quasi-simples on $^{\bot}({^{\bot}\{M_{k},M_{k+1}\}^{\bot}\cap \square_{M_{k},M_{k+1}}})^{\bot}\cap({_{s}Q'_{0}}\cup{_{s}Q'_{1}})$  for  $1\leq k\leq\ell$. Then  $\mathcal{W}$ is  of  the  form   $\{Q(1,m,0)'\mid a_{k}-a_{1}+1\leq m\leq r \}\cup\{Q(0,m,0)'\mid r+1\leq m\leq a_{k+1}-a_{1}-1\}$ for $0\leq r\leq a_{k+1}-a_{1}-1.$
\end{Cor}
Note that we assume that $M_{\ell+1}$ is isomorphic  to $M_{1}$ and $a_{\ell+1}=a_{1}+p$ on Corollary \ref{quasi-tube-simple}.
By Proposition \ref{ortho-quasi-1} and Corollary  \ref{quasi-tube-simple}, we may assume that $\mathcal{S}\cap({_{s}Q_{0}}\cup{_{s}Q_{1}})$ is of the form
\begin{equation}
\begin{aligned}\label{orth-sys-2}
&\{Q(1,1,0)',Q(1,2,0)',\cdots, Q(1,t_{1}-1,0)'\}\cup\{Q(0,t_{1},0)',\cdots,Q(0,a_{2}-a_{1}-1,0)'\}\cup\\
&\{Q(1,a_{2}-a_{1}+1,0)',\cdots,Q(1,t_{2}-1,0)'\}\cup\{Q(0,t_{2},0)',\cdots,Q(0,a_{3}-a_{1}-1,0)'\}\\
&\cup,\cdots,\cup\\
&\{Q(1,a_{r}-a_{1}+1,0)',\cdots,Q(1,t_{r}-1,0)'\}\cup\{Q(0,t_{r},0)',\cdots, Q(0,a_{r+1}-a_{1}-1,0)'\}\\
&\cup,\cdots,\cup\\
&\{Q(1,a_{\ell}-a_{1}+1,0)',\cdots,Q(1,t_{\ell}-1,0)'\}\cup\{Q(0,t_{\ell},0)',\cdots,Q(0,p-1,0)'\}, 
\end{aligned}
\end{equation}
where 
$a_{r}-a_{1}+2\leq t_{r}\leq a_{r+1}-a_{1}-1$ for $1\leq r\leq\ell-1$ and $a_{\ell}-a_{1}+2\leq t_{\ell}\leq p-1.$ 

\begin{Prop}\label{Euclidean-1}
Let $\mathcal{S}$, $\mathcal{M}$, $\mathcal{M'}$ and $\{t_{k}\mid 1\leq k\leq\ell\}$ be as the stated above. Then the following subsets of objects are contained in $\mathcal{F}(\mathcal{S})$.
\begin{enumerate}[$(1)$]
\item $\{(0,a_{1}+j,b_{1})\mid j=0,1,\cdots, t_{1}-1\}\cup\{(0,a_{2}-j,b_{2})\mid j=0,1,\cdots,a_{2}-a_{1}-t_{1}\}$.
	
\item $\{(0,a_{r}+j,b_{r})\mid j=0,1,\cdots,t_{r}-a_{r}+a_{1}-1\}\cup\{(0,a_{r+1}-j,b_{r+1})\mid j=0,1,\cdots,a_{r+1}-a_{1}-t_{r}\}$ for $r=2,3,\cdots,\ell-1$.
	
\item $\{(0,a_{\ell}+j,b_{\ell})\mid j=0,1,\cdots,t_{\ell}-a_{\ell}+a_{1}-1\}\cup\{(0,a_{1}+p-j,b_{1}-q)\mid j=0,1,\cdots,p+1-t_{\ell}\}.$
\end{enumerate}
\end{Prop}
\begin{proof}
Without loss of generality, we only prove  (1), since other cases are similar to handle. 

Consider the following triangles:
\[(0,a_{1},b_{1})\xrightarrow{}(0,a_{1}+1,b_{1})\xrightarrow{} Q(1,1,0)'\xrightarrow{} (1,a_{1}+1,b_{1}+1),\]
\begin{equation}\label{seqss-1}
\begin{aligned}
(0,a_{1}+k,b_{1})\xrightarrow{}(0,a_{1}+k+1,b_{1})\xrightarrow{} Q(1,k+1,0)'\xrightarrow{} (1,a_{1}+k+1,b_{1}+1)
\end{aligned}
\end{equation}
for $1\leq k\leq t_{1}-3$.
Since $(0,a_{1},b_{1})$ and $Q(1,1,0)'$ are in $\mathcal{S}$, $(0,a_{1}+1,b_{1})\in\mathcal{F}(\mathcal{S})$. By  triangle (\ref{seqss-1}), it follows from induction that $(0,a_{1}+j,b_{1})\in\mathcal{F}(\mathcal{S})$ for $j=0,1,\cdots, t_{1}-1.$

Dually, consider triangles as follows.
\[Q(0,a_{2}-a_{1}-1,0)'\xrightarrow{}(0,a_{2}-1,b_{2})\xrightarrow{}(0,a_{2},b_{2})\xrightarrow{} Q(1,a_{2}-a_{1},0)',\]
\begin{equation}\label{seqss-2}
\begin{aligned}
Q(0,a_{2}-a_{1}-k-1,0)'\xrightarrow{}(0,a_{2}-k-1,b_{2})\xrightarrow{}(0,a_{2}-k,b_{2})\xrightarrow{} Q(1,a_{2}-a_{1}-k,0)'
\end{aligned}
\end{equation}
for $1\leq k\leq a_{2}-a_{1}-t_{1}-1$.
Since $(0,a_{2},b_{2})$ and $Q(0,a_{2}-a_{1}-1,0)'$ are in $\mathcal{S}$, $(0,a_{2}-1,b_{2})\in\mathcal{F}(\mathcal{S})$. By the triangle (\ref{seqss-2}), it follows from induction that $\{Q(0,t_{1},0)',\cdots,Q(0,a_{2}-a_{1}-1,0)'\}\subseteq\mathcal{F}(\mathcal{S})$.
\end{proof}
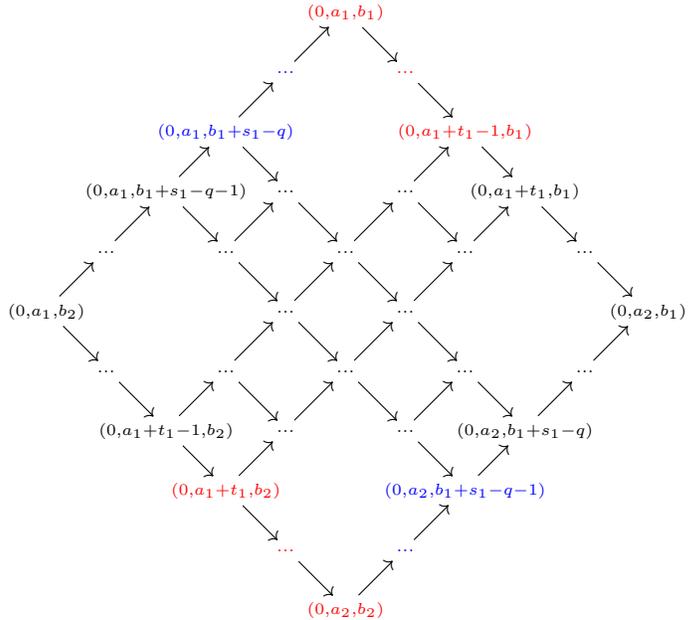
\vspace{-1cm}
\begin{small}
\[\xymatrix@dr@R=16pt@C=16pt@!0{
&&\scriptstyle\old{(0,a_{1},b_{1})}\ar[rr] &&\scriptstyle\scriptstyle\old{\cdots} \ar[rr]& &\scriptstyle \old{(0,a_{1}+t_{1}-1,b_{1})}\ar[rr] && \scriptstyle \scriptstyle(0,a_{1}+t_{1},b_{1})\ar[rr]&& \scriptstyle\scriptstyle\cdots\ar[rr]&&\ \  \scriptstyle(0,a_{2},b_{1})\,.&&\\	
&& \\	
&&\scriptstyle\new{\cdots}\ar[uu] &&\scriptstyle&& \scriptstyle \cdots\ar[rr]\ar[uu] &&\scriptstyle\cdots\ar[uu]&& &&\scriptstyle\scriptstyle\cdots\ar[uu]&&\\	
&&& \\	
&&	\scriptstyle\new{(0,a_{1},b_{1}+s_{1}-q)}\ar[rr] \ar[uu] &&\scriptstyle\cdots\ar[rr]& & \scriptstyle \cdots\ar[uu]\ar[rr]  &&\scriptstyle\cdots\ar[rr]\ar[uu] && \scriptstyle\cdots\ar[rr]  &&\scriptstyle(0,a_{2},b_{1}+s_{1}-q)\ar[uu]\\	
&&& \\
&&\scriptstyle(0,a_{1},b_{1}+s_{1}-q-1)\ar[uu]\ar[rr] && \scriptstyle\cdots\ar[uu]\ar[rr]  && \scriptstyle\cdots\ar[rr] \ar[uu] &&\scriptstyle\cdots\ar[rr]\ar[uu] &&\scriptstyle\cdots\ar[uu]\ar[rr]&&\scriptstyle\new{(0,a_{2},b_{1}+s_{1}-q-1)}\ar[uu]\\
&&& \\
&&\scriptstyle\cdots \ar[uu] && \scriptstyle && \scriptstyle\cdots\ar[uu]\ar[rr] &&\scriptstyle\cdots\ar[uu]&&\scriptstyle &&\scriptstyle\new{\cdots}\ar[uu]\\ 
&&\\
&&\scriptstyle(0,a_{1},b_{2}) \ar[uu] \ar[rr]&&\scriptstyle\cdots\ar[rr]&&\scriptstyle(0,a_{1}+t_{1}-1,b_{2})\ar[uu]\ar[rr]&&\scriptstyle\old{(0,a_{1}+t_{1},b_{2})}\ar[rr]\ar[uu] &&\scriptstyle\old{\cdots}\ar[rr]&&\scriptstyle\old{(0,a_{2},b_{2})} \ar[uu]
}\] 
\end{small}

The above  diagram describes the location of  objects  in (1) of Proposition \ref{Euclidean-1}  and (1) of Proposition \ref{Euclidean-2}. Specifically, the red one is the objects in (1) of Proposition \ref{Euclidean-1} and the blue one is the objects in (1) of Proposition \ref{Euclidean-2}. Note that both objects  $(0,a_{1},b_{1})$ and $(0,a_{2},b_{2})$ are in $\mathcal{S}$.

Dually, we may assume that  $\mathcal{S}\cap({_{s}Q_{0}}\cup{_{s}Q_{1}})$ is of the form: 
\begin{equation}
\begin{aligned}\label{orth-sys-1}
&\{Q(0,q-1,0),Q(0,q-2,0),\cdots, Q(0,s_{1},0)\}\cup\{Q(1,s_{1}-1,0),\cdots,Q(1,q-(b_{1}-b_{2}-1),0)\}\cup\\
&\{Q(0,q-(b_{1}-b_{2}+1),0),\cdots,Q(0,s_{2},0)\}\cup\{Q(1,s_{2}-1,0),\cdots,Q(1,q-(b_{1}-b_{3}-1),0)\}\\
&\cup,\cdots,\cup\\
&\{Q(0,q-(b_{1}-b_{r}+1),0),\cdots,Q(0,s_{r},0)\}\cup\{Q(1,s_{r}-1,0),\cdots, Q(1,q-(b_{1}-b_{r+1}-1),0)\}\\
&\cup,\cdots,\cup\\
&\{Q(0,q-(b_{1}-b_{\ell}+1),0),\cdots,Q(0,s_{\ell},0)\}\cup\{Q(1,s_{\ell}-1,0),\cdots,Q(1,1,0)\}.
\end{aligned}
\end{equation}
where $q-(b_{1}-b_{r+1}+1)\leq s_{r}\leq q-(b_{1}-b_{r}+1)$ for $r=1,\cdots,\ell-1$ and $1\leq s_{\ell}\leq q-(b_{1}-b_{\ell}+1)$. 
Thus 
\begin{Prop}\label{Euclidean-2}
Let $\mathcal{S}$, $\mathcal{M}$ and $\{s_{k}\mid 1\leq k\leq\ell\}$ be as the stated above. Then the following subsets of objects are contained in $\mathcal{F}(\mathcal{S})$.
\begin{enumerate}[$(1)$]
 \item $\{(0,a_{1},b_{1}-j)\mid j=0,1,\cdots, q-s_{1}\}\cup\{(0,a_{2},b_{2}+j)\mid j=0,1,\cdots,b_{1}-b_{2}+s_{1}-q-1\}$,
\item $\{(0,a_{i},b_{i}-j)\mid j=0,1,\cdots,q+b_{i}-b_{1}-s_{i}\}\cup\{(0,a_{i+1},b_{i+1}+j)\mid j=0,1,\cdots,b_{1}-b_{i+1}+s_{i}-q-1\}$ for $i=2,3,\cdots,\ell-1$.
\item $\{(0,a_{\ell},b_{\ell}-j)\mid j=0,1,\cdots,q+b_{\ell}-b_{1}-s_{\ell}\}\cup\{(0,a_{1}+p,b_{1}-q+j)\mid j=0,1,\cdots,s_{i}-1\}.$
\end{enumerate}
\end{Prop}

We hope to generalize Propositions \ref{Euclidean-1} and   \ref{Euclidean-2} for a maximal orthogonal system on $\mathcal{S}\cap({_{s}{Q}_{0}\cup{_{s}{Q}_{1}}})$ or $\mathcal{S}\cap({_{s}{Q'}_{0}\cup{_{s}{Q'}_{1}}})$. 
We study properties of a maximal orthogonal  system  on triangle areas  $\bigtriangleup_{Q(0,j,i)'}\cup\bigtriangleup_{Q(1,j,i)'}$ for $i\leq p-2$. 
Let $\mathcal{W}$ be a maximal orthogonal  system on triangle areas  $\bigtriangleup_{Q(0,j,i)'}\cup\bigtriangleup_{Q(1,j,i)'}$. 
Without loss of generality, we assume 
$\mathcal{W}=\mathcal{U}\cup\mathcal{V}$,
where $\mathcal{U}=\{U_{1},U_{2},\cdots,U_{k}\}\subseteq\bigtriangleup_{Q(0,j,i)'}$, $\mathcal{V}=\{V_{1},V_{2},\cdots,V_{\ell}\}\subseteq\bigtriangleup_{Q(1,j,i)'}$ and  the second coordinates of objects of $\mathcal{U}$ and $\mathcal{V}$ are in an ascending way, respectively. Note that $k+\ell=i+1$ by Theorem \ref{cardinality-of-orth-sys}.
Then we have the following observations.
\begin{Prop}\label{triangle-areas-1}
Let $\mathcal{W}$, $\mathcal{U}$ and $\mathcal{V}$ be subsets on triangle areas $\bigtriangleup_{Q(0,j,i)'}\cup\bigtriangleup_{Q(1,j,i)'}$ as stated above. 
Then 
\begin{enumerate}[$(1)$]
\item  $V_{1}$ in $\mathcal{V}$ is of the form $Q(1,j,m)'$ or $U_{1}$ in $\mathcal{U}$ is of the form $Q(0,j,k)'$ for integers $m$ and $k$. 

\item If $V_{1}$ is of the  form $Q(1,j,m)'$ with the height $m>0$, then $U_{1}$ in $\mathcal{U}$ is of the form $Q(0,j,k)'$ with $k<m$.  

\item If $U_{1}$ is of the  form $Q(0,j,n)'$ with the height $n>0$, then there is an object in $\mathcal{V}$ of the form $Q(0,j+i_{0},n-i_{0})'$.

\item If there is an object $V_{r}=Q(1,j_{1},i_{1})'$ in $\mathcal{V}$ with the height $i_{1}>0$, then there is an object $U_{s}$ in $\mathcal{U}$ of the form $Q(0,j_{1},i'_{1})'$ with $i'_{1}<i_{1}$.

\item If there is an object $U_{e}=Q(0,j_{2},i_{2})'$ in $\mathcal{U}$ with the height $i_{2}>0$, then there is an object $V_{f}$  in $\mathcal{V}$ of the form $Q(1,j_{2}+t,i_{2}-t)'$.
\end{enumerate}
\end{Prop}
\begin{proof}
We prove only (1) and (2), since the proof of (3) is dual to (2) and (4) (resp. (5)) is similar to (2)(resp. (3)).

(1) If there is no object in the set $\{Q(0,j,t)'\mid 0\leq t\leq i\}$ contained in $\mathcal{U}$ and no object in $\{Q(1,j,t)'\mid 0\leq t\leq i\}$ contained in $\mathcal{V}$, then the cardinality of a maximal orthogonal system on triangle areas  $\bigtriangleup_{Q(0,j,i)'}\cup\bigtriangleup_{Q(1,j,i)'}$ equals to 	the cardinality of a maximal orthogonal system on triangle areas  $\bigtriangleup_{Q(0,j+1,i-1)'}\cup\bigtriangleup_{Q(1,j+1,i-1)'}$. It contradicts  Theorem \ref{cardinality-of-orth-sys}.

(2)	By Theorem \ref{cardinality-of-orth-sys}, the cardinality of a maximal orthogonal system on triangle areas  $\bigtriangleup_{Q(1,j,m)'}\cup\bigtriangleup_{Q(1,j,m)'}$ is $m+1$. By equations (\ref{equ-1}) and (\ref{equ-3}), ${^{\bot}Q(1,j,m)'^{\bot}}\cap(\bigtriangleup_{Q(0,j,m)'}\cup\bigtriangleup_{Q(1,j,m)'})=\{Q(0,j,t)'\mid 0\leq t\leq m-1\}\cup(\bigtriangleup_{Q(0,j+1,m-2)'}\cup\bigtriangleup_{Q(1,j+1,m-2)'})$.
By Theorem \ref{cardinality-of-orth-sys}, the cardinality of a maximal orthogonal system on triangle areas  $\bigtriangleup_{Q(0,j+1,m-2)'}\cup\bigtriangleup_{Q(1,j+1,m-2)'}$ is $m-1$.
Since $m>0$,  the set $\{Q(0,j,t)'\mid 0\leq t\leq m-1\}$  is not empty. 
Let $\mathcal{U}_{1}$ be a maximal orthogonal system on triangle areas  $\bigtriangleup_{Q(0,j,m)'}\cup\bigtriangleup_{Q(1,j,m)'}$ containing $Q(1,j,m)'$. 
If $\mathcal{U}_{1}$ contains no object in $\{Q(0,j,t)'\mid 0\leq t\leq m-1\}$, then the cardinality $\mathcal{U}_{1}$ is $1+(m-1)=m$. It contradicts Theorem \ref{cardinality-of-orth-sys}. Therefore there is an object $Q(0,j,k)'$ in $\{Q(0,j,t)'\mid 0\leq t\leq m-1\}$ contained in $\mathcal{V}$ and $k<m$. Thus $V_{1}=Q(0,j,k)'$.
\end{proof}

By equation (\ref{eq-12}), we know that $^{\bot}{\mathcal{M}}^{\bot}\cap(_{s}Q_{0}\cup{_{s}Q_{1}})$ is a disjoint union of triangle areas of some vertices.  Without loss of generality, we only consider a maximal orthogonal system on triangle areas $\bigtriangleup_{Q(0,1,r_{1})'}\cup\bigtriangleup_{Q(1,1,r_{1})'}$, where $r_{1}=a_{2}-a_{1}-2$. By Theorem \ref{cardinality-of-orth-sys}, the cardinality of a maximal orthogonal system on triangle areas  $\bigtriangleup_{Q(0,1,r_{1})'}\cup\bigtriangleup_{Q(1,1,r_{1})'}$ is $r_{1}+1$.
Take   $\mathcal{U}_{1}:=\{U_{1},\cdots,U_{s}\}$ and    $\mathcal{V}_{1}:=\{V_{1},\cdots,V_{t}\}$ such that $\mathcal{U}_{1}$ is a maximal subset of \, $\mathcal{U}$ contained in $\bigtriangleup_{Q(0,1,r_{1})'}$, where $s$ is a positive integer and $t=r_{1}+1-s$.  Note that 
$\mathcal{V}_{1}$ is a maximal subset of $\mathcal{V}$ contained in $\bigtriangleup_{Q(1,1,r_{1})'}$. 

\begin{Prop}\label{triangle-areas-2}
Let $\mathcal{U}_{1}$, $\mathcal{V}_{1}$ and related notations be as stated above. Then there are a subset  $\{Q(1,k_{m},\ell_{m})'\mid1\leq m\leq v_{1}\}$ of $\mathcal{V}_{1}$ and a subset $\{Q(0,k_{n},\ell_{n})'\mid v_{1}+1\leq n\leq v_{2}\}$ of  $\mathcal{U}_{1}$ for non-negative integers $v_{1}$ and $v_{2}$
such that 
$\mathcal{U}_{1}\cap\mathcal{V}_{1}$ is contained in a disjoint union of some triangle areas as follows.
\[\bigcup_{m=1}^{v_{1}}(\bigtriangleup_{Q(0,k_{m},\ell_{m})'}\cup\bigtriangleup_{Q(1,k_{m},\ell_{m})'})\cup\bigcup_{n=v_{1}+1}^{v_{2}}(\bigtriangleup_{Q(0,k_{n},\ell_{n})'}\cup\bigtriangleup_{Q(1,k_{n},\ell_{n})'}).\]	
\end{Prop}
\begin{proof}
{\bf Case one:} There is a maximal non-empty subset $\{U_{1},\cdots, U_{s_{1}}\}$ of $\mathcal{U}_{1}$ such that each object $U_{k}=Q(0,k,0)'$ is a quasi-simple for   $1\leq k\leq s_{1}$.

By equations (\ref{equ-1}) and (\ref{equ-3}), 
\begin{equation}\label{orth-sys-qus-1}
\begin{aligned}
^{\bot}&\{U_{1},\cdots, U_{s_{1}}\}^{\bot}\cap(\bigtriangleup_{Q(0,1,r_{1})'}\cup\bigtriangleup_{Q(1,1,r_{1})'})\\
=&\bigtriangleup_{Q(0,s_{1},r_{1}-s_{1}+1)'}\cup
\bigtriangleup_{Q(1,s_{1}+1,r_{1}-s_{1})'}\cup\{Q(1,1,s_{1}+k)'\mid -1\leq k\leq r_{1}-s_{1}\}.
\end{aligned}
\end{equation}

(1) If  the object $V_{1}$ is in the set $\{Q(1,1,s_{1}+k)'\mid -1\leq k\leq r_{1}-s_{1}\}$. We assume $V_{1}=Q(1,1,\ell)'$, where $s_{1}-1\leq\ell\leq r_{1}$.
By equations (\ref{eq-11}) and (\ref{eq-12}),
\begin{equation}\label{orth-sys-qus-2}
\begin{aligned}
^{\bot}&\{U_{1},\cdots, U_{s_{1}},V_{1}\}^{\bot}\cap(\bigtriangleup_{Q(0,1,r_{1})'}\cup\bigtriangleup_{Q(1,1,r_{1})'})\\
=&(\bigtriangleup_{Q(0,s_{1},\ell-s_{1})'}\cup\bigtriangleup_{Q(0,\ell+2,r_{1}-\ell-1)'})
\cup(\bigtriangleup_{Q(1,s_{1}+1,\ell-s_{1}-1)'}\cup\bigtriangleup_{Q(1,\ell+2,r_{1}-\ell-1)'})\\
=&(\bigtriangleup_{Q(0,s_{1},\ell-s_{1})'}\cup\bigtriangleup_{Q(1,s_{1}+1,\ell-s_{1}-1)'})\cup(\bigtriangleup_{Q(0,\ell+2,r_{1}-\ell-1)'}\cup\bigtriangleup_{Q(1,\ell+2,r_{1}-\ell-1)'}).
\end{aligned}
\end{equation}
Thus ${\mathcal{U}_{1}}\cup{\mathcal{V}_{1}}\subseteq(\bigtriangleup_{Q(0,1,\ell)'}\cup\bigtriangleup_{Q(1,1,\ell)'})\cup(\bigtriangleup_{Q(0,\ell+2,r_{1}-\ell-1)'}\cup\bigtriangleup_{Q(1,\ell+2,r_{1}-\ell-1)'})$ and $Q(1,1,\ell)'\in\mathcal{U}_{1}$. Furthermore, it is easy to see 
that $V_{2}\in\{Q(0,\ell+2,k)'\mid 0\leq k\leq r_{1}-\ell-1\}$ or $U_{s_{1}+1}\in\{Q(1,\ell+2,k)'\mid 0\leq k\leq r_{1}-\ell-1\}$ on $\bigtriangleup_{Q(0,\ell+2,r_{1}-\ell-1)'}\cup\bigtriangleup_{Q(1,\ell+2,r_{1}-\ell-1)'}$.
We may repeat Case one on triangle areas $\bigtriangleup_{Q(0,\ell+2,r_{1}-\ell-1)'}\cup\bigtriangleup_{Q(1,\ell+2,r_{1}-\ell-1)'}$.

(2) If  the object $V_{1}$ is not in the set $\{Q(1,1,s_{1}+k)'\mid -1\leq k\leq r_{1}-s_{1}\}$, then $U_{s_{i}+1}=Q(0,s_{1},\ell)'$  for $1\leq\ell\leq r_{1}-s_{1}+1$.
Thus 
\begin{equation}\label{orth-sys-qus-3}
\begin{aligned}
^{\bot}&(\mathcal{U}_{1}\cup\{U_{s_{1}+1}\})^{\bot}\cap(\bigtriangleup_{Q(0,1,r_{1})'}\cup\bigtriangleup_{Q(1,1,r_{1})'})\\
\subseteq&\bigcup_{k=1}^{s_{i}-1}(\bigtriangleup_{Q(0,k,0)'}\cup\bigtriangleup_{Q(1,k,0)'})\cup(\bigtriangleup_{Q(0,s_{1},\ell)'}\cup\bigtriangleup_{Q(1,s_{1},\ell)'})\\
&\cup(\bigtriangleup_{Q(0,s_{1}+\ell-1,r_{1}-s_{1}-\ell)'}\cup\bigtriangleup_{Q(1,s_{1}+\ell+2,r_{1}-s_{1}-\ell-1)'}).
\end{aligned}
\end{equation}
Hence
\begin{equation}
\begin{aligned} {\mathcal{U}_{1}}\cup{\mathcal{V}_{1}}\subseteq&\bigcup_{k=1}^{s_{i}-1}(\bigtriangleup_{Q(0,k,0)'}\cup\bigtriangleup_{Q(1,k,0)'})\cup(\bigtriangleup_{Q(0,s_{1},\ell)'}\cup\bigtriangleup_{Q(1,s_{1},\ell)'})\\
&\cup(\bigtriangleup_{Q(0,s_{1}+\ell-1,r_{1}-s_{1}-\ell)'}\cup\bigtriangleup_{Q(1,s_{1}+\ell+2,r_{1}-s_{1}-\ell-1)'}),
\end{aligned}
\end{equation}
where $Q(0,k,0)'\in\mathcal{U}_{1}$ for $1\leq k\leq s_{1}-1$ and $U_{s_{1}+1}=Q(0,s_{1},\ell)'\in\mathcal{U}_{1}$. Furthermore, it is easy to see 
that $U_{s_{1}+2}\in\{Q(0,s_{1}+\ell+1,k)'\mid 0\leq k\leq r_{1}-\ell-1\}$. Then we may repeat (2) of Case one  on triangle areas  $\bigtriangleup_{Q(0,s_{1}+\ell-1,r_{1}-s_{1}-\ell)'}\cup\bigtriangleup_{Q(1,s_{1}+\ell+2,r_{1}-s_{1}-\ell-1)'}$. 

\medskip

\noindent{\bf Case two}: There is a maximal non-empty subset $\{V_{1},\cdots, V_{t_{1}}\}$ of $\mathcal{V}_{1}$ such that each $V_{k}=Q(1,k,0)'$ is a quasi-simple for $k=1,\cdots,t_{1}$.

By equations (\ref{equ-1}) and (\ref{equ-3}), 
\begin{equation}\label{orth-sys-qus-4}
\begin{aligned}
^{\bot}&\{V_{1},\cdots, V_{t_{1}}\}^{\bot}\cap(\bigtriangleup_{Q(0,1,r_{1})'}\cup\bigtriangleup_{Q(1,1,r_{1})'})\\
=&\bigtriangleup_{Q(0,t_{1},r_{1}-t_{1}+1)'}\cup
\bigtriangleup_{Q(1,t_{1},r_{1}-t_{1}+1)'}\cup\{Q(0,1,t_{1}+k)'\mid -2\leq k\leq r_{1}-t_{1}\}.
\end{aligned}
\end{equation}

(1) If  the object $U_{1}$ is in the set $\{Q(0,1,t_{1}+k)'\mid -2\leq k\leq r_{1}-t_{1}\}$. We assume $U_{1}=Q(0,1,\ell)'$, where $t_{1}-2\leq\ell\leq r_{1}$.
By equations (\ref{eq-11}) and (\ref{eq-12}),
\begin{equation}\label{orth-sys-qus-5}
\begin{aligned}
^{\bot}&\{V_{1},\cdots, V_{t_{1}},U_{1}\}^{\bot}\cap(\bigtriangleup_{Q(0,1,r_{1})'}\cup\bigtriangleup_{Q(1,1,r_{1})'})\\
=&(\bigtriangleup_{Q(0,t_{1},\ell-t_{1})'}\cup\bigtriangleup_{Q(0,\ell+2,r_{1}-\ell-1)'})
\cup(\bigtriangleup_{Q(1,t_{1},\ell-t_{1}+1)'}\cup\bigtriangleup_{Q(1,\ell+3,r_{1}-\ell-2)'})\\
=&(\bigtriangleup_{Q(0,t_{1},\ell-t_{1})'}\cup\bigtriangleup_{Q(1,t_{1},\ell-t_{1}+1)'})\cup(\bigtriangleup_{Q(0,\ell+2,r_{1}-\ell-1)'}\cup\bigtriangleup_{Q(1,\ell+3,r_{1}-\ell-2)'}).
\end{aligned}
\end{equation}
Therefore ${\mathcal{U}_{1}}\cup{\mathcal{V}_{1}}\subseteq(\bigtriangleup_{Q(0,1,\ell)'}\cup\bigtriangleup_{Q(1,1,\ell)'})\cup(\bigtriangleup_{Q(0,\ell+2,r_{1}-\ell-1)'}\cup\bigtriangleup_{Q(1,\ell+3,r_{1}-\ell-2)'})$, $V_{k}=Q(1,k,0)'\in\mathcal{V}_{1}$ for $1\leq k\leq t_{1}-1$ and $U_{1}=Q(0,1,\ell)'\in\mathcal{U}_{1}$. Furthermore, there is an object  $U_{s'}$ in $\mathcal{U}_{1}$ contained in the set $\{Q(0,\ell+2,k)'\mid 0\leq k\leq r_{1}-\ell-1\}$.
Thus we may repeat (1) of Case one  on $\bigtriangleup_{Q(0,\ell+2,r_{1}-\ell-1)'}\cup\bigtriangleup_{Q(1,\ell+3,r_{1}-\ell-2)'}$.

(2) If the object  $U_{1}$ is not in the set $\{Q(0,1,t_{1}+k)'\mid -2\leq k\leq r_{1}-t_{1}\}$, then $U_{t_{i}+1}=Q(0,t_{1},\ell)'$  for $1\leq\ell\leq r_{1}-t_{1}+1$ or $V_{1}=Q(1,t_{1},\ell')'$  for $1\leq\ell'\leq r_{1}-t_{1}+1$.
Thus 
\begin{equation}\label{orth-sys-qus-6}
\begin{aligned}
^{\bot}&\mathcal{V}_{1}^{\bot}\cap(\bigtriangleup_{Q(0,1,r_{1})'}\cup\bigtriangleup_{Q(1,1,r_{1})'})\\
\subseteq&\bigcup_{k=1}^{t_{i}-1}(\bigtriangleup_{Q(0,k,0)'}\cup\bigtriangleup_{Q(1,k,0)'})
\cup(\bigtriangleup_{Q(0,t_{1},r_{1}-\ell-1)'}\cup\bigtriangleup_{Q(1,t_{1},r_{1}-\ell-1)'}).
\end{aligned}
\end{equation}
Therefore \[{\mathcal{U}_{1}}\cup{\mathcal{V}_{1}}\subseteq\bigcup_{k=1}^{t_{i}-1}(\bigtriangleup_{Q(0,k,0)'}\cup\bigtriangleup_{Q(1,k,0)'})\cup(\bigtriangleup_{Q(0,t_{1},r_{1}-\ell-1)'}\cup\bigtriangleup_{Q(1,t_{1},r_{1}-\ell-1)'}),\] $V_{k}=Q(1,k,0)'\in\mathcal{V}_{1}$ for $1\leq k\leq t_{1}-1$. Furthermore  $U_{1}\in\{Q(0,\ell+2,k)'\mid 0\leq k\leq r_{1}-\ell-1\}$ or $V_{1}\in\{Q(1,\ell+2,k)'\mid 0\leq k\leq r_{1}-\ell-1\}$.
Thus we may repeat Case one or Case two on $\bigtriangleup_{Q(0,t_{1},r_{1}-\ell-1)'}\cup\bigtriangleup_{Q(1,t_{1},r_{1}-\ell-1)'}.$

Summarizing the above two cases, it follows from induction that the set $\mathcal{U}_{1}\cap\mathcal{V}_{1}$ may be divided into a  disjoint union of some subsets as follows.
\[\bigcup_{m=1}^{v_{1}}(\bigtriangleup_{Q(0,k_{m},\ell_{m})'}\cup\bigtriangleup_{Q(1,k_{m},\ell_{m})'})\cup\bigcup_{m=v_{1}+1}^{v_{2}}(\bigtriangleup_{Q(0,k_{m},\ell_{m})'}\cup\bigtriangleup_{Q(1,k_{m},\ell_{m})'}),\]
 where $Q(x,k_{m},\ell_{m})'\notin\bigtriangleup_{Q(0,k_{n},\ell_{n})'}\cup\bigtriangleup_{Q(1,k_{n},\ell_{n})'}$ for any $x=0,1$ and $m\neq n$.
 
For  (1) of Case one and  (2) of Case two, we have 
$Q(1,k_{m},\ell_{m})'\in\mathcal{V}_{1}$ for $1\leq m\leq v_{1}$ and $Q(0,k_{m},\ell_{m})'\in\mathcal{U}_{1}$ for $v_{1}+1\leq m\leq v_{2}$.
For (2) of Case one and (1) if Case two, we have 
$Q(0,k_{m},\ell_{m})'\in\mathcal{V}_{1}$ for  $1\leq m\leq v_{2}$. Thus the  proof is finished.
\end{proof}

Now we consider the triangle areas $\bigtriangleup_{Q(0,k_{m},\ell_{m})'}\cup\bigtriangleup_{Q(1,k_{m},\ell_{m})'}$ with $Q(0,k_{m},\ell_{m})'\in\mathcal{U}_{1}$ and the triangle areas $\bigtriangleup_{Q(0,k_{n},\ell_{n})'}\cup\bigtriangleup_{Q(1,k_{n},\ell_{n})'}$ with $Q(1,k_{n},\ell_{n})'\in\mathcal{V}_{1}$. Without loss of generality, we simplify the above notations and prove the following  Theorem.

\begin{Them}\label{Euclidean-Quasi-tube}
Let $\mathcal{S}$ be an orthogonal system on $A$-$\stmod$ and the subset $\mathcal{X}$ of $\mathcal{S}$  an orthogonal system on triangle areas $\bigtriangleup_{Q(0,m,\ell)'}\cup\bigtriangleup_{Q(1,m,\ell)'}$ with $Q(0,m,\ell)'\in\mathcal{X}$ for $\ell\leq q-2$. If there is an object $(0,c,d)$ contained in  $\mathcal{F}(\mathcal{S})$ such that there is a triangle as follows.
\[Q(0,m+\ell,0)'\xrightarrow{}(0,c-1,d)\xrightarrow{} (0,c,d)\xrightarrow{}
Q(1,m+\ell+1,0)',\]
then the subset $\{(0,c-k,d)\mid 0\leq k\leq\ell\}\subseteq\mathcal{F}(\mathcal{S})$.
\end{Them}
Note that the triangle in Theorem \ref{Euclidean-Quasi-tube} is  from Proposition \ref{sms-quasi-simple}. Dually, 
\begin{Them}\label{Euclidean-Quasi-tube'}
Let $\mathcal{S}$ be an orthogonal system on $A$-$\stmod$ and the subset $\mathcal{X}$ of  $\mathcal{S}$ an orthogonal system on triangle areas $\bigtriangleup_{Q(0,r,s)'}\cup\bigtriangleup_{Q(1,r,s)'}$ with $Q(1,r,s)'\in\mathcal{X}$ for $s\leq q-2$. If there is an object $(0,e,f)$ contained in  $\mathcal{F}(\mathcal{S})$ such that there is a triangle as follows.
\[(0,e,f)\xrightarrow{} (0,e+1,f)\xrightarrow{}Q(1,r,0)'\xrightarrow{}
(1,e+1,f+1),\]
then the subset $\{(0,e+t,f)\mid 0\leq t\leq r\}\subseteq\mathcal{F}(\mathcal{S})$.
\end{Them}

Before proving the above Theorem \ref{Euclidean-Quasi-tube}, we state the following key proposition.
\begin{Prop}\label{Euc-quasi-1} 
Let $\mathcal{M}$ be an orthogonal system on $A$-$\stmod$ and $\mathcal{X}\subseteq\mathcal{M}$ be a maximal orthogonal system on triangle areas $\bigtriangleup_{Q(0,x,y-1)'}\cup\bigtriangleup_{Q(1,x,y)'}$. We assume that  there is a triangle 
\[(0,b-1,c)\xrightarrow{}(0,b,c)\xrightarrow{} Q(1,x,0)'\xrightarrow{} (1,b,c+1).\]
If $(0,b-1,c)\in\mathcal{F}(\mathcal{M})$, then $(0,b+k,c)\in\mathcal{F}(\mathcal{M})$ for integers $k=0,1,\cdots,y.$

Dually, we assume that  $\mathcal{Z}\subseteq\mathcal{M}$ is a maximal orthogonal system on triangle areas  $\bigtriangleup_{Q(0,v,w)'}\cup\bigtriangleup_{Q(1,v+1,w-1)'}$ and that  there is a  triangle 
\[Q(0,v+w,0)'\xrightarrow{}(0,e,f)\xrightarrow{} (0,e+1,f)\xrightarrow{} Q(1,v+w+1,0)'.\] 
If $(0,e+1,f)\in\mathcal{F}(\mathcal{M})$, then $(0,e-k,f)\in\mathcal{F}(\mathcal{M})$ for integers $k=0,1,\cdots,w.$
\end{Prop}
\begin{proof}
We prove it by induction on the height $y$ and $w$. 
We first consider the case $y=1$ and $w=1$. Thus $\bigtriangleup_{Q(0,x,0)'}=\{Q(0,x,0)'\}$ and $\bigtriangleup_{Q(1,x,1)'}=\{Q(1,x,0)',Q(1,x,1)',Q(1,x+1,0)'\}$.
By equations (\ref{eq-11}) and (\ref{eq-12}), there are  two  maximal orthogonal systems on triangle areas  $\bigtriangleup_{Q(0,x,0)'}\cup\bigtriangleup_{Q(1,x,1)'}$ as follows.
\[\mathcal{X}_{1}=\{Q(1,x,0)',Q(1,x+1,0)'\},\ \ \mathcal{X}_{2}=\{Q(1,x,1)',Q(0,x,0)'\}.\]

\noindent{\bf Case one}: $\mathcal{X}=\mathcal{X}_{1}=\{Q(1,x,0)',Q(1,x+1,0)'\}$.

Consider the triangles as follows.
\[(0,b-1,c)\xrightarrow{}(0,b,c)\xrightarrow{} Q(1,x,0)'\xrightarrow{} (1,b,c+1),\]
\[(0,b,c)\xrightarrow{}(0,b+1,c)\xrightarrow{} Q(1,x+1,0)'\xrightarrow{} (1,b+1,c+1).\]
Since $(0,b-1,c)\in\mathcal{F}(\mathcal{M})$ and $Q(1,x,0)'\in\mathcal{X}_{1}\subseteq\mathcal{M}$, $(0,a,b)\in\mathcal{F}(\mathcal{M})$. Since $(0,a,b)\in\mathcal{F}(\mathcal{M})$ and  $Q(1,x+1,0)'\in\mathcal{X}_{1}\subseteq\mathcal{M}$, $(0,b+1,c)\in\mathcal{F}(\mathcal{M}).$

\noindent{\bf Case two}: $\mathcal{X}=\mathcal{X}_{2}=\{Q(1,x,1)',Q(0,x,0)'\}$.

Consider the triangle as follows.
\[(0,b-1,c)\xrightarrow{} (0,b+1,c)\xrightarrow{} Q(1,x,1)'\xrightarrow{} (1,b,c+1),\]
\[Q(0,x,0)'\xrightarrow{}(0,b,c)\xrightarrow{} (0,b+1,c)\xrightarrow{} Q(1,x+1,0)'.\]
Since  $(0,b-1,c)\in\mathcal{F}(\mathcal{M})$ and $Q(1,x,1)'\in\mathcal{X}_{2}\subseteq\mathcal{F}(\mathcal{M})$,  $(0,b+k,c)\in\mathcal{F}(\mathcal{M})$ for $k=0,1$. Thus our conclusion holds for $y=1$.  

Dually,  we show the conclusion is true for $w=1$.
Thus $\bigtriangleup_{Q(0,v,1)'}=\{Q(0,v,0)',Q(0,v,1)',Q(0,v+1,0)'\}$ and $\bigtriangleup_{Q(1,v+1,0)'}=\{Q(1,v+1,0)'\}.$
By equations (\ref{eq-11}) and (\ref{eq-12}), there are two  maximal orthogonal systems on triangle areas  $\bigtriangleup_{Q(0,v,1)'}\cup\bigtriangleup_{Q(1,v+1,0)'}$ as follows.
\[\mathcal{Z}_{1}=\{Q(0,v,0)',Q(0,v+1,0)'\},\ \ \mathcal{Z}_{2}=\{Q(0,v,1)',Q(1,v+1,0)'\}.\]

\noindent{\bf Case one}: $\mathcal{Z}=\mathcal{Z}_{1}=\{Q(0,v,1)',Q(0,v+1,0)'\}$. 

Consider the triangles as follows.
\[Q(0,v+1,0)'\xrightarrow{}(0,e,f)\xrightarrow{} (0,e+1,f)\xrightarrow{} Q(1,v+2,0)'.\]
\[Q(0,v,0)'\xrightarrow{}(0,e-1,f)\xrightarrow{} (0,e,f)\xrightarrow{} Q(1,v+1,0)'.\]
Since $(0,e+1,f)\in\mathcal{F}(\mathcal{M})$ and $Q(0,v+1,0)'\in\mathcal{Z}_{1}\subseteq\mathcal{M}$, $(0,e,f)\in\mathcal{F}(\mathcal{M})$. Since $Q(0,v,0)'\in\mathcal{Z}_{1}$ and $(0,e,f)\in\mathcal{F}(\mathcal{M})$, $(0,e-1,f)\in\mathcal{F}(\mathcal{M}).$

\noindent{\bf Case two}: $\mathcal{Z}=\mathcal{Z}_{2}=\{Q(0,v,1)',Q(1,v+1,0)'\}$.

Consider the triangles as follows.
\[Q(0,v,1)'\xrightarrow{}(0,e-1,f)\xrightarrow{} (0,e+1,f)\xrightarrow{} Q(1,v+1,1)',\]  
\[(0,e-1,f)\xrightarrow{}(0,e,f)\xrightarrow{} Q(1,v+1,0)'\xrightarrow{} (1,e,f+1).\] 
Since $Q(0,v,1)'\in\mathcal{Z}_{2}\subseteq\mathcal{M}$ and $(0,e+1,f)$ is in $\mathcal{F}(\mathcal{M})$, $(0,e-k,f)\in\mathcal{F}(\mathcal{M})$ for $k=0,1$. Thus our conclusion holds for $w=1$. 

We assume our conclusion is true for $y=k$ and $w=k$, and we only prove it is true for $y=k+1$, since the other case is dual. 
By Theorem \ref{cardinality-of-orth-sys}, there is an object $Q(1,x_{1},y_{1})'$ on the set $\{Q(1,x+s,k+1-s)'\mid 0\leq s\leq k+1\}$ belonging to $\mathcal{X}$, where $x_{1}+y_{1}=k+1+x$ and $y_{1}\leq k+1$.
By  equations (\ref{equ-1}) and (\ref{equ-3}), 
\begin{equation}
\begin{aligned}
&(\bigtriangleup_{Q(0,x,k)'}\cup\bigtriangleup_{Q(1,x,k+1)'})\cap\mathcal{X}\\
&=((\bigtriangleup_{Q(0,x_{1},y_{1}-1)'}\cup\bigtriangleup_{Q(1,x_{1}+1,y_{1}-2)'})\cap\mathcal{X})\cup((\bigtriangleup_{Q(0,x,x_{1}-x-1)'}\cup\bigtriangleup_{Q(1,x,x_{1}-x)'})\cap\mathcal{X})
\end{aligned}
\end{equation}
Since  $x_{1}-x<k$ and $(0,b-1,c)\in\mathcal{F}(\mathcal{M})$,  $(0,b+k,c)\in\mathcal{F}(\mathcal{M})$ for $k=0,1,\cdots,x_{1}-x$ by induction. 

Consider the triangle as follows.
\[(0,b+x_{1}-x,0)\xrightarrow{}(0,b+k+1,0)\xrightarrow{} Q(1,x_{1},y_{1})'\xrightarrow{} (1,b-x_{1}-x+1,1).\] 
Since $(0,b-x_{1}-x,0)\in\mathcal{F}(\mathcal{M})$ and $Q(1,x_{1},y_{1})'\in\mathcal{X}\subseteq\mathcal{M}$, $(0,b+k+1,0)\in\mathcal{F}(\mathcal{M})$.
Since $y_{1}-1\leq y$ and $(0,b+x_{1}+y_{1},0)\in\mathcal{F}(\mathcal{M})$, $(0,b+k+1-s,0)\in\mathcal{F}(\mathcal{M})$ for $s=0,\cdots,y_{1}$ by induction. 
Combining $(0,b+k,c)\in\mathcal{F}(\mathcal{M})$ for $k=0,1,\cdots,x_{1}-x$ and $(0,b+k+1-s,0)\in\mathcal{F}(\mathcal{M})$ for $s=0,\cdots,y_{1}$, we have  $(0,b+r,c)\in\mathcal{F}(\mathcal{M})$ for $r=0,\cdots,k+1$. Note that $b+k+1-y_{1}=b+x_{1}-x$.
\end{proof}

\noindent{\bf Proof of Theorem \ref{Euclidean-Quasi-tube}:}
\begin{proof}
Note that there is a non-zero morphism from $Q(0,m,0)'$ to $(0,c-\ell,d)$ by Lemma \ref{irreducible-quasi-simples}.
	
Consider  a triangle as follows.
\begin{align}\label{seqs-E-1}
Q(0,m,\ell)'\xrightarrow{}(0,c-\ell,0)\xrightarrow{} (0,c,d)\xrightarrow{} Q(1,m+1,\ell)'.
\end{align}
Since $(0,c,d)\in\mathcal{F}(\mathcal{S})$ and $Q(0,m,\ell)'\in\mathcal{X}\subseteq\mathcal{S}$, $(0,c-\ell,0)\in\mathcal{F}(\mathcal{S}).$
By  (3) of Proposition \ref{triangle-areas-1}, there is an object $Q(1,m_{1},\ell_{1})'$ in $\mathcal{X}$ such that $m_{1}+\ell_{1}=m+\ell.$
Thus by equations (\ref{equ-1}) and (\ref{equ-3}), 
\begin{equation}\label{tri-area-1}
\begin{aligned}
(\bigtriangleup_{Q(0,m,\ell)'}\cup\bigtriangleup_{Q(1,m,\ell)'})\cap\mathcal{X}
=\{Q(0,m,\ell)'\}\cup(\bigtriangleup_{Q(0,m+1,\ell-2)'}\cup\bigtriangleup_{Q(1,m,\ell-1)'})\cap\mathcal{X}.
\end{aligned}
\end{equation}
Considering an orthogonal system on triangle areas $\bigtriangleup_{Q(0,m+1,\ell-2)'}\cup\bigtriangleup_{Q(1,m,\ell-1)'}$, we have $\{(0,c-k,d)\mid 0\leq k\leq\ell\}\subseteq\mathcal{F}(\mathcal{S})$ by Proposition \ref{Euc-quasi-1}.
\end{proof}

\begin{Them}\label{ortho-Euclidean}
Let $A$ be a 2-domestic Brauer graph algebra and $\mathcal{S}$ a maximal orthogonal system  which contains at least one object for an Euclidean component. We assume $\mathcal{M}=\{M_{1},\cdots,M_{\ell}\}$ is the subset of $\mathcal{S}$ which is contained in the stable Euclidean component $_{s}\Gamma_{0}$, where $M_{k}=(0,a_{k},b_{k})$ for $1\leq k\leq\ell$. 
Then 
\begin{enumerate}[$(a)$]
\item There is a set $\{t_{i}\mid 1\leq i\leq \ell\}$  of integers such that the following subsets are contained in $\mathcal{F}(\mathcal{S})$.
\begin{enumerate}[$(1)$]
\item $\{(0,a_{1}+j,b_{1})\mid j=0,1,\cdots, t_{1}-1\}\cup\{(0,a_{2}-j,b_{2})\mid j=0,1,\cdots,a_{2}-a_{1}-t_{1}\}$.
	
\item $\{(0,a_{r}+j,b_{r})\mid j=0,1,\cdots,t_{r}-a_{r}+a_{1}-1\}\cup\{(0,a_{r+1}-j,b_{r+1})\mid j=0,1,\cdots,a_{r+1}-a_{1}-t_{r}\}$ for $r=2,3,\cdots,\ell-1$.
	
\item $\{(0,a_{\ell}+j,b_{\ell})\mid j=0,1,\cdots,t_{\ell}-a_{\ell}+a_{1}-1\}\cup\{(0,a_{1}+p-j,b_{1}-q)\mid j=0,1,\cdots,p-t_{\ell}\}.$
\end{enumerate}
\item Dually, there is a set $\{s_{i}\mid 1\leq i\leq \ell\}$  of integers such that the following subsets are contained in $\mathcal{F}(\mathcal{S})$.
\begin{enumerate}[$(1)$]
\item $\{(0,a_{1},b_{1}-j)\mid j=0,1,\cdots,  q-s_{1}\}\cup\{(0,a_{2},b_{2}+j)\mid j=0,1,\cdots,b_{1}-b_{2}+s_{1}-q-1\}$,
	
\item $\{(0,a_{i},b_{i}-j)\mid j=0,1,\cdots,q+b_{i}-b_{1}-s_{i}\}\cup\{(0,a_{i+1},b_{i+1}+j)\mid j=0,1,\cdots,b_{1}-b_{i+1}+s_{i}-q-1\}$ for $i=2,3,\cdots,\ell-1$.
	
\item $\{(0,a_{\ell},b_{\ell}-j)\mid j=0,1,\cdots,q+b_{\ell}-b_{1}-s_{\ell}\}\cup\{(0,a_{1}+p,b_{1}-q+j)\mid j=0,1,\cdots,s_{\ell}-1\}.$
\end{enumerate}
\item $\mathcal{S}\cap{_{s}\Gamma_{1}}=\{M'_{k}=(1,a_{1}+t_{k},b_{1}-s_{k}-q)\mid 1\leq k\leq\ell\}$. 
\end{enumerate}
\end{Them}
\begin{proof}
We only prove (1) of (a), (1) of (b) and $(1,a_{1}+t_{1},b_{1}-s_{1}-q)\in\mathcal{S}$ in (c), since other cases are similar.

\noindent(1) of (a) and (1) of (b):

By Proposition \ref{triangle-areas-2},   there are a subset  $\{Q(1,k_{m},\ell_{m})'\mid1\leq m\leq v_{1}\}$ of $_{s}Q'_{1}\cap\mathcal{S}$ and a subset $\{Q(0,k_{n},\ell_{n})'\mid v_{1}+1\leq n\leq v_{2}\}$ of  $_{s}Q'_{0}\cap\mathcal{S}$ for non-negative integers $v_{1}$ and $v_{2}$
such that 
$\mathcal{S}\cap(\triangle_{Q(0,1,a_{2}-a_{1}-2)'}\cup\triangle_{Q(1,1,a_{2}-a_{1}-2)'})$ is contained in a disjoint union of some triangle areas as follows.
\[\bigcup_{m=1}^{v_{1}}(\bigtriangleup_{Q(0,k_{m},\ell_{m})'}\cup\bigtriangleup_{Q(1,k_{m},\ell_{m})'})\cup\bigcup_{n=v_{1}+1}^{v_{2}}(\bigtriangleup_{Q(0,k_{n},\ell_{n})'}\cup\bigtriangleup_{Q(1,k_{n},\ell_{n})'}).\]	
Consider the  triangle areas $\triangle_{Q(0,k_{v_{1}},\ell_{v_{1}})'}\cup\triangle_{Q(1,k_{v_{1}},\ell_{v_{1}})'}$ and $\triangle_{Q(0,k_{v_{1}+1},\ell_{v_{1}+1})'}\cup\triangle_{Q(1,k_{v_{1}+1},\ell_{v_{1}+1})'}$, where $Q(1,k_{v_{1}},\ell_{v_{1}})'\in\mathcal{S}$ and 
$Q(0,k_{v_{1}+1},\ell_{v_{1}+1})'\in\mathcal{S}.$
By the orthogonality of objects on quasi-tubes, $k_{v_{1}}+\ell_{v_{1}}+1=k_{v_{1}+1}$. Take $t_{1}=k_{v_{1}+1}$. It is not hard to see that   $\mathcal{W}_{Q(0,t_{1}-1,0)'}\cap\mathcal{S}=\varnothing$ and $\mathcal{W}_{Q(1,t_{1},0)'}\cap\mathcal{S}=\varnothing$. Note that  $Q(0,t_{1}-1,0)'=\Omega(Q(1,t_{1},0)')$ and  $\mathcal{W}_{Q(0,t_{1}-1,0)'}$ is the wing of object $Q(0,t_{1}-1,0)'$.

By Lemma \ref{irreducible-quasi-simples}, there are two irreducible maps: \[Q(1,t_{1},0)'\rightarrow(1,t_{1}-a_{1},b_{1}),\   (1,t_{1}-a_{1},b_{1})\rightarrow Q(0,t_{1}-1,0)'.\]
Since $\mathcal{W}_{Q(0,t_{1}-1,0)'}\cap\mathcal{S}=\varnothing$ and $\mathcal{W}_{Q(1,t_{1},0)'}\cap\mathcal{S}=\varnothing$, by equation (\ref{equ-7}),   
\begin{equation}\label{gamma1-1''}
\{(1,t_{1}-a_{1},k)\mid k\in\mathbb{Z}\}\subseteq{^{\bot}(\mathcal{S}\cap(_{s}Q'_{0}\cup{_{s}Q'_{1}}))^{\bot}}.
\end{equation}
Since the objects $(0,a_{1},b_{1})$ and $(0,a_{2},b_{2})$ are in $\mathcal{F}(\mathcal{S})$, by Theorem \ref{Euclidean-Quasi-tube} and  Theorem \ref{Euclidean-Quasi-tube'}, 
\[\{(0,a_{1}+j,b_{1})\mid j=0,1,\cdots, t_{1}-1\}\cup\{(0,a_{2}-j,b_{2})\mid j=0,1,\cdots,a_{2}-a_{1}-t_{1}\} \subseteq\mathcal{F}(\mathcal{S}).\] 

Similarly, by Proposition  \ref{triangle-areas-2}, there are a subset  $\{Q(1,x_{m},y_{m})\mid1\leq m\leq u_{1}\}$ of $_{s}Q_{1}\cap\mathcal{S}$ and a subset $\{Q(0,x_{n},y_{n})\mid u_{1}+1\leq n\leq u_{2}\}$ of  $_{s}Q_{0}\cap\mathcal{S}$ for non-negative integers $u_{1}$ and $u_{2}$
such that 
$\mathcal{S}\cap(\triangle_{Q(0,b_{2}-b_{1}+p,b_{1}-b_{2}-1)}\cup\triangle_{Q(1,b_{2}-b_{1}+p,b_{1}-b_{2}-1)})$ is contained in a disjoint union of some triangle areas as follows.
\[\bigcup_{m=1}^{u_{1}}(\bigtriangleup_{Q(0,x_{m},y_{m})}\cup\bigtriangleup_{Q(1,x_{m},y_{m})})\cup\bigcup_{n=u_{1}+1}^{u_{2}}(\bigtriangleup_{Q(0,x_{n},y_{n})}\cup\bigtriangleup_{Q(1,x_{n},y_{n})}).\]	
Consider the  triangle areas $\triangle_{Q(0,x_{u_{1}},y_{u_{1}})}\cup\triangle_{Q(1,x_{u_{1}},y_{u_{1}})}$ and $\triangle_{Q(0,x_{u_{1}+1},y_{u_{1}+1})}\cup\triangle_{Q(1,x_{u_{1}+1},y_{u_{1}+1})}$, where $Q(1,x_{u_{1}},y_{u_{1}})\in\mathcal{S}$ and 
$Q(0,x_{u_{1}+1},y_{u_{1}+1})\in\mathcal{S}.$
By the orthogonality of objects on quasi-tubes, $x_{u_{1}}+y_{u_{1}}+1=x_{u_{1}+1}$. Take $s_{1}=x_{u_{1}+1}$. It is not hard to see that  $\mathcal{W}_{Q(0,s_{1}-1,0)}\cap\mathcal{S}=\varnothing$ and $\mathcal{W}_{Q(1,s_{1},0)}\cap\mathcal{S}=\varnothing$ and $Q(0,s_{1}-1,0)=\Omega(Q(1,s_{1},0))$. 

By Lemma \ref{irreducible-quasi-simples}, there are two irreducible maps: \[Q(1,s_{1},0)\rightarrow(1,a_{1},b_{1}+s_{1}-q), \  (1,a_{1},b_{1}+s_{1}-q)\rightarrow Q(0,s_{1}-1,0). \]
By equation (\ref{equ-7}), 
\begin{equation}\label{gamma1-2''}
\{(1,k,b_{1}+s_{1}-q)\mid k\in\mathbb{Z} \}\subseteq{^{\bot}(\mathcal{S}\cap(_{s}Q_{0}\cup{_{s}Q_{1}}))^{\bot}}.
\end{equation}
Since $(0,a_{1},b_{1})$ and $(0,a_{2},b_{2})$ are in $\mathcal{F}(\mathcal{S})$, by Theorem \ref{Euclidean-Quasi-tube}  and Theorem \ref{Euclidean-Quasi-tube'},  \[\{(0,a_{1},b_{1}-j)\mid j=0,1,\cdots, q-s_{1}\}\cup\{(0,a_{2},b_{2}+j)\mid j=0,1,\cdots,b_{1}-b_{2}+s_{1}-q-1\}\subseteq\mathcal{F}(\mathcal{S}).\]

\medskip
 
\noindent(c) Since  $\{(0,a_{1}+j,b_{1})\mid j=0,1,\cdots, t_{1}-1\}\cup\{(0,a_{2}-j,b_{2})\mid j=0,1,\cdots,a_{2}-a_{1}-t_{1}\}\subseteq\mathcal{F}(\mathcal{S})$,
\begin{equation}
\begin{aligned}
\{(1,a_{1}+j,b_{1})\mid j=0,1,\cdots,& t_{1}-1\}\cup\{(1,a_{2}-j,b_{2})\mid j=0,1,\cdots,a_{2}-a_{1}-t_{1}\}\\
&\nsubseteq{^{\bot}(\mathcal{S}\cap(_{s}\Gamma_{0}\cup {_{s}Q_{0}}\cup {_{s}Q_{1}}))^{\bot}}.
\end{aligned}
\end{equation}
Thus $(1,a_{1}+j,b_{1})$ and 
$(1,a_{2}-k,b_{2})$ are not in $\mathcal{S}$ for any $1\leq j\leq t_{1}-1$ and 
$1\leq k\leq a_{2}-a_{1}-t_{1}$, respectively.
Since $\{(0,a_{1},b_{1}-j)\mid j=0,1,\cdots,  q-s_{1}\}\cup\{(0,a_{2},b_{2}+j)\mid j=0,1,\cdots,b_{1}-b_{2}+s_{1}-q-1\}\subseteq\mathcal{F}(\mathcal{S})$, 
\begin{equation}
\begin{aligned}
\{(1,a_{1},b_{1}-j)\mid j=0,1,\cdots, & q-s_{1}\}\cup\{(1,a_{2},b_{2}+j)\mid j=0,1,\cdots,b_{1}-b_{2}+s_{1}-q-1\}\\
&\nsubseteq{^{\bot}(\mathcal{S}\cap(_{s}\Gamma_{0}\cup {_{s}Q_{0}}\cup{_{s}Q_{1}}))^{\bot}}.
\end{aligned}
\end{equation} 
Thus the objects  $(1,a_{1},b_{1}-j)$ and 
$(1,a_{2},b_{2}+k)$ are not in $\mathcal{S}$ for any $1\leq j\leq q-s_{1}$ and 
$1\leq k\leq b_{1}-b_{2}+s_{1}-q-1$, respectively.

By equations (\ref{gamma1-1''}) and (\ref{gamma1-2''}),
\begin{equation}\label{gamma1-1'}
\{(1,a_{1}+t_{1},k)\mid b_{2}+1\leq k\leq b_{1}\}={^{\bot}(\mathcal{S}\cap(_{s}\Gamma_{0}\cup{_{s}Q'_{0}}\cup {_{s}Q'_{1}}))^{\bot}}\cap\square_{(1,a_{1}+1,b_{1}),(1,a_{2},b_{2}+1)}, 
\end{equation} 
\begin{equation}\label{gamma1-2'}
\{(1,k,b_{1}+s_{1}-q)\mid a_{1}+1\leq k\leq a_{2}\}={^{\bot}(\mathcal{S}\cap(_{s}\Gamma_{0}\cup {_{s}Q_{0}}\cup {_{s}Q_{1}}))^{\bot}}\cap\square_{(1,a_{1}+1,b_{1}),(1,a_{2},b_{2}+1)}.
\end{equation}

By equations (\ref{gamma1-1'}) and (\ref{gamma1-2'}), 
\begin{equation}\label{gamma1-3}
(1,a_{1}+t_{1},b_{1}+s_{1}-q)={^{\bot}(\mathcal{S}\cap(_{s}\Gamma_{0}\cup {_{s}Q_{0}}\cup {_{s}Q_{1}}\cup {_{s}Q'_{0}}\cup {_{s}Q'_{1}}))^{\bot}}\cap\square_{(1,a_{1}+1,b_{1}),(1,a_{2},b_{2}+1)}.
\end{equation}

We claim that $(1,a_{1}+t_{1},b_{1}+s_{1}-q)\in\mathcal{S}$.
Otherwise, we assume that $(1,a_{1}+t_{1},b_{1}+s_{1}-q)$  is not in $\mathcal{S}$. 
Therefore, by equation  (\ref{equ-7}),  $\mathcal{S}\cap{_{s}\Gamma_{1}}\subseteq\square_{(1,a_{2}+1,b_{2}),(1,a_{1}-p,b_{1}+1-q)}$. Thus $(1,a_{1}+t_{1},b_{1}+s_{1}-q)\in {^{\bot}(\mathcal{S}\cap{_{s}\Gamma_{1}})^{\bot}}$. 
By equation (\ref{bi-perp-8}), $(1,a_{1}+t_{1},b_{1}+s_{1}-q)\in {^{\bot}\mathcal{M}^{\bot}}$. 
By equations (\ref{gamma1-1''}) and (\ref{gamma1-2''}), \[(1,a_{1}+t_{1},b_{1}+s_{1}-q)\in {^{\bot}(\mathcal{S}\cap(_{s}Q'_{0}\cup {_{s}Q'_{1}}))^{\bot}}\cup{^{\bot}(\mathcal{S}\cap(_{s}Q_{0}\cup {_{s}Q_{1}}))^{\bot}}.\]
Since $\mathcal{S}=\mathcal{M}\cup(\mathcal{S}\cap{_{s}\Gamma_{1}})\cup(\mathcal{S}\cap(_{s}Q'_{0}\cup {_{s}Q'_{1}}))\cup(\mathcal{S}\cap(_{s}Q_{0}\cup {_{s}Q_{1}})$, $(1,a_{1}+t_{1},b_{1}+s_{1}-q)\in {^{\bot}\mathcal{S}^{\bot}}$. It is a contradiction.  Hence $(1,a_{1}+t_{1},b_{1}+s_{1}-q)\in\mathcal{S}$. 
\end{proof}

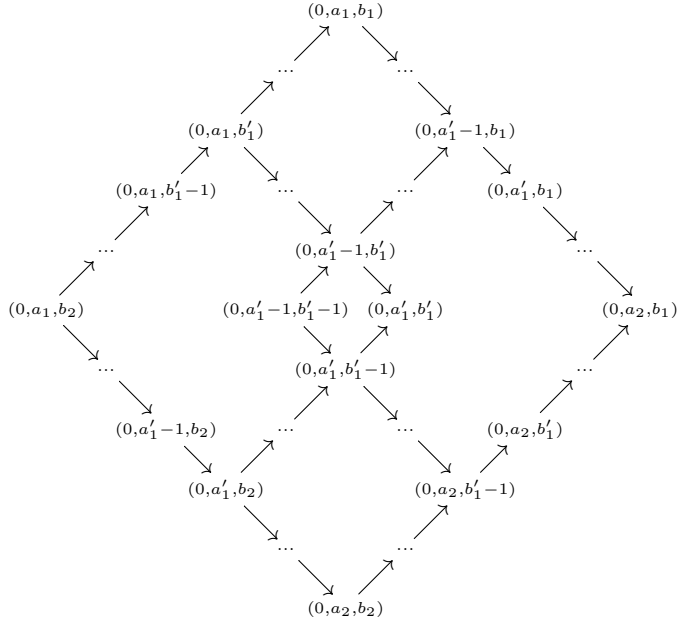
\vspace{-1.3cm}
\begin{small}\[\xymatrix@dr@R=16pt@C=16pt@!0{
&&\scriptstyle(0,a_{1},b_{1})\ar[rr] &&\scriptstyle\scriptstyle\cdots \ar[rr]& &\scriptstyle (0,a'_{1}-1,b_{1})\scriptstyle \ar[rr] && \scriptstyle \scriptstyle(0,a'_{1},b_{1})\ar[rr]&& \scriptstyle\scriptstyle\cdots\ar[rr]&& \scriptstyle(0,a_{2},b_{1})\,.&&\\	
&& \\	
&&\scriptstyle\cdots\ar[uu] &&\scriptstyle&& 	\scriptstyle \cdots\ar[uu] &&\scriptstyle&& &&\scriptstyle\scriptstyle\cdots\ar[uu]&&\\	
&&& \\	
&&	\scriptstyle (0,a_{1},b'_{1})\ar[rr] \ar[uu] &&\scriptstyle\cdots\ar[rr]& & \scriptstyle (0,a'_{1}-1,b'_{1})\ar[rr]\ar[uu]\scriptstyle &&\scriptstyle (0,a'_{1},b'_{1}) && \scriptstyle &&\scriptstyle(0,a_{2},b'_{1})\ar[uu]\\	
&&	& \\
&&\scriptstyle \scriptstyle(0,a_{1},b'_{1}-1)\ar[uu]&& \scriptstyle && \scriptstyle(0,a'_{1}-1,b'_{1}-1)\ar[uu]\ar[rr]\scriptstyle &&\scriptstyle(0,a'_{1},b'_{1}-1)\ar[uu]\ar[rr] &&\scriptstyle\cdots\ar[rr]&&\scriptstyle(0,a_{2},b'_{1}-1)\ar[uu]\\
&&& \\
&&\scriptstyle\cdots \ar[uu] && \scriptstyle && \scriptstyle  &&\scriptstyle\cdots\ar[uu]&&\scriptstyle &&\scriptstyle\cdots\ar[uu]\\ 
&&\\
&&\scriptstyle(0,a_{1},b_{2}) \ar[uu] \ar[rr]&&\scriptstyle\cdots\ar[rr]&&\scriptstyle(0,a'_{1}-1,b_{2})\ar[rr]&&\scriptstyle(0,a'_{1},b_{2})\ar[rr]\ar[uu] &&\scriptstyle\cdots\ar[rr]&&\scriptstyle(0,a_{2},b_{2}) \ar[uu]
}\] 
\end{small}
where $(0,a'_{1},b'_{1})=(0,a_{1}-t_{1},b_{1}+s_{1}-q)=\Omega^{-1}((1,a_{1}-t_{1},b_{1}+s_{1}-q))$. The above  diagram shows us the position of some objects on Theorem \ref{ortho-Euclidean}.

\begin{Prop}\label{square-2'}
Under the above notations on Theorem \ref{ortho-Euclidean}, we have the following results.
\begin{enumerate}[$(1)$]
\item Rectangle areas $\square_{M_{k},(0,a'_{k}-1,b'_{k})}\cup\square_{(0,a'_{k},b'_{k}-1),M_{k+1}}\subseteq\mathcal{F}(\mathcal{S})$ for each $1\leq k\leq\ell$. 
		
\item Rectangle areas $\square_{M'_{k},(1,a_{k+1},b_{k+1}+1)}\cup\square_{(1,a_{k+1}+1,b_{k+1}),M'_{k+1}}\subseteq\mathcal{F}(\mathcal{S})$ for each $1\leq k\leq\ell$. 
		
\item  $\Omega(M_{k})=(1,a_{k},b_{k})$ and  $\Omega^{-1}(M_{k})=(1,a_{k}+1,b_{k}+1)$ are in $\mathcal{F}(\mathcal{S})$ for each  $1\leq k\leq\ell.$
		
\item $\Omega(M'_{k})=(1,a'_{k},b'_{k})$ and  $\Omega^{-1}(M'_{k})=(1,a'_{k}+1,b'_{k}+1)$ are in $\mathcal{F}(\mathcal{S})$ for each $1\leq k\leq\ell.$
\end{enumerate}
\end{Prop}
\begin{proof}
For (1) and (2), we only prove that $k=1$ for (1), Since other cases are similar to prove.
For (3) and (4), we only prove $k=1$ for (3), that is, $\Omega(M_{1})=(1,a_{1},b_{1})$ and  $\Omega^{-1}(M_{1})=(1,a_{1}+1,b_{1}+1)$ are in $\mathcal{F}(\mathcal{S})$. 

\noindent(1) Consider the following triangles:
\begin{align}\label{seqs-1'}
(0,a_{1},b_{1}-j)\xrightarrow{}(0,a_{1},b_{1})\oplus(0,a_{1}+i,b_{1}-j)\xrightarrow{} (0,a_{1}+i,b_{1})\xrightarrow{} (1,a_{1}+1,b_{1}+1-j),
\end{align}
where $j=1,\cdots,q-s_{1}$ and $i=1,\cdots,t_{1}-1$.
\begin{align}\label{seqs-2'}
(0,a_{2}-i,b_{2})\xrightarrow{}(0,a_{2},b_{2})\oplus(0,a_{1}-i,b_{2}+j)\xrightarrow{} (0,a_{2},b_{2}+j)\xrightarrow{} (1,a_{2}-i+1,b_{2}+1),
\end{align}
where $j=1,\cdots,b_{1}-b_{2}+s_{1}-q-1$ and $i=1,\cdots,a_{2}-a_{1}-t_{1}$.

By (1) of (b) and (1) of (a) in Theorem \ref{ortho-Euclidean}, $(0,a_{1},b_{1}-j), (0,a_{1},b_{1})$ and $(0,a_{1}+i,b_{1})$ are contained in $\mathcal{F}(\mathcal{S})$ in triangle (\ref{seqs-1'}). Since $\mathcal{F}(\mathcal{S})$ is closed under direct summands, $(0,a_{1}-i,b_{2}+j)\in\mathcal{F}(\mathcal{S})$ for $j=1,\cdots,q-s_{1}$ and $i=1,\cdots,t_{1}-1$. Thus $\square_{M_{1},Y_{1}}\subseteq\mathcal{F}(\mathcal{S})$, where $Y_{1}=(0,a_{1}+t_{1}-1,b_{1}+s_{1}-q)$.
	
By (1) of (b) and (1) of (a) in Theorem \ref{ortho-Euclidean}, $(0,a_{2}-i,b_{2}), (0,a_{2},b_{2})$ and $(0,a_{2},b_{2}+j)$ are contained in $\mathcal{F}(\mathcal{S})$ in  triangle (\ref{seqs-2'}). Since $\mathcal{F}(\mathcal{S})$ is closed under direct summands,  $(0,a_{1}-i,b_{2}+j)\in\mathcal{F}(\mathcal{S})$ for
 $j=1,\cdots,b_{1}-b_{2}+s_{1}-q-1$ and $i=1,\cdots,a_{2}-a_{1}-t_{1}$. Thus   $\square_{X_{1},M_{2}}\subseteq\mathcal{F}(\mathcal{S})$, where $X_{1}=(0,a_{1}+t_{1},b_{1}-s_{1}-q-1)$.	
By (c) of  Theorem \ref{ortho-Euclidean}, \[\mathcal{S}\cap\square_{\Omega(M_{1}),\Omega(M_{2})}=\mathcal{S}\cap\square_{(1,a_{1},b_{1}),(1,a_{2},b_{2})}=M'_{1}=(1,a_{1}+t_{1},b_{1}+s_{1}-q)=(1,a'_{1},b'_{1}).\]
Therefore $Y_{1}=(0,a_{1}+t_{1}-1,b_{1}+s_{1}-q)=(0,a'_{1}-1,b'_{1})$. $X_{1}=(0,a_{1}+t_{1},b_{1}-s_{1}-q-1)=(0,a'_{1},b'_{1}-1)$.	
Thus rectangle areas $\square_{M_{1},(0,a'_{1}-1,b'_{1})}\cup\square_{(0,a'_{1},b'_{1}-1),M_{2}}\subseteq\mathcal{F}(\mathcal{S}).$

\noindent(3) Consider the following triangle:
\begin{align}\label{seqs-3'}
(0,a_{1},b_{1})\xrightarrow{}(1,a_{1},b_{1})\xrightarrow{} (1,a_{1},b_{1}+1)\oplus(1,a_{1}+1,b_{1}) \xrightarrow{} (1,a_{1}+1,b_{1}+1),
\end{align}
By conclusion  (2), $(1,a_{1},b_{1}+1)\oplus(1,a_{1}+1,b_{1})\in\mathcal{F}(\mathcal{S})$. Thus 
$\Omega(M_{1})=(1,a_{1},b_{1})\in\mathcal{F}(\mathcal{S})$.

Dually, consider the triangle as follows.
\begin{align}\label{seqs-4'}
(1,a_{1},b_{1}+1)\oplus(1,a_{1}+1,b_{1})\xrightarrow{} (1,a_{1}+1,b_{1}+1)\xrightarrow{}(0,a_{1},b_{1}) \xrightarrow{} (0,a_{1},b_{1}+1)\oplus(0,a_{1}+1,b_{1}). 
\end{align}
It is easy to know that $\Omega^{-1}(M_{1})=(1,a_{1}+1,b_{1}+1)\in\mathcal{F}(\mathcal{S})$.
\end{proof}
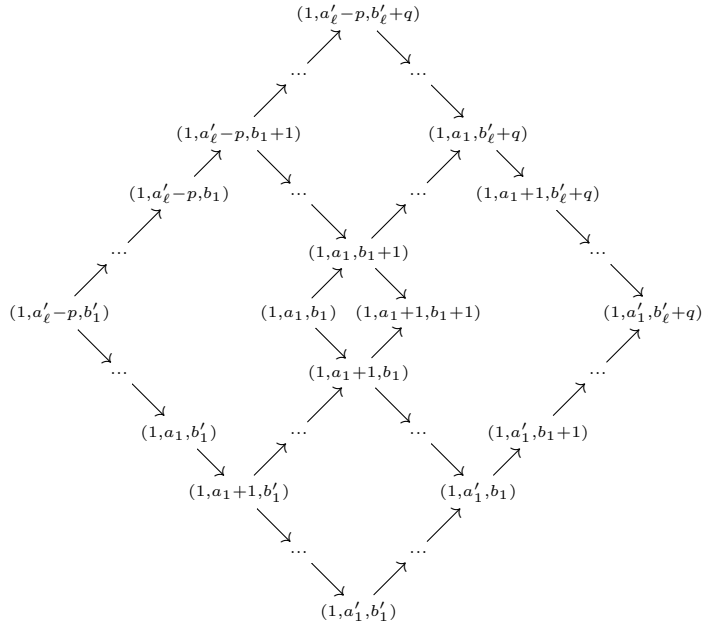
\vspace{-1.3cm}
\begin{small}\[\xymatrix@dr@R=16pt@C=16pt@!0{
&&\scriptstyle(1,a'_{\ell}-p,b'_{\ell}+q)\ar[rr] &&\scriptstyle\scriptstyle\cdots \ar[rr]& &\scriptstyle (1,a_{1},b'_{\ell}+q)\scriptstyle \ar[rr] && \scriptstyle \scriptstyle(1,a_{1}+1,b'_{\ell}+q)\ar[rr]&& \scriptstyle\scriptstyle\cdots\ar[rr]&&\scriptstyle(1,a'_{1},b'_{\ell}+q) \,.&&\\	
&& \\	
&&\scriptstyle\cdots\ar[uu] &&\scriptstyle&& 	\scriptstyle \cdots\ar[uu] &&\scriptstyle&& &&\scriptstyle\scriptstyle\cdots\ar[uu]&&\\	
&&& \\	
&&	\scriptstyle(1,a'_{\ell}-p,b_{1}+1)\ar[rr] \ar[uu] &&\scriptstyle\cdots\ar[rr]& & \scriptstyle (1,a_{1},b_{1}+1)\ar[rr]\ar[uu]\scriptstyle &&\scriptstyle (1,a_{1}+1,b_{1}+1) && \scriptstyle &&\scriptstyle(1,a'_{1},b_{1}+1)\ar[uu]\\	
&&	& \\
&&\scriptstyle(1,a'_{\ell}-p,b_{1})\ar[uu]&& \scriptstyle && \scriptstyle(1,a_{1},b_{1})\ar[uu]\ar[rr]\scriptstyle &&\scriptstyle(1,a_{1}+1,b_{1})\ar[uu]\ar[rr] &&\scriptstyle\cdots\ar[rr]&&\scriptstyle(1,a'_{1},b_{1})\ar[uu]\\
&&& \\
&&\scriptstyle\cdots \ar[uu] && \scriptstyle && \scriptstyle  &&\scriptstyle\cdots\ar[uu]&&\scriptstyle &&\scriptstyle\cdots\ar[uu]\\ 
&&\\
&&\scriptstyle(1,a'_{\ell}-p,b'_{1})  \ar[uu] \ar[rr]&&\scriptstyle\cdots\ar[rr]&&\scriptstyle(1,a_{1},b'_{1})\ar[rr]&&\scriptstyle(1,a_{1}+1,b'_{1})\ar[rr]\ar[uu] &&\scriptstyle\cdots\ar[rr]&&\scriptstyle(1,a'_{1},b'_{1}) \ar[uu]
}\] 
\end{small}
The above  diagram show us the position of some rectangle areas on (2) of Proposition  \ref{square-2'}. 
\begin{Prop}\label{square-1}
Let $A$ be a 2-domestic Brauer graph algebra and $\mathcal{S}$ a maximal orthogonal system  which contains at least one object on an Euclidean component. We assume $\mathcal{M}=\{M_{1},\cdots,M_{\ell}\}$ {(\rm}resp, $\mathcal{M'}=\{M'_{1},\cdots,M'_{\ell}\}${\rm)} is the subset of $\mathcal{S}$ which is contained in the stable Euclidean component $_{s}\Gamma_{0}$ {\rm(}resp. $_{s}\Gamma_{1}${\rm)}, where $M_{k}=(0,a_{k},b_{k})$ {\rm(}resp. $M'_{k}=(1,a'_{k},b'_{k})${\rm)} for $1\leq k\leq\ell$. Then 
\begin{enumerate}[$(a)$]
\item The following subsets on $_{s}\Gamma_{1}$ are contained $\mathcal{F}(\mathcal{S})$.
\begin{enumerate}[$(1)$]
\item $\{(1,a_{1}+1+j,b_{1}+1)\mid j=0,1,\cdots, a'_{1}-a_{1}-1\}\cup\{(1,a'_{\ell}-p+j,b_{1})\mid j=0,1,\cdots, a_{1}+p-a'_{\ell}\}\cup\{(1,a_{1},b'_{1}+j)\mid j=0,1,\cdots, b_{1}-b'_{1}\}\cup\{(1,a_{1}+1,b_{1}+1+j)\mid j=0,1,\cdots, b'_{\ell}+q-b_{1}-1\}$.
\item $\{(1,a_{k}+1+j,b_{k}+1)\mid j=0,1,\cdots, a'_{k}-a_{k}-1\}\cup\{(1,a'_{k-1}+j,b_{k})\mid j=0,1,\cdots, a_{k}-a'_{k-1}\}\cup\{(1,a_{k},b_{k}-j)\mid j=0,1,\cdots, b_{k}-b'_{k}\}\cup\{(1,a_{k}+1,b_{k}+1+j)\mid j=0,1,\cdots, b'_{k-1}-b_{k}-1\}$ for $2\leq k\leq\ell$. 
\item Rectangle areas $\square_{M'_{k},(1,a_{k+1}+1,b_{k+1})}\cup\square_{(1,a_{k+1},b_{k+1}+1),M'_{k+1}}$ for any integer $1\leq k\leq\ell$.
\end{enumerate} 
\item Dually, the following subsets on $_{s}\Gamma_{0}$ are contained in $\mathcal{F}(\mathcal{S})$.
\begin{enumerate}[$(1)$]
\item $\{(0,a'_{1}+j,b'_{1})\mid j=0,1,\cdots, a_{2}-a'_{1}\}\cup\{(0,a_{1}+j,b'_{1}-1)\mid j=0,1,\cdots, a'_{1}-1-a_{1}\}\cup\{(1,a'_{1},b'_{1}+j)\mid j=0,1,\cdots, b_{1}-b'_{1}\}\cup\{(1,a'_{1}-1,b_{2}+j)\mid j=0,1,\cdots, b'_{1}-1-b_{2}\}$.
\item $\{(0,a'_{k}+j,b'_{k})\mid j=0,1,\cdots, a_{k+1}-a'_{k}\}\cup\{(0,a_{k}+j,b'_{k}-1)\mid j=0,1,\cdots, a'_{k}-1-a_{k}\}\cup\{(1,a'_{k},b'_{k}+j)\mid j=0,1,\cdots, b_{k}-b'_{k}\}\cup\{(1,a'_{k}-1,b_{k+1}+j)\mid j=0,1,\cdots, b'_{k}-1-b_{k+1}\}$  for $2\leq k\leq\ell$. Note that $a_{\ell+1}=a_{1}-p$ and $b_{\ell+1}=b_{1}+q$ {\rm(}resp. $a'_{\ell+1}=a'_{1}-p$ and $b'_{\ell+1}=b'_{1}+q${\rm)}.
\item Rectangle areas $\square_{M_{k},(0,a'_{k},b'_{k}-1)}\cup\square_{(0,a'_{k}-1,b'_{k}),M_{k+1}}$ for any integer $1\leq k\leq\ell$. 
\end{enumerate}
\end{enumerate}
\end{Prop}
\begin{proof}
We only prove (1) of (a), since (2) of (a) (resp. (b)) are similar to prove, (3) of (a) (resp. (b)) is a direct consequence of (1) and (2) of (a) (resp. (b)). We divide (1) of (a) into four claims as follows.

\noindent{\bf Claim one:} $\{(1,a_{1}+1+j,b_{1}+1)\mid j=0,1,\cdots, a'_{1}-a_{1}-1\}\subseteq\mathcal{F}(\mathcal{S})$.

Consider the triangle as follows.
\begin{align}\label{seqs-2''}
(1,A_{j},B)\oplus(1,A_{j}+1,B-1)\xrightarrow{}(1,A_{j}+1,B)\xrightarrow{}(0,A_{j},B-1)\xrightarrow{}(0,A_{j},B)\oplus(0,A_{j}+1,B-1).
\end{align}
where $a_{1}+j=A_{j}$ and $b_{1}+1=B$ for $0\leq j\leq a'_{1}-a_{1}-1.$

By (1) of (a) in Theorem \ref{ortho-Euclidean},  
\[(0,A_{j},B-1)=(0,a_{1}+j,b_{1})\in\mathcal{F}(\mathcal{S})\] for $0\leq j\leq t_{1}-1$. 
Note that $t_{1}=a'_{1}-a_{1}$. 
By (2) of Proposition \ref{square-2'},  \[(1,A_{j}+1,B-1)=(1,a_{1}+j+1,b_{1})\in\square_{(1,a_{1}+1,b_{1}),M'_{1}}\subseteq\mathcal{F}(\mathcal{S})\] for $0\leq j\leq a'_{1}-a_{1}-1$. 
By  (3) of Proposition \ref{square-2'}, $\Omega^{-1}(M_{1})=(1,a_{1}+1,b_{1}+1)\in\mathcal{F}(\mathcal{S})$. By the triangle (\ref{seqs-2''}) and by induction, \[(1,A_{j}+1,B)=(1,a_{1}+j+1,b_{1}+1)\in\mathcal{F}(\mathcal{S})\] for  $j=0,1,\cdots, a'_{1}-a_{1}-1$.

\noindent{\bf Claim two:} $\{(1,a'_{\ell}-p+j,b_{1})\mid j=0,1,\cdots, a_{1}+p-a'_{\ell}\}\subseteq\mathcal{F}(\mathcal{S})$.

Note that $\{(1,a'_{\ell}-p+j,b_{1})\mid j=0,1,\cdots, a_{1}+p-a'_{\ell}\}
=\{(1,a_{1}-j,b_{1})\mid j=0,1,\cdots, a_{1}-a'_{\ell}+p\}$.
Consider the triangle as follows.
\begin{align}\label{seqs-1}
(0,a_{1}-j,b_{1})\xrightarrow{}(1,a_{1}-j,b_{1})\xrightarrow{}(1,a_{1}-j,b_{1}+1)\oplus(1,a_{1}-j+1,b_{1})\xrightarrow{}(1,a_{1}-j+1,b_{1}+1)
\end{align}
 for $0\leq j\leq a_{1}-a'_{\ell}+p$.
 
By (3) of (a) in Theorem \ref{ortho-Euclidean}, $(0,a_{1}-j,b_{1})$ is in 
$\mathcal{F}(\mathcal{S})$ for $0\leq j\leq a_{1}-a'_{\ell}+p$.
By  (2) of Proposition \ref{square-2'}, \[(1,a_{1}-j,b_{1}+1)\in\square_{M'_{\ell},(1,a_{1},b_{1}+1)}\subseteq\mathcal{F}(\mathcal{S})\]
for $0\leq j\leq a_{1}-a'_{\ell}+p$. By  (3) of Proposition \ref{square-2'}, $\Omega(M_{1})=(1,a_{1},b_{1})\in\mathcal{F}(\mathcal{S}).$
Thus, by induction and triangle (\ref{seqs-1}), $(1,a_{1}-j,b_{1})\in\mathcal{F}(\mathcal{S})$ for $j=0,1,\cdots, a_{1}-a'_{\ell}+p.$

\noindent{\bf Claim three:} $\{(1,a_{1},b'_{1}+j)\mid j=0,1,\cdots, b_{1}-b'_{1}\}\subseteq\mathcal{F}(\mathcal{S})$.

Note that $\{(1,a_{1},b'_{1}+j)\mid j=0,1,\cdots, b_{1}-b'_{1}\}=\{(1,a_{1},b_{1}-j)\mid j=0,1,\cdots, b_{1}-b'_{1}\}$.
Consider the triangle as follows.
\begin{align}\label{seqs-2}
(0,a_{1},b_{1}-j)\xrightarrow{}(1,a_{1},b_{1}-j)\xrightarrow{}(1,a_{1},b_{1}-j-1)\oplus(1,a_{1}+1,b_{1}-j)\xrightarrow{}(1,a_{1}+1,b_{1}-j+1).
\end{align}
for $0\leq j\leq b_{1}-b'_{1}$.

By (1) of (b) in Theorem \ref{ortho-Euclidean}, $(0,a_{1},b_{1}-j)$ is in 
$\mathcal{F}(\mathcal{S})$ for $0\leq j\leq b_{1}-b'_{1}$.
By  (2) of Proposition \ref{square-2'}, \[(1,a_{1}+1,b_{1}-j)\in\square_{(1,a_{1}+1,b_{1}),M'_{1}}\subseteq\mathcal{F}(\mathcal{S})\]
 for $0\leq j\leq b_{1}-b'_{1}$. By  (3) of Proposition \ref{square-2'}, $\Omega^(M_{1})=(1,a_{1},b_{1})\in\mathcal{F}(\mathcal{S})$. Thus, by induction and triangle (\ref{seqs-2}), $(1,a_{1},b_{1}-j)\in\mathcal{F}(\mathcal{S})$ for $j=0,1,\cdots, b_{1}-b'_{1}.$

\noindent{\bf Claim four:} $\{(1,a_{1}+1,b_{1}+1+j)\mid j=0,1,\cdots, b'_{\ell}+q-b_{1}-1\}\subseteq\mathcal{F}(\mathcal{S})$.

Consider the triangle as follows.
\begin{align}\label{seqs-3}
(1,A,B_{j})\oplus(1,A-1,B_{j}+1)\xrightarrow{}(1,A,B_{j}+1)\xrightarrow{}(0,A-1,B_{j})\xrightarrow{}(0,A,B_{j})\oplus(0,A-1,B_{j}+1).
\end{align}
for $0\leq j\leq b'_{\ell}+q-b_{1}-1$, where $A=a_{1}+1,B_{j}=b_{1}+j$.

By   (2) of Proposition \ref{square-2'}, \[(1,A-1,B_{j}+1)=(1,a_{1},b_{1}+j+1)\in\square_{M'_{\ell},(1,a_{1},b_{1}+1)}\subseteq\mathcal{F}(\mathcal{S})\] for $0\leq j\leq b'_{\ell}+q-b_{1}-1$.
By Proposition \ref{square-2'} (1),  \[(0,A-1,B_{j})=(0,a_{1},b_{1}+j)\in\square_{(0,a'_{\ell}-p,b'_{\ell}+q-1),M_{1}}\subseteq\mathcal{F}(\mathcal{S})\] for $0\leq j\leq b'_{\ell}+q-b_{1}-1$. 

By  (3) of Proposition \ref{square-2'}, $\Omega^{-1}(M_{1})=(1,a_{1}+1,b_{1}+1)\in\mathcal{F}(\mathcal{S})$. By the triangle (\ref{seqs-3}) and by induction,  \[(1,A,B_{j}+1)=(1,a_{1}+1,b_{1}+j+1)\in\mathcal{F}(\mathcal{S})\] for  $j=0,1,\cdots, b'_{\ell}+q-b_{1}-1$. 
\end{proof}

\begin{Cor}\label{square-2}
Under the above notations on Proposition \ref{square-1}, we have the following 	results.
\begin{enumerate}[$(1)$]
\item Rectangle areas $\square_{(0,a'_{k}-1,b_{k+1}),(0,a_{k+1},b'_{k+1}-1)}\cup\square_{(0,a_{k+1},b'_{k}),(0,a'_{k+1},b_{k+1})}\subseteq\mathcal{F}(\mathcal{S})$ for any integer $1\leq k\leq\ell$. 

\item Rectangle areas $\square_{(1,a_{k},b'_{k}),(1,a'_{k},b_{k+1})}\cup\square_{(1,a'_{k},b_{k}+1),(1,a_{k+1}+1,b'_{k})}\subseteq\mathcal{F}(\mathcal{S})$ for any integer $1\leq k\leq\ell$. 

\item Rectangle areas $\square_{(0,a'_{k}-1,b'_{k}),(0,a'_{k+1},b'_{k+1}-1)}\cup\square_{(1,a_{k},b_{k}+1),(1,a_{k+1}+1,b_{k+1})}\subseteq\mathcal{F}(\mathcal{S})$ for any integer $1\leq k\leq\ell$.
\end{enumerate}
\end{Cor}
\begin{proof}
We only prove conclusion (1) for $k=\ell$, since conclusion (2) is similar to prove and conclusion (3) is a direct consequence of (1) and (2).

\noindent{\bf Claim one:} Rectangle area $\square_{(0,a'_{\ell}-1,b_{1}),(0,a_{1},b'_{1}-1)}\subseteq\mathcal{F}(\mathcal{S})$.

Consider the triangle as follows.
\begin{align}\label{seqs-3''}
(1,a_{1},b_{1})\xrightarrow{}(0,A_{m},B_{n})\xrightarrow{}(0,A_{m},b_{1})\oplus(0,a_{1},B_{n})\xrightarrow{} (0,a_{1},b_{1}).
\end{align}
where $A_{m}=a'_{\ell}-p-1+m$  for $0\leq m\leq a_{1}+p-a'_{\ell}+1$ and $B_{n}=b'_{1}-1+n$ for $0\leq n\leq b_{1}-b'_{1}+1$.
 
By  (3) of Proposition \ref{square-2'}, $\Omega(M_{1})=(1,a_{1},b_{1})\in\mathcal{F}(\mathcal{S})$. By (b) of Proposition \ref{square-1}, 
\[(0,A_{m},b_{1})=(0,a'_{\ell}-p-1+m,b_{1})\in\square_{(0,a'_{\ell}-p,b'_{\ell}+q-1),M_{1}}\subseteq\mathcal{F}(\mathcal{S})\] for $0\leq m\leq a_{1}+p-a'_{\ell}+1$. 
By  (b) of Proposition \ref{square-1},  
\[(0,a_{1},B_{n})=(0,a_{1},b'_{1}-1+n)\in\square_{M_{1},(0,a'_{1},b'_{1}-1)}\subseteq\mathcal{F}(\mathcal{S})\] for $0\leq n\leq b_{1}-b'_{1}+1$. Therefore, by triangle (\ref{seqs-3''}),  \[(0,A_{m},B_{n})=(0,a'_{\ell}-p-1+m,b'_{1}-1+n)\in\mathcal{F}(\mathcal{S})\] for $0\leq m\leq a_{1}+p-a'_{\ell}+1$ and  $0\leq n\leq b_{1}-b'_{1}+1$. Thus $\square_{(0,a'_{\ell}-1,b_{1}),(0,a_{1},b'_{1}-1)}\subseteq\mathcal{F}(\mathcal{S}).$
 
\noindent{\bf Claim two:} Rectangle area $\square_{(0,a_{1},b'_{\ell}),(0,a'_{1},b_{1})}\subseteq\mathcal{F}(\mathcal{S})$.
 
 Consider the following triangle.
\begin{align}\label{seqs-4''}
(0,a_{1},b_{1})\xrightarrow{} (0,a_{1}+m,b_{1})\oplus(0,a_{1},b_{1}+n)\xrightarrow{}(0,a_{1}+m,b_{1}+n)\xrightarrow{} (1,a_{1}+1,b_{1}+1).
\end{align}
where  $0\leq m\leq a'_{1}-a_{1}$ and $0\leq n\leq b'_{\ell}+p-b_{1}$. 

By  (3) of Proposition \ref{square-2'}, \[\Omega^{-1}(M_{1})=(1,a_{1}+1,b_{1}+1)\in\mathcal{F}(\mathcal{S}).\]
By   (b) of Proposition \ref{square-1},  \[(0,a_{1}+m,b_{1})\in\square_{M_{1},(0,a'_{1},b'_{1}-1)}\subseteq\mathcal{F}(\mathcal{S})\] $0\leq m\leq a'_{1}-a_{1}$.
By  (b) of Proposition \ref{square-1}, \[(0,a_{1},b_{1}+n)\in\square_{(0,a'_{\ell}-1,b'_{\ell}),M_{1}}\subseteq\mathcal{F}(\mathcal{S})\] for $0\leq n\leq b'_{\ell}+p-b_{1}$. 
Thus $\square_{(0,a_{1},b'_{\ell}),(0,a'_{1},b_{1})}\subseteq\mathcal{F}(\mathcal{S})$.
\end{proof}

The following diagrams indicate the location of some rectangle areas in Corollary \ref{square-2}. 
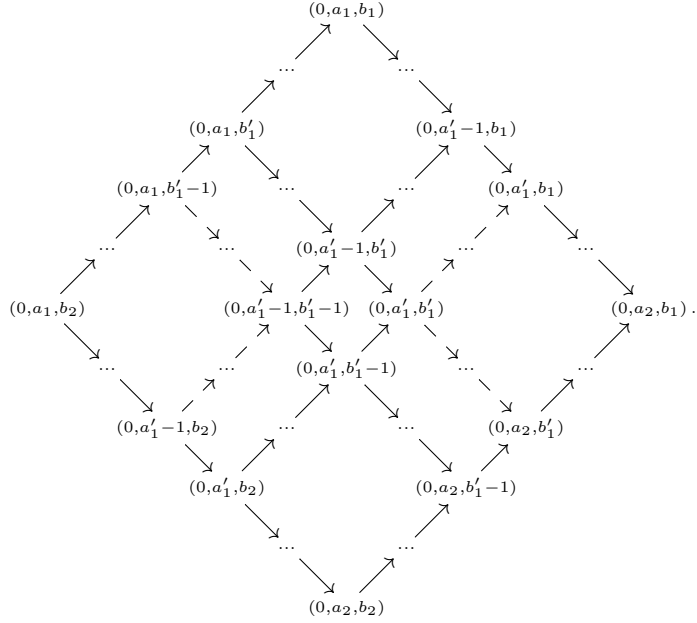
\vspace{-0.5cm}
\begin{small}\[\xymatrix@dr@R=16pt@C=16pt@!0{
&&\scriptstyle(0,a_{1},b_{1})\ar[rr] &&\scriptstyle\scriptstyle\cdots \ar[rr]& &\scriptstyle (0,a'_{1}-1,b_{1})\scriptstyle \ar[rr] && \scriptstyle (0,a'_{1},b_{1})\ar[rr]&& \scriptstyle\cdots\ar[rr]&&\ \  \scriptstyle (0,a_{2},b_{1})\,.&&\\	
&& \\	
&&\scriptstyle\cdots\ar[uu] &&\scriptstyle&& 	\scriptstyle \cdots\ar[uu] &&\scriptstyle\cdots\ar@{-->}[uu]&& &&\scriptstyle\scriptstyle\cdots\ar[uu]&&\\	
&&& \\	
&&\scriptstyle(0,a_{1},b'_{1})\ar[rr] \ar[uu] &&\scriptstyle\cdots\ar[rr]& & \scriptstyle (0,a'_{1}-1,b'_{1})\ar[rr]\ar[uu]\scriptstyle &&\scriptstyle (0,a'_{1},b'_{1})\ar@{-->}[uu] \ar@{-->}[rr]&& \scriptstyle \cdots\ar@{-->}[rr] &&\scriptstyle (0,a_{2},b'_{1})\ar[uu]\\	
&&	& \\
&& \scriptstyle(0,a_{1},b'_{1}-1)\ar@{-->}[rr]\ar[uu]&& \scriptstyle \cdots\ar@{-->}[rr] && \scriptstyle(0,a'_{1}-1,b'_{1}-1)\ar[uu]\ar[rr]\scriptstyle &&\scriptstyle(0,a'_{1},b'_{1}-1)\ar[uu]\ar[rr] &&\scriptstyle\cdots\ar[rr]&&\scriptstyle (0,a_{2},b'_{1}-1)\ar[uu]\\
&&& \\
&&\scriptstyle\cdots \ar[uu] && \scriptstyle  && \scriptstyle \cdots\ar@{-->}[uu] &&\scriptstyle\cdots \ar[uu]&&\scriptstyle &&\scriptstyle \cdots\ar[uu]\\ 
&&\\
&&\scriptstyle(0,a_{1},b_{2}) \ar[uu] \ar[rr]&&\scriptstyle\cdots\ar[rr]&&\scriptstyle(0,a'_{1}-1,b_{2})\ar@{-->}[uu]\ar[rr]&&\scriptstyle(0,a'_{1},b_{2})\ar[rr]\ar[uu] &&\scriptstyle\cdots\ar[rr]&&\scriptstyle(0,a_{2},b_{2}) \ar[uu]}\] 
\end{small}
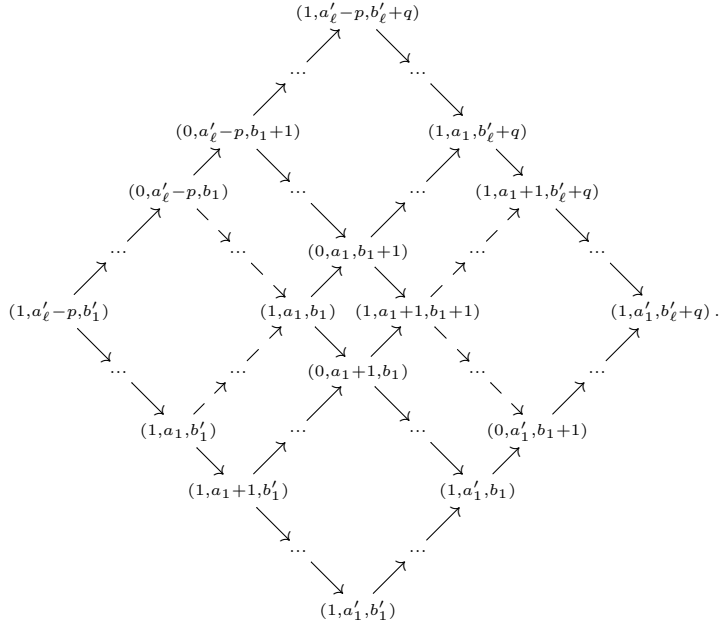
\vspace{-0.5cm}
\begin{small}\[\xymatrix@dr@R=16pt@C=16pt@!0{
&&\scriptstyle(1,a'_{\ell}-p,b'_{\ell}+q)\ar[rr] &&\scriptstyle\scriptstyle\cdots \ar[rr]& &\scriptstyle (1,a_{1},b'_{\ell}+q)\scriptstyle \ar[rr] && \scriptstyle (1,a_{1}+1,b'_{\ell}+q)\ar[rr]&& \scriptstyle\cdots\ar[rr]&&\ \  \scriptstyle (1,a'_{1},b'_{\ell}+q)\,.&&\\	
&& \\	
&&\scriptstyle\cdots\ar[uu] &&\scriptstyle&& 	\scriptstyle \cdots\ar[uu] &&\scriptstyle\cdots\ar@{-->}[uu]&& &&\scriptstyle\scriptstyle\cdots\ar[uu]&&\\	
&&& \\	
&&\scriptstyle(0,a'_{\ell}-p,b_{1}+1)\ar[rr] \ar[uu] &&\scriptstyle\cdots\ar[rr]& & \scriptstyle (0,a_{1},b_{1}+1)\ar[rr]\ar[uu]\scriptstyle &&\scriptstyle (1,a_{1}+1,b_{1}+1)\ar@{-->}[uu] \ar@{-->}[rr]&& \scriptstyle \cdots\ar@{-->}[rr] &&\scriptstyle (0,a'_{1},b_{1}+1)\ar[uu]\\	
&&	& \\
&& \scriptstyle(0,a'_{\ell}-p,b_{1})\ar@{-->}[rr]\ar[uu]&& \scriptstyle \cdots\ar@{-->}[rr] && \scriptstyle(1,a_{1},b_{1})\ar[uu]\ar[rr]\scriptstyle &&\scriptstyle(0,a_{1}+1,b_{1})\ar[uu]\ar[rr] &&\scriptstyle\cdots\ar[rr]&&\scriptstyle (1,a'_{1},b_{1})\ar[uu]\\
&&& \\
&&\scriptstyle\cdots \ar[uu] && \scriptstyle  && \scriptstyle \cdots\ar@{-->}[uu] &&\scriptstyle\cdots \ar[uu]&&\scriptstyle &&\scriptstyle \cdots\ar[uu]\\ 
&&\\
&&\scriptstyle(1,a'_{\ell}-p,b'_{1}) \ar[uu] \ar[rr]&&\scriptstyle\cdots\ar[rr]&&\scriptstyle(1,a_{1},b'_{1})\ar@{-->}[uu]\ar[rr]&&\scriptstyle(1,a_{1}+1,b'_{1})\ar[rr]\ar[uu] &&\scriptstyle\cdots\ar[rr]&&\scriptstyle(1,a'_{1},b'_{1}) \ar[uu]
	}\] 
\end{small}

\begin{Them}\label{sms-BGA}
Let $A$ be a 2-domestic Brauer graph algebra and $\mathcal{S}$ a maximal orthogonal system  which contains at least one object on an Euclidean component. Then $\Omega^{-1}(\mathcal{S})\subseteq\mathcal{F}(\mathcal{S})$, in particular, $\mathcal{S}$ is a simple-minded system on $A$-$\stmod.$
\end{Them}
\begin{proof}
Take $\mathcal{M}=\mathcal{S}\cap{_{s}\Gamma_{0}}$  and   $\mathcal{M'}=\mathcal{S}\cap{_{s}\Gamma_{1}}$. We assume that $\mathcal{M}=\{M_{1},\cdots,M_{\ell}\}$ and $\mathcal{M'}=\{M'_{1},\cdots,M'_{\ell}\}$, where $M_{i}=(0,a_{i},b_{i})$ and $M'_{i}=(1,a'_{i},b'_{i})$. Without loss of generality, we assume that $(0,a_{1},b_{1})=(0,a,b)$. 
By (3) and (4) of Proposition \ref{square-2'}, $\Omega(M_{k})=(1,a_{k},b_{k})$ and  $\Omega^{-1}(M_{k})=(1,a_{k}+1,b_{k}+1)$ are in $\mathcal{F}(\mathcal{S})$ for any $1\leq k\leq\ell.$ and $\Omega(M'_{k})=(1,a'_{k},b'_{k})$ and  $\Omega^{-1}(M'_{k})=(1,a'_{k}+1,b'_{k}+1)$ are in $\mathcal{F}(\mathcal{S})$ for any $1\leq k\leq\ell.$ 

Now we prove that, for every $\tau$-periodic module $T\in\mathcal{S}$, $\Omega^{-1}(T)\in\mathcal{F}(\mathcal{S}).$ 
By Theorem \ref{ortho-Euclidean}, for $\tau$-periodic module $T\in\mathcal{S}$, there are triangles as follows. 
\begin{align}\label{seqs-5'}
T\xrightarrow{} W_{1}\xrightarrow{}W_{2} \xrightarrow{} \Omega^{-1}(T),
\end{align}
or 
\begin{align}\label{seqs-6'}
W_{1}\xrightarrow{}W_{2}\xrightarrow{} T\xrightarrow{}\Omega^{-1}(W_{1})
\end{align}
satisfying that both $W_{1}$ and $W_{2}$ are in $\mathcal{F}(\mathcal{S})$.

Rotating triangle (\ref{seqs-5'}) to the right once, we have
\begin{align}\label{seqs-9'}
W_{1}\xrightarrow{}W_{2} \xrightarrow{} \Omega^{-1}(T)\xrightarrow{}\Omega^{-1}(W_{1}),
\end{align}
Acting the functor $\Omega^{-1}$ on the triangle (\ref{seqs-6'}), we have 
\begin{align}\label{seqs-10'}
\Omega^{-1}(W_{1})\xrightarrow{}\Omega^{-1}(W_{2})\xrightarrow{} \Omega^{-1}(T)\xrightarrow{}\Omega^{-2}(W_{1}).
\end{align}
We only consider the case that $W_{1}$ and $W_{2}$ are in (1) of (a) over Theorem \ref{ortho-Euclidean}, since other cases are similar to handle.

\noindent{\bf Case one}: $W_{1}$ and $W_{2}$ are in the set  $\{(0,a_{1}+j,b_{1})\mid j=0,1,\cdots, t_{1}-1\}$.

In triangle (\ref{seqs-9'}),
by  (a) of Proposition \ref{square-1}, \[\Omega^{-1}(W_{1})\in\square_{(1,a_{1},b_{1}+1),M'_{1}}\subseteq\square_{(1,a_{1},b_{1}+1),(1,a_{2}+1,b_{2})}\subseteq\mathcal{F}(\mathcal{S}).\] 
Since  $W_{2}\in\mathcal{F}(\mathcal{S})$, $\Omega^{-1}(T)\in\mathcal{F}(\mathcal{S})$.

In triangle (\ref{seqs-10'}), 
by  (a) of Proposition \ref{square-1}, $\Omega^{-1}(W_{2})\in\mathcal{F}(\mathcal{S})$. By Corollary \ref{square-2},  
\[\Omega^{-2}(W_{1})\in\{(0,a_{1}+j+1,b_{1}+1)\mid j=0,1,\cdots, t_{1}-1\}\subseteq\square_{(0,a'_{\ell}-1,b'_{\ell}),(0,a'_{1}-1,b'_{1})}\subseteq\mathcal{F}(\mathcal{S}).\]
Thus $\Omega^{-1}(T)\in\mathcal{F}(\mathcal{S})$.

\noindent{\bf Case two}: $W_{1}$ and $W_{2}$ are in the set  $\{(0,a_{2}-j,b_{2})\mid j=0,1,\cdots,a_{2}-a_{1}-t_{1}\}$.

For every object  $W\in\{(0,a_{2}-j,b_{2})\mid j=0,1,\cdots,a_{2}-a_{1}-t_{1}\}$, \[\Omega^{-1}(W)\in\square_{(1,a_{1},b_{1}+1),(1,a_{2}+1,b_{2})}\subseteq\mathcal{F}(\mathcal{S}).\]
Thus $\Omega^{-1}(W_{1})$ and $\Omega^{-1}(W_{2})$ are in $\mathcal{F}(\mathcal{S}).$

In triangle (\ref{seqs-9'}),
 $\Omega^{-1}(W_{1})\in\mathcal{F}(\mathcal{S})$. Since  $W_{2}\in\mathcal{F}(\mathcal{S})$, $\Omega^{-1}(T)\in\mathcal{F}(\mathcal{S})$.

In triangle (\ref{seqs-10'}), 
 $\Omega^{-1}(W_{2})\in\mathcal{F}(\mathcal{S})$. By Corollary \ref{square-2}, 
\[\Omega^{2}(W_{1})\in\{(0,a_{2}-j+1,b_{2}+1)\mid 0\leq j\leq a_{2}-a_{1}-t_{1}\}\subseteq\square_{(0,a'_{1}-1,b'_{1}),(0,a'_{2},b'_{2}-1)}\subseteq\mathcal{F}(\mathcal{S}).\] 
 Thus $\Omega^{-1}(T)\in\mathcal{F}(\mathcal{S})$. 

Therefore, for every $M\in\mathcal{S}$, we have  $\Omega^{-1}(M)\in\mathcal{F}(\mathcal{S})$.  Thus $\mathcal{S}$ is a simple-mined system on $A$-$\stmod$ by Theorem \ref{BGA-sms}.
\end{proof}

\begin{Cor}
Let $A$ be a 2-domestic Brauer graph algebra and $\mathcal{S}$ an orthogonal system on an Euclidean component. Then $\mathcal{S}$ extends to a simple-minded system on $A$-$\stmod$. In particular, $\mathcal{F}(\mathcal{S})$ is functorially finite on $A$-$\stmod$.
\end{Cor}
\begin{proof}
It follows from Theorem \ref{sms-BGA} and Theorem \ref{subset-of-sms}.
\end{proof}

\begin{Cor}
Let $A$ be a 2-domestic Brauer graph algebra and $\mathcal{S}$  a simple-minded system on $A$-$\stmod$. Then the cardinality of $\mathcal{S}$ is the number of non-projective simple $A$-modules.
\end{Cor}
\begin{proof}
It is a direct consequence of  Theorem \ref{sms-BGA} and Corollary \ref{cardinality-of-ortho-sys-4}.
\end{proof}

\begin{Cor}
Let $A$ be a 2-domestic Brauer graph algebra  and $\mathcal{S}$ a weakly simple-minded system with a finite cardinality. Then $\mathcal{S}$ is a simple-minded system on $A$-$\stmod$.
\end{Cor}
\begin{proof}
By Corollary \ref{2-BGA-quasi-tube-bricks}, any quasi-simple in a homogeneous tube is mutual stable orthogonal with each quasi-simple on ${_{s}Q_{0}}\cup{_{s}Q_{1}}\cup{_{s}Q'_{0}}\cup{_{s}Q'_{1}}$. Thus if $\mathcal{S}$  does not contain any object on  Euclidean components, then there are quasi-simples of infinitely many homogeneous tubes contained in $\mathcal{S}$. Thus the cadrinality of $\mathcal{S}$ is infinite. Hence every weakly simple-minded system with a finite cardinality contains at least one object for an Euclidean component.  By Theorem  \ref{subset-of-sms}, $\mathcal{S}$ is a simple-minded system on $A$-$\stmod$.

\end{proof}

We \cite{Z0} provided a class of weakly simple-minded systems which are not simple-minded systems over a 2-domestic Brauer graph algebra.    Note that those weakly simple-minded systems contains objects in homogeneous tubes and they have a infinite cardinality. We conjecture that  a weakly simple-minded system $\mathcal{S}$ with a finite cardinality is a simple-midned system over a self-injective algerbra.

\begin{Conj}
Let $A$ be a self-injective algebra  and $\mathcal{S}$ a weakly simple-minded system with a finite cardinality. Then $\mathcal{S}$ is a simple-minded system on $A$-$\stmod$.	
\end{Conj}


\begin{thebibliography}{99}	 
\addcontentsline{toc}{chapter}{\numberline{}{\bf 参考文献}}

\bibitem{AAC}{{ T. Adachi, T. Aihara and A. Chan,} Classification of two-term tilting complexes over Brauer graph algebras. Math. Z. 290 (2018), 1-36.}		
		
\bibitem{Ap}{{ M. Antipov,} Invariants of the stable equivalence of symmetric special biserial algebras. J. Math. Sci. 140 (2007), 611–621.}
		
\bibitem{ARS}{{ M. Auslander, I. Reiten and S. Smal\o,}  Representation theory of artin algebras. Cambridge University Press, 1995.}
		
\bibitem{AZ}{{ M. Antipov and  A. Zvonareva,} On stably biserial algebras and the Auslander-Reiten conjecture for special biserial algebras.  J. Math. Sci. 240 (2019), 375–394.}

\bibitem{AZ1}{{ M. Antipov and  A. Zvonareva,} Brauer graph algebras are closed under derived equivalence.  Math. Z. 301 (2022), 1963–1981.}
		
	
\bibitem{BoS}{{ R. Bocian and A. Skowro\'{n}ski,} Symmetric special biserial algebras of Euclidean type. Colloq. Math. 96 (2003), 121-148.} 
		
\bibitem{CLZ}{{ A. Chan, Y. Liu and Z. Zhang,} On Simple-minded systems and $\tau$-periodic modules of self-injective algebras. J. Algebra 560 (2020), 416-441.}
		

\bibitem{CB}{{ W. Crawley-Boevey,} On tame algebras and bocses. Proc. London Math. Soc. 56 (1988), 451-483.}
		
\bibitem{D1}{{ D. Duffield,} Auslander-Reiten components of Brauer graph algebras. J. Algebra  508 (2018), 475-511. }

\bibitem{Dugas}{{ A. Dugas,} Torsion pairs and simple-minded systems in triangulated categories. Appl. Categ. Structures  23  (2015), 507-526.}

\bibitem{D}{{ Y. Drozd, } Tame and wild matrix problems, in Representation Theory II. Lecture Notes in Math.  832, Springer Verlag, 242-258, 1980.}		

\bibitem{DF}{{ P. Donovan and M. Freislich}, The indecomposable modular representations of certain groups with dihedral Sylow subgroup. Math. Ann. 238 (1978), 207-216.}		
		
\bibitem{E}{{ K. Erdmann,}  Blocks of tame representation type and related algebras. Lecture Notes in Math.  1428, Springer, Berlin, 1990.} 
		
\bibitem{EK}{{ K. Erdmann and O. Kerner,} On the stable module category of a self-injective algebra. Trans. Amer. Math. Soc. 352 (2006), 2389--2405.}

\bibitem{GLYZ2}{{ J. Guo, Y. M.  Liu, Y. Ye and Z. Zhang,} An explicit  construction of simple-minded systems over self-injective Nakayama algebras.  Colloq. Math.  164 (2021), 185-210.}

\bibitem{GLYZ}{{J. Guo, Y. M.  Liu, Y. Ye and Z. Zhang}, On simple-minded systems over representation-finite self-injective algebras. Algebr Represent Theor 25 (2022), 983–1002.}

\bibitem{JKS}{{ A. Jaworska-Pastuszak, M. Kwiecie\'{n} and A. Skowro\'{n}ski,} Self-injective algebras having a generalized standard family of quasi-tubes maximally saturated by simple and projective modules.  J. Algebra  520 (2019), 367-399.}
		
\bibitem{K}{{ H. Krause},  The kernel of an irreducible map. Proc. Amer. Math. Soc. 121 (1994), 57-66.}	

\bibitem{KL}{{ S. Koenig and Y. M. Liu,} Simple-minded systems in stable module categories. Quart. J. Math. 63  (2012), 653-674.}

\bibitem{MS}{{ P. Malicki, and A. skowronski,}  Almost cyclic coherent components of an Auslander-Reiten quiver. J. Algebra 229  (2000), 695-749.}	

\bibitem{OZ}{{S. Opper and A. Zvonareva,} Derived equivalence classification of Brauer graph algebras. Adv. Math. 402 (2022), 108341.}

\bibitem{P}{{Z. Pogorza\l{}y,} Properties that are left invariant under stable equivalence. Comm. in Algebra 18 (1990), 4141-4169.}

\bibitem{P1}{{Z. Pogorza\l{}y,} Algebra stably equivalence to self-injective special biserial algebras. Comm. in Algebra 22 (1994), 509-517.}

\bibitem{Ro}{{K. Roggenkamp,}  Biserial algebras and graphs. (English summary) Algebras and modules, II (Geiranger, 1996),  CMS Conf. Proc., 24, Amer. Math. Soc., Providence, RI, 481-496, 1998.}

\bibitem{S0}{{S. Schroll}, Brauer graph algebras.  { Homological methods, representation theory, and cluster algebras.} CRM Short Courses. Springer,  177-223, 2018.}

\bibitem{S1}{{A. Skowro\'{n}ski,} Generalized standard Auslander-Reiten components. J. Math. Soc. Japan 46 (1994), 517-543.}	 

\bibitem{S2}{{A. Skowro\'{n}ski,}  Self-injective algebras:  finite and tame type. Trends in presentation theory of algebras and related topics.   Contemp. Math., 406, Amer. Math. Soc., Province, RI, 169-238, 2006.} 

\bibitem{SK}{{A. Skowro\'{n}ski,} Group algebras of polynomial growth. Manuscr. Math. 59 (1987), 499-516.}		

\bibitem{SY}{{ A.  Skowro\'{n}ski, and K. Yamagata,} Galois coverings of self-injective algebras by repetitive algebras. Trans. Amer. Math. Soc. 351 (1999), 715-734.}	

\bibitem{Z0}{{ Z. Zhang,}  On simple-minded systems over domestic Brauer graph algebras. Eur. J.  Math. 10 (2024), https://doi.org/10.1007/s40879-023-00717-x.}

\bibitem{Z}{{ Z. Zhang}, Simple-minded systems and weakly simple-minded systems  over self-injective algebras. Preprint.}
\end{thebibliography}
\end{document}